\documentclass{amsart}
\usepackage{amsmath,amssymb, amsbsy}
\usepackage{subfigure}
\usepackage{color,psfrag}
\usepackage{enumerate}
\usepackage[dvips]{graphicx}
\usepackage[a4paper,width={15cm},left=3cm,bottom=2.5cm, top=2.5cm,
includeheadfoot]{geometry}
\usepackage[utf8]{inputenc}

\newcommand{\R}{{\mathbb R}}
\newcommand{\N}{{\mathbb N}}

\newcommand{\SN}{{\mathbb S}^{N-1}}
\newcommand{\weakly}{\rightharpoonup}
\newcommand{\e }{\varepsilon}

\newcommand{\alchi}{\raisebox{1.7pt}{$\chi$}}

\renewcommand{\geq }{\geqslant}

\renewcommand{\leq }{\leqslant}

\newenvironment{pf}{\noindent{\sc Proof}.\enspace}{\hfill\qed\medskip}
\newenvironment{pfn}[1]{\noindent{\bf Proof of
    {#1}.\enspace}}{\hfill\qed\medskip}
\newtheorem{Theorem}{Theorem}[section]
\newtheorem{Corollary}[Theorem]{Corollary}
\newtheorem{Lemma}[Theorem]{Lemma}
\newtheorem{Proposition}[Theorem]{Proposition}
\theoremstyle{definition}

\newtheorem{remark}[Theorem]{Remark}

\begin{document}
\title[Singularity of eigenfunctions]{Singularity of eigenfunctions at
  the junction of shrinking tubes. Part I.}

\author[V. Felli]{Veronica Felli}
\address{\hbox{\parbox{5.7in}{\medskip\noindent{Universit\`a di Milano
        Bicocca,\\
        Dipartimento di Ma\-t\-ema\-ti\-ca e Applicazioni, \\
        Via Cozzi
        53, 20125 Milano, Italy. \\[3pt]
        \em{E-mail address: \tt veronica.felli@unimib.it}.}}}}
        
        \author[S. Terracini]{Susanna Terracini}
\address{\hbox{\parbox{5.7in}{\medskip\noindent{Universit\`a di Torino,\\
        Dipartimento di Ma\-t\-ema\-ti\-ca ``G. Peano'', \\
        Via Carlo Alberto 10, 10123 Torino, Italy. \\[3pt]
        \em{E-mail address: \tt susanna.terracini@unito.it}.}}}}

\date{March 30, 2013}

\thanks{
\indent Supported by the PRIN2009 grant ''Critical Point Theory and 
Perturbative Methods for Nonlinear Differential\\
\indent Equations''.}

 \begin{abstract}
   Consider two domains connected by a thin tube: it can be shown
   that, generically, the mass of a given eigenfunction of the
   Dirichlet Laplacian concentrates in only one of them. The
   restriction to the other domain, when suitably normalized, develops
   a singularity at the junction of the tube, as the channel section
   tends to zero.  Our main result states that, under a nondegeneracy
   condition, the normalized limiting profile has a singularity of
   order $N-1$, where $N$ is the space dimension. We give a precise
   description of the asymptotic behavior of eigenfunctions at the
   singular junction, which provides us with some important
   information about its sign near the tunnel entrance.
   More precisely,
   the solution is shown to be one-sign in a neighborhood of the
   singular junction. In
   other words, we prove that the nodal set does not enter inside the
   channel.
 \end{abstract}

\maketitle

\section{Introduction and statement of the main
  results}\label{sec:intr-stat-main}

We are concerned with the behavior of eigenfunctions of the Dirichlet
Laplacian on dumbbell domains depending on a parameter and
disconnecting in some limit process. More precisely, let us consider
two slightly different domains which are connected by a thin tube so
that the mass of a given eigenfunction is concentrated in one of the
two domains.  Then the restriction of the eigenfunction to the other
domain develops a singularity right at the junction of the tube, as
the section of the channel shrinks to zero. The purpose of this paper
is to describe the features of this singularity formation.

A strong motivation for the interest in the spectral analysis of thin branching
domains comes from the theory of quantum graphs modeling waves in thin
graph-like structures (narrow waveguides, quantum wires, photonic
crystals, blood vessels, lungs) and having applications in
nanotechnology, optics, chemistry, medicine, see e.g. \cite{kuchment,
  CF} and references therein.

The behavior of the eigenvalues and eigenfunctions of the Laplace
operator in varying domains has been intensively studied in the
literature starting from \cite{BV,CH,grigorieff,RT,stummel} and more recently in
\cite{AD,arrieta, AK,BZ, DD,daners}, where spectral continuity is discussed
under different kind of perturbations and boundary conditions (of
either Dirichlet or Neumann type).  The problem of rate of convergence
for eigenvalues of elliptic systems was investigated in \cite{taylor},
while in \cite{BHM} estimates of the splitting between the first two
eigenvalues of elliptic operators under Dirichlet boundary conditions
are provided. We also mention that some results on the behavior of
eigenfunctions of the Laplace operator under singular perturbation
adding a thin handle to a compact manifold have been obtained in
\cite{anne}.  As far as the nonlinear counterpart of the problem is concerned, the
effect of the domain shape on the number of positive solutions to some
nonlinear Dirichlet boundary value problems has been investigated in
\cite{Dancer1,Dancer2}, where domains constructed as connected approximations to a finite number
of disjoint or touching balls have been considered, proving that the number of
positive solutions which are not ``large'' grows with the number of
the balls.

When dealing with a dumbbell domain which is
going to disconnect, the spectral continuity proved
e.g. in \cite{daners} implies that eigenfunctions of the approximating
problem converge to the eigenfunction of some limit eigenvalue problem
on a domain with two connected components, whose spectrum is
therefore the union of the spectra on the two components; as a
consequence, if an eigenfunction of the limit problem is supported in
one of the two domains, then
 the corresponding eigenfunction of the approximating problem is going
 to vanish on the other domain. We are going to show that a suitable
 normalization of 
 such eigenfunction  develops a singularity at
 the junction of the tube, whose rate is related to the order of the zero that the limit
eigenfunction has at the other junction  (see Theorem \ref{t:main}). The description
of the behavior of eigenfunctions at the junction will also provide
informations about nodal sets; more precisely we
will prove in Corollary \ref{c:cor} that if the limit
eigenfunction has at one junction of the tube a zero of order one,
then the nodal regions of the corresponding eigenfunctions on the
dumbbell stay away from the other junction. 

In this paper we set up a strategy to evaluate the rate to the
singularity at the junction, based upon a sharp control of the
transversal frequencies along the connecting tube.  To this aim, we
shall exploit the monotonicity method introduced by Almgren
\cite{almgren} in 1979 and then extended by Garofalo and Lin \cite{GL}
to elliptic operators with variable coefficients in order to prove
unique continuation properties. We mention that monotonicity methods
were recently used in \cite{FFT,FFT2,FFT3} to prove not only unique
continuation but also precise asymptotics near singularities of
solutions to linear and semilinear elliptic equations with singular
potentials, by extracting such precious information from the behavior
of the quotient associated with the Lagrangian energy.

As a paradigmatic example, let us consider the following dumbbell domain
in $\R^N=\R\times \R^{N-1}$, $N\geq 3$,
$$
\Omega^\e=D^-\cup \mathcal C_\e\cup D^+
$$
where    $\e\in (0,1)$,
\begin{align*}
D^-&=\{(x_1,x')\in \R\times \R^{N-1}:x_1<0\},\\
\mathcal C_\e&=\big\{(x_1,x')\in \R\times \R^{N-1}:0\leq x_1\leq1,\ \tfrac{x'}{\e}\in\Sigma\big\},\\
D^+&=\{(x_1,x')\in \R\times \R^{N-1}:x_1>1\},
\end{align*}
and $\Sigma\subset \R^{N-1}$ is an open bounded  set with
$C^{2,\alpha}$-boundary containing $0$. For simplicity of notation,
without loss of generality, we
assume that $\Sigma$ satisfies 
\begin{equation}\label{eq:condsigma}
\big\{x'\in\R^{N-1}:|x'|\leq\tfrac1{\sqrt{2}}\big\}\subset \Sigma\subset
\{x'\in\R^{N-1}:|x'|<1\}.
\end{equation}

\begin{figure}[h]
 \centering
   \begin{psfrags}
     \psfrag{D-}{$D^-$}
     \psfrag{D+}{$D^+$}
\psfrag{e}{${\scriptsize{\e}}$}
\psfrag{1}{${\scriptsize{1}}$}
\psfrag{C}{$\mathcal C_\e$}
     \includegraphics[width=7cm]{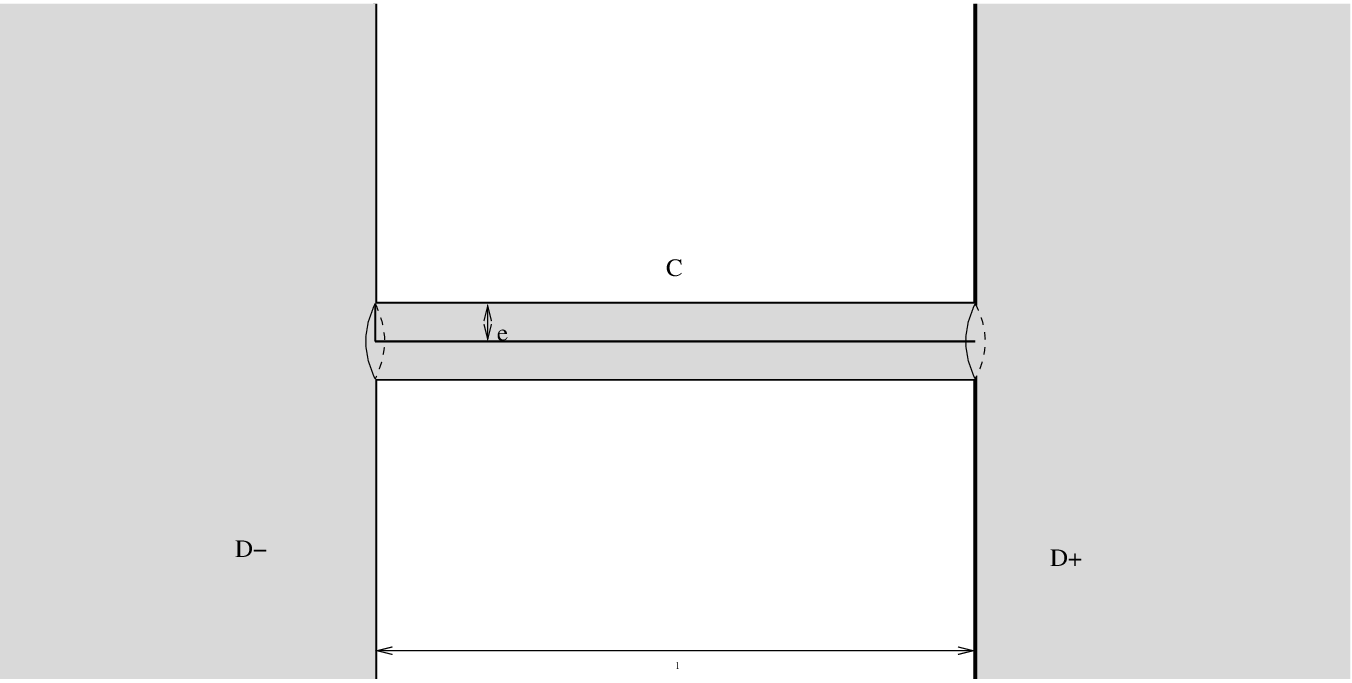}
   \end{psfrags}
 \caption{The domain $\Omega^\e$.}\label{fig:dd}
\end{figure}

\noindent We also denote, for all $t>0$,
$$
B^+_t:=D^+\cap B({\mathbf e}_1,t),\quad
B^-_t:=D^-\cap B({\mathbf 0},t),
$$
where ${\mathbf e}_1
=(1,0,\dots,0)\in \R^N$ and $B(P,t):=\{x\in\R^N:|x-P|<t\}$ denotes the
ball of radius $t$ centered at $P$.
Let $p\in C^1(\R^N,\R)\cap L^{\infty}(\R^N)$ satisfying 
\begin{align}
\label{eq:p} & p\geq 0\text{ a.e. in }\R^N,\ 
p\in L^{N/2}(\R^N),\ \nabla p(x)\cdot x\in L^{N/2}(\R^N),
\ \frac{\partial p}{\partial x_1}\in L^{N/2}(\R^N),\\
\label{eq:p2} &  
\begin{cases}
p\not\equiv 0\text{ in }D^-,\quad
 p\not\equiv 0\text{ in }D^+,\\
 p(x)=0\text { for all } 
x\in \{(x_1,x')\in \R\times \R^{N-1}:1/2\leq x_1\leq1,\ x'\in\Sigma\}
\cup B^+_3.
\end{cases}
\end{align}
While assumption \eqref{eq:p} makes the problem consistent with the usual spectral theory,   \eqref{eq:p2}  is introduced for technical reasons; we don't believe it is necessary:  its only
use is in section \ref{sec:estimates-u_e-right},
to prove some uniform estimates for approximating eigenfunctions close to the
right junction uniformly with respect to the parameter $\e$.
 Possible weakening of assumption (\ref{eq:p2}) is the 
object of a current elaboration.

By classical spectral theory, for every open set $\Omega\subset\R^N$
such that $p\not\equiv 0$ in $\Omega$, 
 the weighted eigenvalue problem 
$$
\begin{cases}
-\Delta \varphi=\lambda p \varphi,&\text{in }\Omega,\\
\varphi=0,&\text{on }\partial \Omega,
\end{cases}
$$
admits a sequence of diverging eigenvalues $\{\lambda_k(\Omega)\}_{k\geq
  1}$; in the enumeration 
$$
\lambda_1(\Omega)\leq
\lambda_2(\Omega)\leq\cdots\leq\lambda_k(\Omega)\leq \cdots
$$
 we repeat each
eigenvalue as many times as its multiplicity. We
denote $\sigma_p(\Omega)=\{\lambda_k(\Omega):k\geq 1\}$. 
For all $\e\in(0,1)$, we also denote 
$$
\lambda^\e_k=\lambda_k(\Omega^\e),
\quad \sigma^\e_p=\sigma_p(\Omega^\e).
$$
It is easy to verify that $\sigma_p(D^-\cup D^+)=\sigma_p(D^-)\cup \sigma_p(D^+)$. 
Let us assume that there exists $k_0\geq 1$ such that 
\begin{align}
  \label{eq:53} &\lambda_{k_0}(D^+)\text{ is simple and the
    corresponding eigenfunctions have in ${\mathbf e}_1
    $ a zero of order $1$},\\
  \label{eq:54} &\lambda_{k_0}(D^+)\not\in \sigma_p(D^-).
  \end{align}
In view of \cite{Micheletti}, these non degeneracy assumptions hold generically with respect to domain (and weight) variations. We can then fix an eigenfunction $\varphi_{k_0}^+\in{\mathcal
    D}^{1,2}(D^+)\setminus\{0\}$ associated to $\lambda_{k_0}(D^+)$,
  i.e. solving
$$
\begin{cases}
-\Delta \varphi_{k_0}^+=\lambda_{k_0}(D^+) p \varphi_{k_0}^+,&\text{in }D^+,\\
\varphi_{k_0}^+=0,&\text{on }\partial D^+,
\end{cases}
$$
such that 
\begin{equation}\label{eq:13}
\frac{\partial \varphi_{k_0}^+}{\partial x_1}({\mathbf e}_1
)>0.
\end{equation}
Here and in the sequel, for every open set $\Omega\subseteq\R^N$, ${\mathcal
    D}^{1,2}(\Omega)$ denotes the functional space obtained as completion of 
$C^\infty_{\rm c}(\Omega)$ 
with respect to the Dirichlet norm $\big(\int_{\Omega}|\nabla u|^2dx\big)^{1/2}$.

We refer to \cite[Example 8.2, Corollary 4.7, Remark 4.3]{daners} for
the proof of the following lemma.
\begin{Lemma}\label{l:conv_eigen}
Let 
\begin{align*}
\bar k&=\mathop{\rm card}\big\{
j\in\N\setminus\{0\}: \lambda_j(D^-\cup D^+)\leq \lambda_{k_0}(D^+)\}\\
&=k_0+
\mathop{\rm card}\big\{
j\in\N\setminus\{0\}:\lambda_j(D^-)\leq \lambda_{k_0}(D^+)\},
\end{align*}
so that $\lambda_{k_0}(D^+)=\lambda_{\bar k}(D^-\cup D^+)$. Then
\begin{equation}\label{eq:52}
\lambda^\e_{\bar k}\to \lambda_{k_0}(D^+)\quad\text{as }\e\to0^+.
\end{equation}
Furthermore, for every $\e$ sufficiently small, $\lambda^\e_{\bar k}$
is simple and there exists an eigenfunction $\varphi_{\bar k}^\e$
associated to $\lambda^\e_{\bar k}$, i.e. satisfying
$$
\begin{cases}
-\Delta \varphi_{\bar k}^\e=\lambda^\e_{\bar k} p \varphi_{\bar k}^\e,&\text{in }\Omega^\e,\\
\varphi_{\bar k}^\e=0,&\text{on }\partial \Omega^\e,
\end{cases}
$$
such that
$$
\varphi_{\bar k}^\e\to \varphi_{k_0}^+\quad\text{in }{\mathcal D}^{1,2}(\R^N)
\quad\text{as }\e\to0^+,
$$
where in the above formula we mean the functions $\varphi_{\bar
  k}^\e,\varphi_{k_0}^+$ to be trivially extended to the whole $\R^N$.
\end{Lemma}

\noindent 
We mention that uniform convergence of eigenfunctions has been
established in \cite[\S 5.2]{bucur2006}.

Henceforward, for simplicity of notation, we denote
\begin{equation}
  \label{eq:193}
u_\e=\varphi_{\bar k}^\e
\quad\text{and}\quad 
u_0=\varphi_{k_0}^+.  
\end{equation}
Hence, for small $\e$, $u_\e$ solves
\begin{align}\label{eq:24}
\begin{cases}
-\Delta u_\e=\lambda^\e_{\bar k} p u_\e,&\text{in }\Omega^\e,\\
u_\e=0,&\text{on }\partial \Omega^\e.
\end{cases}
\end{align}
The main result of the present paper is the following theorem
describing the behavior as $\e\to 0^+$ of $u_\e$  at the junction
${\mathbf 0}=(0,\dots,0)$.  
For all $t>0$, let us denote 
$$
\mathcal D_t^-:=\{v\in C^\infty(D^-\setminus B_t^-):\mathop{\rm supp}v\Subset 
D^-\}
$$
and let $\mathcal H_t^-$ be the completion of $\mathcal D_t^-$ with
respect to the norm $\big(\int_{D^-\setminus {B_t^-}}|\nabla
v|^2dx\big)^{1/2}$, i.e.  $\mathcal H_t^-$ is the space of functions
with finite energy in $D^-\setminus \overline{B_t^-}$ vanishing on
$\partial D^-$. 
We also define, for all $t>0$, 
\begin{equation}\label{eq:defGamma_r-}
\Gamma_t^-=D^-\cap \partial  B^-_{t}.
\end{equation}
Let
\begin{align}\label{eq:Y1}
  Y_1: {\mathbb S}^{N-1}_-\to\R, \quad Y_1 (\theta_1,\theta_2,\dots,\theta_N)=-\frac{\theta_1}{\Upsilon_N},
\end{align}
where 
\begin{align}\label{eq:upsilonN}
{\mathbb
  S}^{N-1}_-:=\{\theta=(\theta_1,\theta_2,\dots,\theta_N)\in{\mathbb
  S}^{N-1}:\theta_1<0\},\quad
\Upsilon_N=\sqrt{{\textstyle{\int}}_{{\mathbb
        S}^{N-1}_-}\theta_1^2d\sigma(\theta)},
\end{align}
being ${\mathbb S}^{N-1}$ the unit $(N-1)$-dimensional sphere.  Here
and in the sequel, the notation $d\sigma$ is used to denote the volume
element on $(N-1)$-dimensional surfaces.  We notice that $Y_1$ is the
first positive $L^2({\mathbb S}^{N-1}_-)$-normalized eigenfunction of
$-\Delta_{{\mathbb S}^{N-1}}$ on ${\mathbb S}^{N-1}_-$ under null
Dirichlet boundary conditions and satisfies $-\Delta_{{\mathbb
    S}^{N-1}}Y_1=(N-1) Y_1$ on ${\mathbb S}^{N-1}_-$, where
$\Delta_{\SN}$ is the Laplace-Beltrami operator on the unit sphere
$\SN$.

\begin{Theorem}\label{t:main}
  Let us assume
  \eqref{eq:p}--\eqref{eq:13} hold
and let $u_\e$ as in \eqref{eq:193}. Then there exists $\tilde
h\in(0,1)$ such that, 
for every sequence $\e_{n}\to 0^+$, there exist a subsequence
  $\{\e_{n_j}\}_j$, $U\in
  C^2(D^-)\cup\big(\bigcup_{t>0}\mathcal H_t^-\big)$, $U\not\equiv 0$,
 and $\beta<0$ such
  that
\begin{align*}
i)\quad&
\frac{u_{\e_{n_j}}}{\sqrt{\int_{\Gamma^-_{\tilde
        h}}u_{\e_{n_j}}^2d\sigma}}\to U  \quad\text{as }j\to+\infty
\quad\ \ \text{strongly in }\mathcal H_t^- \text{ for every }t>0\\[-10pt]
& \hskip5.7cm\text{and in }
C^2(\overline{B_{t_2}^-\setminus B_{t_1}^-}) \text{ for all }0<t_1<t_2;\\[5pt]
ii)\quad&
\lambda^{N-1} U(\lambda x)\to \beta\,\frac{x_1}{|x|^N} \quad\text{as }\lambda \to 0^+
\quad\text{strongly in }\mathcal H_t^- \text{ for every }t>0\\[-5pt]
& \hskip5.7cm\text{and in }
C^2(\overline{B_{t_2}^-\setminus B_{t_1}^-}) \text{ for all }0<t_1<t_2;\\[5pt]
iii) \quad&\beta=-\frac{
\int_{{\mathbb
  S}^{N-1}_-}UY_1\,d\sigma
-\tfrac{\lambda_{k_0}(D^+)}N\int_{D^-}p(x)U(x)Y_1(\tfrac{x}{|x|})
\Big(|x|\alchi_{B^-_1}(x)+\frac{\chi_{D^-\setminus B^-_1}(x)}{|x|^{N-1}}\Big)
dx}{\Upsilon_N}.
\end{align*}
%\end{itemize}
 \end{Theorem}
In the forthcoming paper \cite{aft}, some improvements of Theorem
\ref{t:main} will be obtained; more precisely, the dependence on the
subsequence will be removed 
and the exact asymptotic behavior of the normalization $\sqrt{\int_{\Gamma^-_{\tilde
        h}}u_{\e_{n_j}}^2d\sigma}$ will be derived.

 The description of the behavior of eigenfunctions at the junction
 given by  Theorem
\ref{t:main} provides us with some important information
 about the sign of $u_\e$ near the left junction. More precisely, the
 nondegeneracy condition \eqref{eq:53} on the right junction implies
 that the solution is one-sign in a neighborhood of the left one. In
 other words, the nodal set of $u_\e$ does not enter inside the
 channel.

\begin{Corollary}\label{c:cor}
  Let us assume
  \eqref{eq:p}--\eqref{eq:13} hold
and let $u_\e$ as in \eqref{eq:193}. Then there
    exists $R>0$ such that 
\begin{align*}
  \text{for every }r\in(0,R)\text{ there exists }\e_r>0\text{ such
    that }u_\e>0\text{ in } \Gamma_r^-\text{
    for all }\e\in(0,\e_r).
\end{align*}
\end{Corollary}

The paper is organized as follows. In section
\ref{sec:estimates-u_e-right} we prove some estimates from above and
from below of eigenfunctions of the approximating problem close to the
right junction uniformly with respect to the parameter $\e$.  In
section \ref{sec:freq-funct} we introduce a frequency function
associated to the approximating problem and study its behavior at the
left, in the corridor, and at the right of the domain.  Sections
\ref{sec:limiting-problem} and \ref{sec:blow-up-at-left} contain a
blow-up analysis (at the right and at the left junction respectively)
leading to some uniform bounds of the frequency function which allow
describing, in section \ref{sec:asymptotics-at-left}, the asymptotic
behavior of the eigenfunctions (suitably normalized) close to the left
junction of the tube, thus proving Theorem \ref{t:main} and Corollary
\ref{c:cor}.

\section{Estimates on $u_\e$ on the right}\label{sec:estimates-u_e-right}

This section collects some estimates of eigenfunctions $u_\e$ close to
the right junction, which will be crucial to control the frequency
function at the right.
\begin{Lemma}\label{l:U_e_ud}
There exist $0<r_0<3$, $\e_0\in(0,r_0/2)$, and $C_0>0$ such that 
$$
\frac1{C_0}(x_1-1)\leq u_\e(x)\leq C_0(x_1-1)
\quad \text{for all }x\in D^+\cap\partial B_{r_0}^+\text{ and }\e\in(0,\e_0).
$$
\end{Lemma}
\begin{pf}
From Lemma \ref{l:conv_eigen} and classical elliptic regularity theory,
\begin{equation}\label{eq:14}
u_\e\to u_0 \text{ in }C^2_{\rm loc}(\overline{D^+}\setminus\{{\mathbf e}_1
\})
\text{ and } 
\nabla u_\e\to \nabla u_0 \text{ in }C^1_{\rm loc}(\overline{D^+}\setminus\{{\mathbf e}_1
\}).
\end{equation} 
Furthermore (\ref{eq:13}) implies that there exist $C>0$ and
$r_0\in(0,3)$ such that
\begin{equation}\label{eq:15}
\frac{\partial u_0}{\partial x_1}(x)\geq C,
\quad u_0(x)>0,\quad\text{for all }x\in B_{r_0}^+.
\end{equation}
Let $t_0\in(1,1+r_0/4)$ such that,  if $x=(x_1,x')\in \mathcal A_0:=\big(B_{r_0}^+\setminus
B_{(3r_0)/4}^+\big)\cap \{1<x_1<t_0\}$, then $(1,x')\in B_{r_0}^+\setminus
B_{r_0/2}^+$. By (\ref{eq:15}) and continuity of $u_0$, there exist $c>0$ such that 
\begin{align}\label{eq:18}
u_0(x)\geq c\quad \text{for all }x\in \big(B_{r_0}^+\setminus
B_{(3r_0)/4}^+\big)\setminus\mathcal A_0.
\end{align}
From (\ref{eq:14}), there exists $\e_0\in(0,r_0/2)$ such that 
equation \eqref{eq:24} is satisfied for $\e\in(0,\e_0)$ and
\begin{align}\label{eq:16}&\bigg|\frac{\partial u_\e}{\partial x_1}(x)-\frac{\partial
  u_0}{\partial x_1}(x)\bigg|\leq \frac C2
\quad\text{for all }x\in B_{r_0}^+\setminus B_{r_0/2}^+\text{ and }\e\in(0,\e_0),\\
&\label{eq:17}|u_\e(x)-u_0(x)|\leq \frac c2
\quad\text{for all }x\in \big(B_{r_0}^+\setminus
B_{(3r_0)/4}^+\big)\setminus\mathcal A_0\text{ and }\e\in(0,\e_0).
\end{align}
Estimate (\ref{eq:17}) together with (\ref{eq:18}) implies that 
\begin{align}\label{eq:20}
u_\e(x)\geq \frac c2\quad \text{for all }x\in \big(B_{r_0}^+\setminus
B_{(3r_0)/4}^+\big)\setminus\mathcal A_0\text{ and }\e\in(0,\e_0).
\end{align}
On the other hand, (\ref{eq:16}) together with (\ref{eq:15}) implies 
\begin{align}\label{eq:19}
\frac{\partial u_\e}{\partial x_1}(x)\geq \frac C2
\quad\text{for all }x\in B_{r_0}^+\setminus B_{r_0/2}^+\text{ and }\e\in(0,\e_0).
\end{align}
We notice that, if $x\in\mathcal A_0$ then from (\ref{eq:19}) it follows that  
\begin{align}\label{eq:21}
u_\e(x_1,x')=u_\e(1,x')+\int_1^{x_1}
\frac{\partial u_\e}{\partial x_1}(s,x')\,ds>0.
\end{align}
Combining (\ref{eq:20}) and (\ref{eq:21}) we conclude that
\begin{align}\label{eq:22}
u_\e(x)>0 \quad \text{for all }x\in B_{r_0}^+\setminus
B_{(3r_0)/4}^+\text{ and }\e\in(0,\e_0).
\end{align}
If $x\in D^+\cap\partial B_{r_0}^+$ and $\e\in(0,\e_0)$, from 
(\ref{eq:19}) and (\ref{eq:22}) we have that 
\begin{align*}
u_\e(x)=u_\e\Big(x-\frac{x_1-1}{4}{\mathbf e}_1
\Big)+\int_0^1\frac{\partial u_\e}{\partial x_1}
\Big(x-\frac{(1-t)(x_1-1)}{4}{\mathbf e}_1
\Big)\frac{x_1-1}{4}\,dt
\geq \frac C2\,\frac{x_1-1}{4}
\end{align*}
thus proving the stated lower bound.
The upper bound follows combining (\ref{eq:16}), (\ref{eq:17}), and (\ref{eq:21}).\end{pf}

\noindent The following  iterative  Brezis-Kato type argument yields a uniform $L^\infty$-bound for
$\{u_\e\}_\e$.
\begin{Lemma}\label{l:bk}
There exists $C_1>0$ such that
$$
|u_\e(x)|\leq C_1\quad\text{for all }x\in\Omega^\e\text{ and }\e\in(0,\e_0).
$$
\end{Lemma}
\begin{pf}
Since $u_\e\to u_0$ in ${\mathcal D}^{1,2}(\R^N)$, we have that 
\begin{equation}\label{eq:25}
\sup_{\e\in(0,\e_0)}\|u_\e\|_{L^{2^*}(\R^N)}<\infty.
\end{equation}
We claim that 
\begin{gather}\label{eq:23}
\begin{array}{ll}
  &\text{there exists a positive constant $C>0$ independent of $\e$ and $q$ 
  }\\
  & \text{such that if }u_\e\in L^q(\R^N)
  \text{ for some $q\geq 2^*$ and all $\e\in(0,\e_0)$ then }
\end{array}\\
\notag \|u_\e\|_{L^{\frac{q2^*}2}(\R^N)}\leq C^{\frac 1q}(q-2)^{\frac 1q} 
\|u_\e\|_{L^{q}(\R^N)}.
\end{gather}
The claim can be proved by following the Brezis-Kato procedure
\cite{BrezisKato}. For every $n\in\N$, we set $u_\e^n=\min\{n,|u_\e|\}$ and test 
(\ref{eq:24}) with $u_\e(u_\e^n)^{q-2}$ thus obtaining
$$
(q-2)\int_{\Omega^\e}|\nabla u_\e^n|^2(u_\e^n)^{q-2}\,dx+
\int_{\Omega^\e}|\nabla u_\e|^2(u_\e^n)^{q-2}dx=
\lambda_\e\int_{\Omega^\e}p u_\e^2(u_\e^n)^{q-2}dx.
$$
Letting $C(q)=\min\big\{\frac2{q-2},\frac12\big\}$, we then obtain
\begin{align*}
  C(q)\int_{\Omega^\e}|\nabla ((u_\e^n)^{\frac q2-1}u_\e)|^2dx&\leq
  C(q)\int_{\Omega^\e}\bigg(\frac{(q-2)^2}{2}(u_\e^n)^{q-2}|\nabla
  u_\e^n|^2+2(u_\e^n)^{q-2}|\nabla u_\e|^2\bigg)dx\\
&\leq \lambda_\e\int_{\Omega^\e}p u_\e^2(u_\e^n)^{q-2}dx\leq {\rm const\,}
\int_{\Omega^\e}|u_\e|^qdx
\end{align*}
for some ${\rm const\,}>0$ independent of $\e$ and $q$, which, letting
$n\to+\infty$, implies claim (\ref{eq:23}) by Sobolev inequality.
Starting from $q=2^*$ and iterating the estimate of claim
(\ref{eq:23}), we obtain that, for all $n\in\N$, $n\geq1$, letting
$q_n=2\big(\frac{2^*}{2}\big)^n$, there holds
\begin{align*}
\|u_\e\|_{L^{q_{n+1}}(\R^N)}\leq \|u_\e\|_{L^{2^*}(\R^N)}C^{\sum_{k=1}^n\frac 1{q_k}}
\prod_{k=1}^n(q_k-2)^{\frac 1{q_k}}\leq  
{\rm const\,}\|u_\e\|_{L^{2^*}(\R^N)}
\end{align*}
for some ${\rm const\,}>0$ independent of $\e$ and $n$. Letting
$n\to\infty$, (\ref{eq:25}) yields the conclusion.~\end{pf}

\noindent We denote
\begin{align}
\label{eq:94}  &T_1^-=\{(x_1,x'):x'\in\Sigma,x_1\leq1\},\quad \widetilde D=D^+\cup T_1^-,\\
\notag &T_\e^-=\big\{(x_1,x'):\tfrac{x'}{\e}\in\Sigma,\,x_1\leq1\big\},\\
\notag&T_1=\{(x_1,x'):x_1\in\R,\,x'\in\Sigma\},
\end{align}
and, for $r\in \R\setminus (1,2)$,
\begin{align}\label{eq:115}
\widetilde \Omega_r=
\begin{cases}
\{(x_1,x')\in T_1:x_1<r\},&\text{if }r\leq1,\\
T_1^-\cup B^+_{r-1},&\text{if } r\geq 2,
\end{cases}\qquad
\widetilde \Gamma_r=
\begin{cases}
\{(x_1,x')\in T_1:x_1=r\},&\text{if }r\leq1,\\
\Gamma^+_{r-1},&\text{if } r\geq 2,
\end{cases}
\end{align}
where, for all $t>0$, we denote 
\begin{equation}\label{eq:gamma+}
\Gamma_t^+=D^+\cap \partial  B^+_{t}.
\end{equation}
Let us define
\begin{equation}\label{eq:29}
f:T_1\to\R,\quad 
f(x_1,x')=e^{-\sqrt{\lambda_1(\Sigma)}(x_1-1)}\psi_1^\Sigma(x'),
\end{equation}
where $\lambda_1(\Sigma)$ is the first eigenvalue of the Laplace
operator on $\Sigma$ under null Dirichlet boundary conditions and
$\psi_1^\Sigma(x')$ is the corresponding positive
$L^2(\Sigma)$-normalized eigenfunction, so that 
$$
\begin{cases}
-\Delta_{x'}\psi_1^\Sigma(x')=
\lambda_1(\Sigma)\psi_1^\Sigma(x'),&\text{in }\Sigma,\\
\psi_1^\Sigma=0,&\text{on }\partial\Sigma,
\end{cases}
$$
being $\Delta_{x'}=\sum_{j=2}^N\frac{\partial^2}{\partial x_j^2}$,
$x'=(x_2,\dots,x_N)$. In particular $f\in  C^2(\overline{T_1})$ 
and satisfies 
$$
\begin{cases}
-\Delta f=0,&\text{in }T_1,\\
f=0,&\text{on }\partial T_1.
\end{cases}
$$
Lemma \ref{l:sup_right} below shows how harmonic functions in $D^+$ can be
extended (up to a finite energy perturbation) to harmonic functions in $\widetilde D$ with finite energy at
$-\infty$.  In order to prove it, the following Poincar\'e type
inequality is needed.

\begin{Lemma}\label{l:poincare}
  There exists a constant $C_{P}=C_P(N)$ depending only on the dimension $N$
  such that for every function $v:D^+\setminus B_1^+\to \R$ satisfying
$$
v\in \bigcap_{R>1}H^1(B_R^+\setminus B_1^+)\quad\text{and}
\quad v=0\text{ on }\{x_1=1,|x'|>1\},
$$
there holds 
$$
\int_{B_{2R}^+\setminus B^+_{R}}v^2(x)\,dx\leq C_PR^2
\int_{B_{2R}^+\setminus B^+_{R}}|\nabla v(x)|^2\,dx
\quad\text{for all }R>1.
$$
\end{Lemma}
\begin{pf}
  It follows by scaling of the Poincar\'e inequality for functions
  vanishing on a portion of the boundary.
\end{pf}

\begin{Lemma}\label{l:sup_right}
For every  $\psi\in C^2(D^+)\cap C^1(\overline{D^+})$ such that
$$
\begin{cases}
-\Delta \psi=0,&\text{in }D^+,\\
\psi=0,&\text{on }\partial D^+,
\end{cases}
$$
there exists a unique function $u=\mathcal T (\psi)$ such that 
\begin{align}
\label{eq:2}&\int_{\widetilde \Omega_R}\Big(|\nabla u(x)|^2
+|u(x)|^{2^*}\Big)
\,dx<+\infty\text{ for all }R>2,\\
\label{eq:1}&-\Delta u=0\text{ in a distributional sense in }\widetilde D,
\quad u=0\text{ on }\partial \widetilde D,\\
\label{eq:3}&\int_{D^+}|\nabla (u-\psi)(x)|^2\,dx<+\infty.
\end{align}
Furthermore
\begin{equation}\label{eq:38}
\mathcal T (\psi)-\widetilde \psi\in {\mathcal D}^{1,2}(\widetilde D),
\quad 
\text{where}\quad 
\widetilde \psi:=
\begin{cases}
\psi&\text{in }D^+\\
0&\text{in }T_1^-.
\end{cases}
\end{equation}
\end{Lemma}
\begin{pf}
Let us define $J_\psi:{\mathcal D}^{1,2}(\widetilde D)\to\R$ as 
\begin{align}\label{eq:10}
  J_\psi(\varphi)= \frac12\int_{\widetilde
    D}|\nabla\varphi(x)|^2\,dx-\int_{\Sigma}
 \varphi(1,x')\Big(\frac{\partial\psi}{\partial x_1}\Big)_+(1,x')\,dx'
\end{align}
where $\big(\frac{\partial\psi}{\partial
  x_1}\big)_+(1,x'):=\lim_{t\to0^+}\frac{\psi(1+t,x')}{t}$.
By standard minimization methods it is easy to prove that there exists 
$w\in {\mathcal D}^{1,2}(\widetilde D)$ such that 
$J_\psi(w)=\min_{{\mathcal D}^{1,2}(\widetilde D)}J_\psi$.
In particular $w$ satisfies 
$$
0=dJ_\psi(w)[\varphi]=\int_{\widetilde D}\nabla w(x)\cdot\nabla \varphi(x)\,dx
-\int_{\Sigma}
 \varphi(1,x')\Big(\frac{\partial\psi}{\partial x_1}\Big)_+(1,x')\,dx'
$$
for all $\varphi \in {\mathcal D}^{1,2}(\widetilde D)$. Hence the function 
$u:\widetilde D\to \R$,
$$
u=\begin{cases}
w+\psi,&\text{in }D^+,\\
w,&\text{in }T_1^-,
\end{cases}
$$
satisfies (\ref{eq:2}), (\ref{eq:3}), and, for every $\varphi\in
C^\infty_{\rm c}(\widetilde D)$, 
\begin{align*}
\int_{\widetilde D}\nabla u(x)\cdot\nabla \varphi(x)\,dx
&=\int_{\widetilde D}\nabla w(x)\cdot\nabla \varphi(x)\,dx
+\int_{D^+}\nabla \psi(x)\cdot\nabla \varphi(x)\,dx\\
&=
\int_{\Sigma}
\varphi(1,x')\Big(\frac{\partial\psi}{\partial x_1}\Big)_+(1,x')\,dx'
-\int_{\Sigma}
\varphi(1,x')\Big(\frac{\partial\psi}{\partial x_1}\Big)_+(1,x')\,dx'
=0
\end{align*}
thus implying (\ref{eq:1}). To prove uniqueness, let us assume that 
 $u_1$ and $u_2$ both satisfy
(\ref{eq:1}--\ref{eq:3}); then the difference $u=u_1-u_2$ solves 
\begin{equation}\label{eq:4}
-\Delta u=0\text{ in a distributional sense in }\widetilde D,
\quad u=0\text{ on }\partial \widetilde D,\\
\end{equation}
and satisfies 
\begin{equation}\label{eq:5}
\int_{D^+}|\nabla u(x)|^2dx=
\int_{D^+}|\nabla (u_1-\psi)(x)
-\nabla (u_2-\psi)(x)|^2dx<+\infty.
\end{equation}
For all $R>2$ let $\eta_R$ be a cut-off function satisfying 
\begin{align*}
&\eta_R\in C^\infty(\widetilde D),\quad \eta_R\equiv1\text{ in
}\widetilde \Omega_R,\quad \eta_R\equiv0\text{ in }D^+\setminus B_{2(R-1)}^+,\quad 
|\nabla \eta_R(x)|\leq \frac2{R-1}\text{ in }\widetilde D.
\end{align*}
Multiplying (\ref{eq:4}) with $\eta_R^2u$ and integrating by parts
over $\widetilde D$ we obtain
\begin{align*}
\int_{\widetilde D}|\nabla u(x)|^2\eta_R^2(x)\,dx&=-2
\int_{\widetilde D}u(x)\eta_R(x)\nabla u(x)\cdot\nabla \eta_R(x)\,dx\\
&\leq \frac12\int_{\widetilde D}|\nabla u(x)|^2\eta_R^2(x)\,dx+2
\int_{\widetilde D}u^2(x)|\nabla \eta_R(x)|^2\,dx
\end{align*}
thus implying, in view of Lemma \ref{l:poincare},
\begin{align*}
  \frac12\int_{\widetilde \Omega_R}|\nabla u(x)|^2\,dx
  &\leq \frac12 \int_{\widetilde D}|\nabla u(x)|^2\eta_R^2(x)\,dx\\
  &\leq 2 \int_{\widetilde D}u^2(x)|\nabla \eta_R(x)|^2\,dx
  \leq \frac{8}{(R-1)^2}\int_{B^+_{2(R-1)}\setminus B^+_{R-1}}u^2(x)\,dx\\
  &\leq 8C_P \int_{B_{2(R-1)}^+\setminus B^+_{R-1}}|\nabla
  u(x)|^2\,dx.
\end{align*}
Letting $R\to+\infty$, from (\ref{eq:5}) we deduce that
$\int_{\widetilde D}|\nabla u|^2dx=0$ and hence $u$ must be constant on
$\widetilde D$.  Since $u$ vanishes on $\partial \widetilde D$, we
deduce that $u\equiv 0$ and then $u_1=u_2$ in  $\widetilde D$ thus proving 
uniqueness.\end{pf}

\noindent Henceforward we denote 
\begin{equation}\label{eq:Phi1}
\Phi_1=\mathcal T(x_1-1).
\end{equation}
Since in the case $\psi(x)=x_1-1$ we have that
$\big(\frac{\partial\psi}{\partial x_1}\big)_+(1,x')=1>0$, the minimum
of the functional $J_{x_1-1}$ defined in \eqref{eq:10} is attained by a
nonnegative  function $w$. Hence we deduce that 
\begin{equation}\label{eq:Phi1geq}
\Phi_1(x_1,x')\geq (x_1-1)^+\quad \text{for all }(x_1,x')\in\widetilde D.
\end{equation}
Hence, from the Strong Maximum Principle we deduce that 
\begin{equation}\label{eq:Phi1_pos}
\Phi_1(x_1,x')>0\quad \text{for all }(x_1,x')\in\widetilde D.
\end{equation}
For all $r\in\R$, let us  denote 
\begin{equation}\label{eq:141}
T_{1,r}=:\{(x_1,x'):x'\in\Sigma,\ x_1\leq r\},
\quad 
\Gamma_r:=\{(x_1,x'):x'\in\Sigma,\ x_1= r\},
\end{equation}
and define  
$\mathcal E_r$  as the completion of 
$C^\infty_{\rm c}(T_{1,r})$ 
with respect to the norm $\big(\int_{T_{1,r}}|\nabla v|^2dx\big)^{1/2}$ (which is actually
equivalent to the norm $\big(\int_{T_{1,r}}|\nabla
v|^2dx+\int_{\Gamma_r}v^2d\sigma \big)^{1/2}$), i.e.
$\mathcal E_r$ is the space of finite energy  functions  in
$T_{1,r}$ vanishing on $\{(x_1,x'):x_1\leq r\text{ and }x'\in\partial\Sigma\}$.

The following Lemma associate an Almgren type  frequency function to harmonic
functions in $\mathcal E_R$ and describe its behavior at $-\infty$.

\begin{Lemma}\label{l:limNphi1gen}
Let  $R\in\R$ and $\phi\in\mathcal E_R\setminus\{0\}$ satisfying  
\begin{align*}
\begin{cases}
-\Delta \phi=0,&\text{ in }T_{1,R},\\
\phi=0,&\text{ on }\{(x_1,x'):x_1\leq R\text{ and }x'\in\partial\Sigma\},
\end{cases}
\end{align*}
in a weak sense, and let  $N_\phi:(-\infty,R)\to\R$ be defined as 
$$
N_\phi(r):=\frac{\int_{T_{1,r}}|\nabla \phi(x)|^2dx}
{\int_{\Gamma_r}\phi^2(x)\,d\sigma}.
$$
Then 
\begin{itemize}
\item[i)]
$N_\phi$ is non decreasing in $(-\infty,R)$;
\item[ii)]
there exists $K_0\in\N$, $K_0\geq1$, such that 
$$
\lim_{r\to-\infty}N_\phi(r)=
\sqrt{\lambda_{K_0}(\Sigma)},
$$
where $\lambda_{K_0}(\Sigma)$ is the $K_0$-th eigenvalue of the Laplace
operator on $\Sigma$ under null Dirichlet boundary conditions;
\item[iii)] if $N_\phi\equiv \gamma$ for some $\gamma\in\R$ then
  $\gamma= \sqrt{\lambda_{K_0}(\Sigma)}$ and $\phi(x_1,x')=e^{ \sqrt{\lambda_{K_0}(\Sigma)}x_1}\psi(x')$ for some 
eigenfunction $\psi$ of  
$-\Delta_{x'}$ in $\Sigma$ associated to the eigenvalue $\lambda_{K_0}(\Sigma)$;
\item[iv)] if $\phi>0$  in $T_{1,R}$, then $K_0=1$.
\end{itemize}
\end{Lemma}
\begin{pf}
It is easy to prove that  $N_\phi\in C^1(-\infty,R)$ and,  for all $r\in(-\infty,R)$, 
$$
N_\phi'(r)= 2\frac{\Big(\int_{\Gamma_r}\big|\frac{\partial
    \phi}{\partial x_1}\big|^2
  d\sigma\Big)\Big(\int_{\Gamma_r}\phi^2\,d\sigma\Big)-\Big(\int_{\Gamma_r}\phi\frac{\partial
    \phi}{\partial
    x_1}\,d\sigma\Big)^2}{\Big(\int_{\Gamma_r}\phi^2\,d\sigma\Big)^2}.
$$
Hence, Schwarz's inequality implies that $N_\phi'(r)\geq 0$ for all
$r<R$. Therefore $N_\phi$ is non-decreasing in $(-\infty,R)$ and statement i) is proved.
By monotonicity, there exists 
\begin{align}\label{eq:136}
\gamma:=\lim_{r\to-\infty}N_\phi(r)\in [0,+\infty).
\end{align}
For every $\lambda>0$ let us define
$$
\phi_\lambda(x_1,x'):=\frac{\phi(x_1-\lambda,x')}{\sqrt{\int_{\Gamma_{R-\lambda}}\phi^2d\sigma}}.
$$
We have that $\phi_\lambda\in\mathcal E_{R+\lambda}$,
\begin{equation}\label{eq:130}
\int_{\Gamma_R}\phi_\lambda^2d\sigma=1,
\end{equation}
and $\phi_\lambda$ weakly solves 
\begin{align}\label{eq:131}
\begin{cases}
-\Delta \phi_\lambda=0,&\text{ in }T_{1,R+\lambda},\\
\phi_\lambda=0,&\text{ on }\{(x_1,x'):x_1\leq R+\lambda\text{ and }x'\in\partial\Sigma\}.
\end{cases}
\end{align}
Moreover, the change of
variable $(x_1,x')=(y_1-\lambda,y')$ yields
\begin{align}\label{eq:135}
N_{\phi}(r-\lambda)=\frac{\int_{T_{1,r}}|\nabla\phi_\lambda(y)|^2dy}{\int_{\Gamma_r}\phi_\lambda^2d\sigma}
\quad\text{for all }r<R+\lambda.
\end{align}
In particular we have that 
$$
N_{\phi}(R-\lambda)=\int_{T_{1,R}}|\nabla\phi_\lambda(y)|^2dy\leq N_\phi\Big(\frac R2\Big)
\quad\text{for every }\lambda\geq\frac R2,
$$
and hence $\{\phi_\lambda\}_{\lambda\geq R/2}$ is bounded in $\mathcal
E_R$. Therefore there exist a sequence $\lambda_n\to +\infty$ and some
$\tilde\phi \in\mathcal E_R$ such that $\phi_{\lambda_n}\weakly
\tilde\phi$ weakly in $\mathcal E_R$ and a.e. in $T_{1,R}$.
 From
compactness of the embedding $\mathcal E_R\hookrightarrow
L^2(\Gamma_R)$ and (\ref{eq:130}) we deduce that
$\int_{\Gamma_R}\tilde\phi^2d\sigma=1$; in particular
$\tilde\phi\not\equiv 0$.  Passing to the weak limit in
\eqref{eq:131} as $\lambda_n\to+\infty$ we have that
\begin{align}\label{eq:132}
\begin{cases}
-\Delta \tilde\phi=0,&\text{ in }T_{1,R},\\
\tilde\phi=0,&\text{ on }\{(x_1,x'):x_1\leq R\text{ and }x'\in\partial\Sigma\}.
\end{cases}
\end{align}
By classical elliptic regularity estimates, we also have that $\phi_{\lambda_n}\to
 \tilde \phi$ in $C^2(T_{1,r_2}\setminus
 T_{1,r_1})$ for all $r_1<r_2<R$. Therefore, multiplying
(\ref{eq:132}) by $\tilde \phi$ and integrating over
$T_{1,r}$ with $r<R$, we obtain
\begin{align}\label{eq:133}
\int_{\Gamma_r}\frac{\partial \phi_{\lambda_n}}{\partial x_1}\phi_{\lambda_n}d\sigma
\to
\int_{\Gamma_r}\frac{\partial \tilde\phi}{\partial x_1}\tilde \phi\,d\sigma
=\int_{T_{1,r}}|\nabla \tilde \phi(x)|^2dx.
\end{align}
On the other hand,  multiplication of
(\ref{eq:131}) by $\phi_{\lambda_n}$ and integration by parts over
$T_{1,r}$ yield
\begin{align}\label{eq:134}
\int_{T_{1,r}}|\nabla  \phi_{\lambda_n}(x)|^2dx=
\int_{\Gamma_r}\frac{\partial  \phi_{\lambda_n}}{\partial x_1} \phi_{\lambda_n}\,d\sigma.
\end{align}
From \eqref{eq:133} and \eqref{eq:134}, we deduce that
$\| \phi_{\lambda_n}\|_{\mathcal E_r}\to\| \tilde \phi\|_{\mathcal
  E_r}$ and then $\phi_{\lambda_n}\to\tilde \phi$ strongly in
$\mathcal E_r$ for every $r<R$. Therefore, for every $r<R$, passing to
the limit as $\lambda_n\to+\infty$ in \eqref{eq:135} and letting
$\gamma$ as in \eqref{eq:136}, we obtain that
\begin{align}\label{eq:138}
  N_{\tilde\phi}(r)= \gamma \quad\text{for all }r<R,
\end{align}
 where 
$$
 N_{\tilde\phi}(r)=\frac{\int_{T_{1,r}}|\nabla\tilde\phi(y)|^2dy}{\int_{\Gamma_r}\tilde\phi^2d\sigma}.
$$
Then 
$$
N_{\tilde\phi}'(r)=
2\frac{\Big(\int_{\Gamma_r}\big|\frac{\partial
    \tilde\phi}{\partial x_1}\big|^2
  d\sigma\Big)\Big(\int_{\Gamma_r}\tilde\phi^2\,d\sigma\Big)-\Big(\int_{\Gamma_r}
\tilde\phi\frac{\partial
    \tilde\phi}{\partial
    x_1}\,d\sigma\Big)^2}{\Big(\int_{\Gamma_r}\tilde\phi^2\,d\sigma\Big)^2}=0
 \quad\text{for all }r<R.
$$
Since equality in the Schwarz's inequality holds only for parallel
vectors, we infer that $\frac{\partial \tilde\phi}{\partial x_1}$ and
$\tilde\phi$ must be parallel as vectors in $L^2(\Gamma_r)$, hence
there exists some function $\eta:(-\infty,R)\to\R$ such that 
$$
\frac{\partial \tilde\phi}{\partial x_1}(x_1,x')=\eta(x_1)\tilde\phi(x_1,x')
\quad\text{for all }x_1\in (-\infty,R)\text{ and }x'\in\Sigma.
$$
Integration with respect to $x_1$ yields
\begin{equation}\label{eq:137}
\tilde\phi(x_1,x')=\varphi(x_1)\psi(x')\quad\text{for all }x_1\in
(-\infty,R)\text{ and }x'\in\Sigma,
\end{equation}
where $\varphi(x_1)=e^{\int_R^{x_1} \eta(s)ds}$, $\psi(x')=\tilde\phi(R,x')$.
From \eqref{eq:132} and \eqref{eq:137}, we derive
$$
\varphi''(x_1)\psi(x')+\varphi(x_1)\Delta_{x'}\psi(x')=0.
$$
Taking $x_1$ fixed, we deduce that $\psi$ is an eigenfunction of  
$-\Delta_{x'}$ in $\Sigma$ under homogeneous Dirichlet boundary conditions. 
If $\lambda_{K_0}(\Sigma)$ is the corresponding
eigenvalue then $\varphi(x_1)$ solves the equation
$$
\varphi''(x_1)-\lambda_{K_0}(\Sigma)\varphi(x_1)=0
$$
and hence $\varphi$ is of the form
$$
\varphi(x_1)=c_1 e^{\sqrt{\lambda_{K_0}(\Sigma)}(x_1-R)}+c_2e^{-\sqrt{\lambda_{K_0}(\Sigma)}(x_1-R)}
\quad\text{for some $c_1,c_2\in\R$.}
$$
Since the function
 $e^{-\sqrt{\lambda_{K_0}(\Sigma)}(x_1-R)}\psi(x')\notin \mathcal E_R$, then $c_2=0$ and  
$\varphi(x_1)=c_1 e^{\sqrt{\lambda_{K_0}(\Sigma)}(x_1-R)}$. Since $\varphi(R)=1$, we obtain that $c_1=1$ and
then
\begin{equation}\label{eq:139}
  \tilde\phi(x_1,x')=e^{\sqrt{\lambda_{K_0}(\Sigma)}(x_1-R)}\psi(x'),
  \quad 
  \text{for all }x_1\in
  (-\infty,R)\text{ and }x'\in\Sigma.
\end{equation}
Substituting \eqref{eq:139} into \eqref{eq:138} we obtain that 
$\gamma=\sqrt{\lambda_{K_0}(\Sigma)}$. Hence statement ii) is proved. We notice that the above argument 
of classification of harmonic functions  $\tilde\phi$ with constant frequency
$N_{\tilde\phi}$  also proves statement iii).

In order to prove iv), let us assume that  $\phi>0$ in $T_{1,R}$. Then $\phi_\lambda>0$ in $T_{1,R+\lambda}$.
Hence a.e.  convergence implies that $\tilde\phi\geq 0$ in
$T_{1,R}$.  From the Strong Maximum Principle we obtain that
$\tilde\phi> 0$ in $T_{1,R}$, which necessarily implies that $\psi>0$ in $\Sigma$.
Then $\psi$ must be the eigenfunction associated to the first
eigenvalue, i.e. $\lambda_{K_0}(\Sigma)=\lambda_{1}(\Sigma)$.~\end{pf}

\noindent The previous lemma allows describing the behavior of the
Almgren type frequency quotient naturally associated to the function $\Phi_1$
introduced in \eqref{eq:Phi1}.  For all $r\in \R\setminus (1,2)$, let $\widetilde{\mathcal
  N}(r)=\widetilde{\mathcal N}_{\Phi_1}(r)$ be the frequency function
associated to $\Phi_1$, i.e.
\begin{equation}\label{eq:freq_Phi1}
  \widetilde{\mathcal N}(r)=
\widetilde{\mathcal N}_{\Phi_1}(r)=\frac{\Lambda_N(r)
    \int_{\widetilde\Omega_r}|\nabla \Phi_1(x)|^2dx}
  {\int_{\widetilde\Gamma_r}\Phi_1(x)\,d\sigma},
\end{equation}
where 
\begin{equation}\label{eq:L_N}
\Lambda_N(r)=
\begin{cases}
  \big(\frac2{\omega_{N-1}}\big)^{\!\frac1{N-1}}
  |\widetilde \Gamma_r|^{\frac1{N-1}}=r-1,&\text{if }r\geq2,\\[3pt]
  \big(\frac {N-1}{\omega_{N-2}}\big)^{\!\frac1{N-1}}|\widetilde
  \Gamma_r|^{\frac1{N-1}}=1,&\text{if }r\leq1,
\end{cases}  
\end{equation}
$|\widetilde \Gamma_r|$ denotes the $(N-1)$-dimensional volume of
$\widetilde \Gamma_r$, and $\omega_{N-1}$
is the volume of the unit sphere ${\mathbb S}^{N-1}$,
i.e. $\omega_{N-1}=\int_{{\mathbb S}^{N-1}}d\sigma(\theta)$.

An immediate consequence of Lemma \ref{l:limNphi1gen} and
\eqref{eq:Phi1_pos} is the following corollary. 
\begin{Corollary}\label{l:limNphi1}
$\lim_{r\to-\infty}\widetilde{\mathcal N}(r)=\sqrt{\lambda_1(\Sigma)}$.
\end{Corollary}

As a left counterpart of Lemma \ref{l:sup_right}, we now construct a
harmonic extension to
$\widetilde D$ of the function $f$ defined in (\ref{eq:29}) (up to a
finite energy perturbation in the tube)
having finite energy at the right.

\begin{Lemma}\label{l:Phi2}
There exists a unique function $\Phi_2:\widetilde D\to\R$ such that 
\begin{align}
\label{eq:6}&\int_{D^+}\Big(|\nabla \Phi_2(x)|^2+|\Phi_2(x)|^{2^*}\Big)\,dx<+\infty,\\
\label{eq:9}&-\Delta \Phi_2=0\text{ in a distributional sense in }\widetilde D,
\quad \Phi_2=0\text{ on }\partial \widetilde D,\\
\label{eq:7}&\int_{T_1}|\nabla (\Phi_2-f)(x)|^2\,dx<+\infty,
\end{align}
where $f$ is defined in (\ref{eq:29}).
Furthermore 
\begin{equation}\label{eq:12}
\Phi_2\geq f\quad\text{in }T_1\quad\text{and}\quad
\Phi_2\geq 0\quad\text{in }\widetilde D.
\end{equation}
\end{Lemma}
\begin{pf}
Let us define $J:{\mathcal D}^{1,2}(\widetilde D)\to\R$ as 
\begin{align*}
  J(\varphi)&= \frac12\int_{\widetilde
    D}|\nabla\varphi(x)|^2\,dx+\int_{(1,+\infty)\times\partial\Sigma}
\varphi \frac{\partial f}{\partial\nu}\,d\sigma \\
&=\frac12\int_{\widetilde
    D}|\nabla\varphi(x)|^2\,dx+\int_{(1,+\infty)\times\partial\Sigma}
\varphi(x_1,x') 
e^{-\sqrt{\lambda_1(\Sigma)}(x_1-1)}\frac{\partial \psi_1^\Sigma}{\partial \nu_{x'}}(x')
\,d\sigma,
\end{align*}
where $\nu$ denotes the normal external unit vector to $\partial T_1$
and $\nu_{x'}$ the normal external unit vector to $\partial \Sigma$. It is easy to prove that
$J(\varphi)\geq c_1\|\varphi\|_{\mathcal D^{1,2}(\widetilde D)}^2-c_2$
for some constants $c_1,c_2>0$ and all $\varphi\in \mathcal
D^{1,2}(\widetilde D)$ and that $J$ is weakly lower
semi-continuous. Hence there exists $w\in {\mathcal
  D}^{1,2}(\widetilde D)$ such that $J(w)=\min_{{\mathcal
    D}^{1,2}(\widetilde D)}J$.  Since, by the Hopf Lemma,
$\frac{\partial \psi_1^\Sigma}{\partial \nu_{x'}}<0$ on
$\partial\Sigma$, we can assume that $w\geq 0$ (otherwise we take
$|w|$ which is still a minimizer). The minimizer $w$ satisfies 
$$
0=dJ(w)[\varphi]=\int_{\widetilde D}\nabla w(x)\cdot\nabla \varphi(x)\,dx
+\int_{(1,+\infty)\times\partial\Sigma}
\varphi \frac{\partial f}{\partial\nu}\,d\sigma 
$$
for all $\varphi \in {\mathcal D}^{1,2}(\widetilde D)$. Hence the function 
$\Phi_2:\widetilde D\to \R$,
$$
\Phi_2=\begin{cases}
w+f,&\text{in }T_1,\\
w,&\text{in }\widetilde D\setminus T_1,
\end{cases}
$$
satisfies (\ref{eq:6}), (\ref{eq:7}), (\ref{eq:12}), and, for every $\varphi\in
C^\infty_{\rm c}(\widetilde D)$, 
\begin{align*}
\int_{\widetilde D}\nabla \Phi_2(x)\cdot\nabla \varphi(x)\,dx
&=\int_{\widetilde D}\nabla w(x)\cdot\nabla \varphi(x)\,dx
+\int_{T_1}\nabla f(x)\cdot\nabla \varphi(x)\,dx\\
&=
-\int_{(1,+\infty)\times\partial\Sigma}
\varphi \frac{\partial f}{\partial\nu}\,d\sigma +\int_{(1,+\infty)\times\partial\Sigma}
\varphi \frac{\partial f}{\partial\nu}\,d\sigma =0
\end{align*}
thus implying (\ref{eq:9}).  To prove uniqueness, let us assume that 
 $u_1$ and $u_2$ both satisfy
(\ref{eq:9}--\ref{eq:7}); then the difference $u=u_1-u_2$ solves 
\begin{equation}\label{eq:8}
-\Delta u=0\text{ in a distributional sense in }\widetilde D,
\quad u=0\text{ on }\partial \widetilde D,\\
\end{equation}
and satisfies 
\begin{equation}\label{eq:11}
\int_{T_1}|\nabla u(x)|^2dx=
\int_{T_1}|\nabla (u_1-f)(x)
-\nabla (u_2-f)(x)|^2dx<+\infty.
\end{equation}
For all $t<1$ let $\eta_t$ be a cut-off function satisfying 
\begin{align*}
&\eta_t\in C^\infty(\widetilde D),\quad \eta_1(x_1,x')=1\text{ if }x_1>t,
\quad \eta_t(x_1,x')=0\text{ if }x_1<t-1,\quad 
|\nabla \eta_t(x)|\leq 2\text{ in }\widetilde D.
\end{align*}
Multiplying (\ref{eq:8}) with $\eta_t^2u$ and integrating by parts
over $\widetilde D$ we obtain
\begin{align*}
\int_{\widetilde D}|\nabla u(x)|^2\eta_t^2(x)\,dx&=-2
\int_{\widetilde D}u(x)\eta_t(x)\nabla u(x)\cdot\nabla \eta_t(x)\,dx\\
&\leq \frac12\int_{\widetilde D}|\nabla u(x)|^2\eta_t^2(x)\,dx+2
\int_{\widetilde D}u^2(x)|\nabla \eta_t(x)|^2\,dx
\end{align*}
thus implying
\begin{align*}
  \frac12\int_{\widetilde D\cap \{x_1>t\}}|\nabla u(x)|^2\,dx
  &\leq \frac12 \int_{\widetilde D}|\nabla u(x)|^2\eta_t^2(x)\,dx\\
  &\leq 2 \int_{\widetilde D}u^2(x)|\nabla \eta_t(x)|^2\,dx
  \leq 8\int_{\widetilde D\cap \{t-1<x_1<t\}}u^2(x)\,dx\\
  &\leq 8{\widetilde C}_P \int_{\widetilde D\cap \{t-1<x_1<t\}}|\nabla
  u(x)|^2\,dx
\end{align*}
where the constant ${\widetilde C}_P>0$ depends only on the dimension and is 
the best constant of the Poincar\'e inequality for functions on $(-1,0)\times\Sigma$ 
  vanishing on $\partial\Sigma$.
Letting $t\to-\infty$, from (\ref{eq:11}) we deduce that
$\int_{\widetilde D}|\nabla u|^2=0$ and hence $u$ must be constant on
$\widetilde D$.  Since $u$ vanishes on $\partial \widetilde D$, we
deduce that $u\equiv 0$ and then $u_1=u_2$ in  $\widetilde D$, thus proving 
uniqueness.
\end{pf}

\begin{remark}\label{rem:phi2pos}
From \eqref{eq:12} and the Strong Maximum Principle we deduce that 
\begin{equation*}
\Phi_2(x_1,x')>0\quad \text{for all }(x_1,x')\in\widetilde D.
\end{equation*}
\end{remark}

\noindent 
The functions $\Phi_1$ and $\Phi_2$ can be estimated as follows.
\begin{Lemma}\label{l:stima_w}
\quad\\[-7pt]
\begin{itemize}
\item[(i)] For every $\delta>0$ there exists $c(\delta)>0$ such that 
$$
\Big|\Phi_1(x)-(x_1-1)^+\Big|\leq c(\delta)\,\frac{x_1-1}{|x-{\mathbf e}_1|^{N}}
\quad \text{and}\quad
\Phi_2(x)\leq c(\delta)\,\frac{x_1-1}{|x-{\mathbf e}_1|^{N}}
$$
for all $x\in D^+\setminus B_{1+\delta}^+$.\\[-8pt]
\item[(ii)] There exists $C_2>0$ such that 
$$
\Phi_1(x)\leq C_2 e^{\sqrt{\lambda_1(\Sigma)}\frac{x_1-1}2}
\quad\text{for all }x\in T_1^-.
$$
\end{itemize}
\end{Lemma}
\begin{pf}
  Let us first prove (i) for the function
  $w=\Phi_1(x)-(x_1-1)^+=\mathcal T(x_1-1)(x)-(x_1-1)^+$ (the analogous
  estimate for $\Phi_2$ can be proved in a similar way). We observe
  that $w$ belongs to ${\mathcal D}^{1,2}(\widetilde D)$ in view of
  (\ref{eq:38}) and weakly solves $-\Delta w=0$ in $D^1\setminus
  \overline{B_1^+}$ by (\ref{eq:1}); moreover $w(x)=0$ for all
  $x\in\{(x_1,x'):x_1=1,|x -{\mathbf e}_1|>1\}$.  Therefore, its
  Kelvin transform
$$
\widetilde w(x)=|x-{\mathbf e}_1|^{-(N-2)}w\bigg(
\frac{x-{\mathbf e}_1}{|x-{\mathbf e}_1|^2}+{\mathbf e}_1\bigg)
$$
belongs to $H^1(B_1^+)$ and weakly satisfies
$$
\begin{cases}
-\Delta \widetilde w(x)=0,&\text{in }B_1^+,\\
\widetilde w(x)=0,&\text{on }\{(x_1,x'):x_1=1,|x -{\mathbf e}_1|<1\}.
\end{cases}
$$
By classical elliptic estimates, for any $\delta>0$ there exists
$c(\delta)>0$ such that $\big|\frac{\partial \widetilde w}{\partial
  x_1}\big|\leq c(\delta)$ in $\overline{B^+_{1/(\delta+1)}}$, thus implying 
$$
|\widetilde w(x_1,x')|=\bigg|\widetilde w(1,x')+\int_1^{x_1}
\frac{\partial \widetilde w}{\partial
  x_1}(s,x')\,ds\bigg|\leq \int_1^{x_1}
\bigg|\frac{\partial \widetilde w}{\partial
  x_1}(s,x')\bigg|\,ds\leq c(\delta)(x_1-1)
$$
for all $(x_1,x')\in \overline{B^+_{1/(\delta+1)}}$, which implies (i).
To prove (ii), it is enough to observe that, in view of \eqref{eq:condsigma}, the function
$v(x_1,x')=e^{\sqrt{\lambda_1(\Sigma)}\frac{x_1-1}2}\psi_1^\Sigma\big({x'}/2\big)$
is well-defined,  harmonic and strictly positive in $T_1^-$, bounded from below away
from $0$ on $\{(x_1,x')\in T_1^-:x_1=1\}$, and $\int_{T_1^-}(|\nabla
v|^2+|v|^{2^*})<+\infty$. Hence, from the Maximum Principle we deduce
that $\Phi_1(x)\leq {\rm const\,}v(x)$ in $T_1^-$, thus implying statement (ii).
\end{pf}

\noindent 
In order to control $u_\e$ with suitable sub/super-solutions and
obtain the needed upper and lower estimates, let us introduce the
following functions:
\begin{align}\label{eq:26}
&\Phi^\e:D^+\cup T_\e^-\to\R,\quad
\Phi^\e(x)=\e\Phi_1\Big({\mathbf e}_1+\frac{x-{\mathbf e}_1}{\e}\Big)+
2\gamma_\e \e\Phi_2\Big({\mathbf e}_1+\frac{x-{\mathbf e}_1}{2\e}\Big),\\
\label{eq:35}
&\widetilde \Phi^\e:D^+\cup T_\e^-\to\R,\quad
\widetilde \Phi^\e(x)=\e\Phi_1\Big({\mathbf e}_1+\frac{x-{\mathbf e}_1}{\e}\Big)-
\sqrt2\widetilde \gamma_\e \e\Phi_2\Big({\mathbf e}_1+\frac{x-{\mathbf e}_1}{\sqrt2\e}\Big),
\end{align}
where 
\[
\gamma_\e=\bigg(2\e
\exp\Big({\frac{\sqrt{\lambda_1(\Sigma)}}{4\e}}\Big)\bigg)^{-1}, \quad
\widetilde \gamma_\e=\bigg(\sqrt 2\e
\exp\Big({\frac{\sqrt{\lambda_1(\Sigma)}}{2\sqrt2 \e}}\Big)\bigg)^{-1}.
\]
We notice that $\Phi^\e, \widetilde \Phi^\e$ are well-defined in view
of \eqref{eq:condsigma}.
\begin{Lemma}\label{l:sup_sol_bound}
There exists $C_3>0$ such that 
$$
|u_\e(x)|\leq C_3 \Phi^\e(x)\quad\text{for all }\e\in(0,\e_0)\text{ and }
x\in \mathcal B_\e,
$$
where 
\begin{equation}\label{eq:B-eps}
\mathcal B_\e=B_{r_0}^+\cup
\Big\{(x_1,x')\in\R^N:\tfrac{x'}\e\in\Sigma,\,\tfrac12<x_1\leq 1\Big\}.
\end{equation}
\end{Lemma}
\begin{pf}
Let us first observe that 
\begin{equation}\label{eq:33}
-\Delta \Phi^\e=0,\quad\text{in }D^+\cup T_\e^-.
\end{equation}
Moreover, if $x\in \Gamma_{r_0}^+=\partial B_{r_0}^+\cap D^+$ and $\e\in (0,\e_0)$, then 
Lemma \ref{l:U_e_ud} implies that 
\begin{equation}\label{eq:27}
|u_\e(x)|\leq C_0 (x_1-1),
\end{equation}
while (\ref{eq:Phi1geq}) and (\ref{eq:12}) ensure 
\begin{equation}\label{eq:28}
\Phi^\e(x)\geq \e\,\frac{(x_1-1)^+}{\e}=(x_1-1)^+\quad\text{in
}D^+\cup T_\e^+.
\end{equation}
From (\ref{eq:27}--\ref{eq:28}) we deduce that 
\begin{equation}\label{eq:30}
|u_\e(x)|\leq C_0 \Phi^\e(x)\quad\text{for all }x\in \Gamma_{r_0}^+ \text{ and }\e\in (0,\e_0).
\end{equation}
On the other hand, if $x=(x_1,x')\in T_\e^-$ and $x_1=\frac12$, then 
from (\ref{eq:Phi1geq}), (\ref{eq:12}), (\ref{eq:29}), and \eqref{eq:condsigma}, it follows that 
$$
\Phi^\e(x)\geq 2\gamma_\e \e e^{\frac{\sqrt{\lambda_1(\Sigma)}}{4\e}}
\psi_1^\Sigma\Big(\frac{x'}{2\e}\Big)
\geq \min_{\substack{y'\in\R^{N-1}\\|y'|\leq 1/2}}\psi_1^\Sigma(y')=C_4>0
$$
thus implying, in view of Lemma \ref{l:bk}, that 
\begin{equation}\label{eq:31}
|u_\e(x)|\leq \frac{C_1}{C_4}\Phi^\e(x)
\quad\text{for all $x=(x_1,x')\in T_\e^-$ such that $x_1=\frac12$}.
\end{equation}
From (\ref{eq:30}) and (\ref{eq:31}) we conclude that 
\begin{equation}\label{eq:32}
|u_\e(x)|\leq \max\bigg\{C_0,\frac{C_1}{C_4}\bigg\}\Phi^\e(x)
\quad \text{for all }x\in \partial \mathcal B_\e.
\end{equation}
Since, from \eqref{eq:p2} and Kato's inequality, $-\Delta|u_\e|\leq 0$ in $\mathcal B_\e$, from 
(\ref{eq:33}), (\ref{eq:32}), and the Maximum Principle we reach the conclusion.
\end{pf}

\noindent Let us define
\begin{align}\label{eq:84}
\widetilde u_\e:\widetilde \Omega^\e\to\R,\quad
\widetilde u_\e(x)=\frac1\e u_\e\big({\mathbf e}_1
+\e(x-{\mathbf e}_1
)\big),
\end{align}
where 
\begin{equation}\label{eq:85}
\widetilde \Omega^\e:={\mathbf e}_1+\frac{\Omega^\e-{\mathbf
    e}_1}{\e} =\{x\in\R^N:{\mathbf e}_1 +\e(x-{\mathbf e}_1)\in\Omega^\e\}.
\end{equation} 
We observe that $\widetilde u_\e$ solves
\begin{equation}\label{eq:eqtilde}
\begin{cases}
-\Delta \widetilde u_\e(x)=\e^2\lambda^\e_{\bar k} p\big({\mathbf e}_1
+\e(x-{\mathbf e}_1
)\big) \widetilde u_\e(x),&\text{in }\widetilde \Omega^\e,\\
\widetilde u_\e=0,&\text{on }\partial\widetilde \Omega^\e.
\end{cases}
\end{equation}
From Lemma \ref{l:sup_sol_bound}, the following uniform estimate on
the gradient of $u_\e$ on half-annuli with radius of order $\e$ can de derived.
\begin{Lemma}\label{l:grad_est}
For every $1<r_1<r_2<\frac{r_0}{\e_0}$ there exists $C_{r_1,r_2}>0$  such that 
$$
|\nabla u_\e(x)|\leq C_{r_1,r_2}\quad\text{for all }x\in B_{\e
  r_2}^+\setminus B_{\e r_1}^+\text{ and }\e\in (0,\e_0).
$$
\end{Lemma}
\begin{pf}
  From Lemma \ref{l:sup_sol_bound} and (\ref{eq:26}), it follows that,
  letting $\widetilde u_\e$ as in (\ref{eq:84}--\ref{eq:85}),
\begin{align}\label{eq:34}
|\widetilde u_\e(x)|&\leq \frac{C_3}{\e} 
\Phi^\e\big({\mathbf e}_1+\e(x-{\mathbf e}_1)\big)\\
\notag&=C_3\bigg(\Phi_1(x)+2\gamma_\e\Phi_2\Big(\frac{x+{\mathbf e}_1}{2}\Big)\bigg)
\quad\text{for all }x\in B_{r_0/\e}^+,\ \e\in(0,\e_0).
\end{align}
Let us fix $R_1,R_2$ such that $1<R_1<r_1<r_2<R_2<\frac{r_0}{\e_0}$. From (\ref{eq:34}) it follows that 
\begin{align*}
|\widetilde u_\e(x)|&\leq {\rm const}\quad\text{for all }x\in B_{R_2}^+\setminus
 B_{R_1}^+,\ \e\in(0,\e_0),
\end{align*}
for some ${\rm const}>0$ independent of $\e$ (but depending on $R_1,R_2$). 
Hence, from (\ref{eq:eqtilde}) and classical elliptic estimates, we deduce that 
\begin{align*}
  |\nabla \widetilde u_\e(x)|\leq C_{r_1,r_2}\quad\text{for all }x\in B_{r_2}^+\setminus
 B_{r_1}^+,\ \e\in(0,\e_0),
\end{align*}
thus proving the statement.
\end{pf}

\noindent A lower bound for $u_\e$ can be given in terms of the
function $\widetilde\Phi^\e$ defined in (\ref{eq:35}).
\begin{Lemma}\label{l:sub_sol_bound}
There exist $C_5>0$ and $\e_1\in(0,\e_0)$  such that 
$$
u_\e(x)\geq C_5 \widetilde\Phi^\e(x)\quad\text{for all }\e\in(0,\e_1)\text{ and }
x\in \mathcal B_\e,
$$
where $\mathcal B_\e$ is defined in (\ref{eq:B-eps}) and
$\widetilde\Phi^\e$ in (\ref{eq:35}). 
\end{Lemma}
\begin{pf}
Let us first observe that 
\begin{equation}\label{eq:36}
-\Delta \widetilde \Phi^\e=0,\quad\text{in }D^+\cup T_\e^-.
\end{equation}
Moreover, if $x\in \Gamma^+_{r_0}$ and $\e\in (0,\e_0)$, then 
Lemma \ref{l:U_e_ud} implies that 
\begin{equation}\label{eq:37}
u_\e(x)\geq \frac1{C_0} (x_1-1).
\end{equation}
Furthermore, from (\ref{eq:12}) and (\ref{eq:35}), we have that 
\begin{equation}\label{eq:39}
 \widetilde \Phi^\e(x)\leq \e\Phi_1
\Big({\mathbf e}_1+\frac{x-{\mathbf e}_1}{\e}\Big)
\quad\text{for all }x\in D^+\cup T_\e^-.
\end{equation}
From Lemma \ref{l:stima_w}, there exist $C_6,C_7>0$ such that 
\begin{equation}\label{eq:40}
\Phi_1(x)\leq (x_1-1)\bigg(1+\frac{C_6}{|x-{\mathbf e}_1|^N}\bigg)
\leq C_7(x_1-1)\quad\text{for all }x\in D^+\setminus B^+_{2}.
\end{equation}
Combining (\ref{eq:39}) and (\ref{eq:40}), we obtain that 
$$
\widetilde \Phi^\e(x)\leq C_7 (x_1-1)\quad \text{for all }x\in
D^+\setminus B^+_{2\e},
$$
which, together with (\ref{eq:37}), yields
\begin{equation}\label{eq:41}
\widetilde \Phi^\e(x)\leq C_0C_7 u_\e(x)\quad \text{for all }
x\in \Gamma_{r_0}^+\text{ and }0<\e<\e_0.
\end{equation}
On the other hand, if $x=(x_1,x')\in T_\e^-$ and $x_1=\frac12$, then 
(\ref{eq:35}), 
(\ref{eq:39}), Lemma \ref{l:stima_w}(ii),
(\ref{eq:12}),  (\ref{eq:29}), and Lemma \ref{l:bk}  yield
\begin{equation}\label{eq:42}
 \widetilde \Phi^\e(x)\leq C_2\e e^{-\frac{\sqrt{\lambda_1(\Sigma)}}{4\e}}
-\min_{\substack{y'\in\R^{N-1}\\|y'|\leq 1/\sqrt2}}\psi_1^\Sigma(y')\leq -\frac12 
\min_{\substack{y'\in\R^{N-1}\\|y'|\leq 1/\sqrt2}}\psi_1^\Sigma(y')
\leq 
\frac{u_\e(x)}{2C_1}
\min_{\substack{y'\in\R^{N-1}\\|y'|\leq 1/\sqrt2}}\psi_1^\Sigma(y')
\end{equation}
provided $\e$ is sufficiently small.  Estimates (\ref{eq:41}) and
(\ref{eq:42}) imply the existence of some $C_5>0$ and $\e_1>0$ such that 
$$
u_\e(x)\geq C_5 \widetilde\Phi^\e(x)\quad\text{for all }\e\in(0,\e_1)\text{ and }
x\in \partial\mathcal B_\e,
$$
which, together with (\ref{eq:36}) and the  Maximum Principle,
yields the conclusion.
\end{pf}

\begin{Lemma}\label{l:usot}
There exists $\e_2\in(0,\e_1)$ such that 
$$
u_\e(x)\geq \frac{C_5}{2}(x_1-1)\quad\text{for all }
x\in B_{r_0}^+\setminus B_{2\e}^+,\ \e\in(0,\e_2).
$$
\end{Lemma}
\begin{pf}
From (\ref{eq:35}), (\ref{eq:Phi1geq}), and Lemma \ref{l:stima_w},
it follows that, for all $x\in B_{r_0}^+\setminus B_{2\e}^+$,
\begin{align}\label{eq:44}
  \widetilde \Phi^\e(x)&=\e\Phi_1\Big({\mathbf e}_1+\frac{x-{\mathbf e}_1}{\e}\Big)-
\sqrt2\widetilde \gamma_\e \e\Phi_2\Big({\mathbf e}_1+\frac{x-{\mathbf e}_1}{\sqrt2\e}\Big)\\
&\notag\geq (x_1-1)- {\rm const\,}\widetilde\gamma_\e (x_1-1)
\geq \frac12(x_1-1)
\end{align}
provided $\e$ is sufficiently small. The conclusion  follows
from Lemma \ref{l:sub_sol_bound} and (\ref{eq:44}).
\end{pf}

\section{The frequency function}\label{sec:freq-funct}

In this section we introduce an Almgren type quotient associated to
problem \eqref{eq:24} and study its monotonicity properties with the
aim of uniformly controlling the transversal frequencies along the
connecting tube. 

For every $\e>0$, let $\xi_\e:\R\setminus
\big((-\e,0)\cup(1,1+\e) \big) \to \R$ be such that 
$$
\begin{cases}
\xi_\e(r)=-r,&\text{if }r\leq -\e,\\
\xi_\e(r)=r,&\text{if }0\leq r\leq1,\\
\xi_\e(r)=r-1,&\text{if } r\geq \e+1.
\end{cases}
$$
For $r\in \R\setminus
\big((-\e,0)\cup(1,1+\e) \big)$, 
we define
\begin{align*}
\Omega_r^\e&=
\begin{cases}
D^-\setminus \overline{B^-_{\xi_\e(r)}},&\text{if }r\leq -\e,\\
D^-\cup\{(x_1,x')\in \mathcal C_\e:x_1<r\},&\text{if }0\leq r\leq1,\\
D^-\cup\mathcal C_\e\cup B^+_{\xi_\e(r)},&\text{if } r\geq \e+1,
\end{cases}\\
\Gamma_r^\e&=
\begin{cases}
D^-\cap\partial B^-_{\xi_\e(r)},&\text{if }r\leq -\e,\\
\{(x_1,x')\in \mathcal C_\e:x_1=r\},&\text{if }0\leq r\leq1,\\
D^+\cap \partial  B^+_{\xi_\e(r)},&\text{if } r\geq \e+1.
\end{cases}
\end{align*}
We also denote 
\begin{equation}\label{eq:defOmega_r}
\Omega_r:=D^-\setminus \overline{B^-_{-r}}\quad \text{for all }r<0
\end{equation}
and notice that 
$$
\Omega_r^{\e}=\Omega_r \quad \text{for all }r\leq -\e.
$$

\begin{figure}[h]
  \centering \subfigure [$\Omega^\e_{r}$ with $r\leq-\e$]{
     \includegraphics[width=4cm]{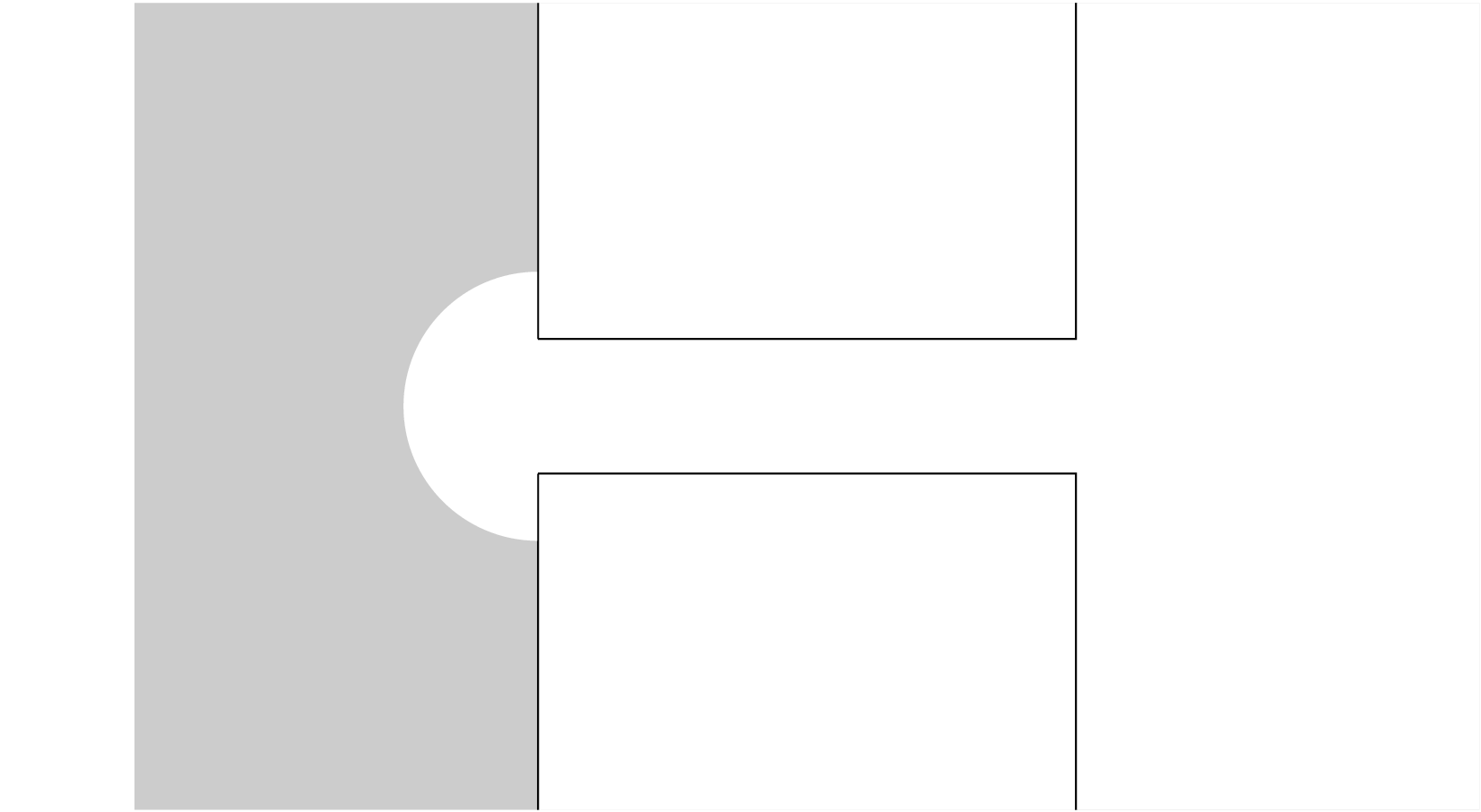}
     } 
\quad \subfigure [$\Omega^\e_{r}$ with $0\leq r\leq1$]{
     \includegraphics[width=4cm]{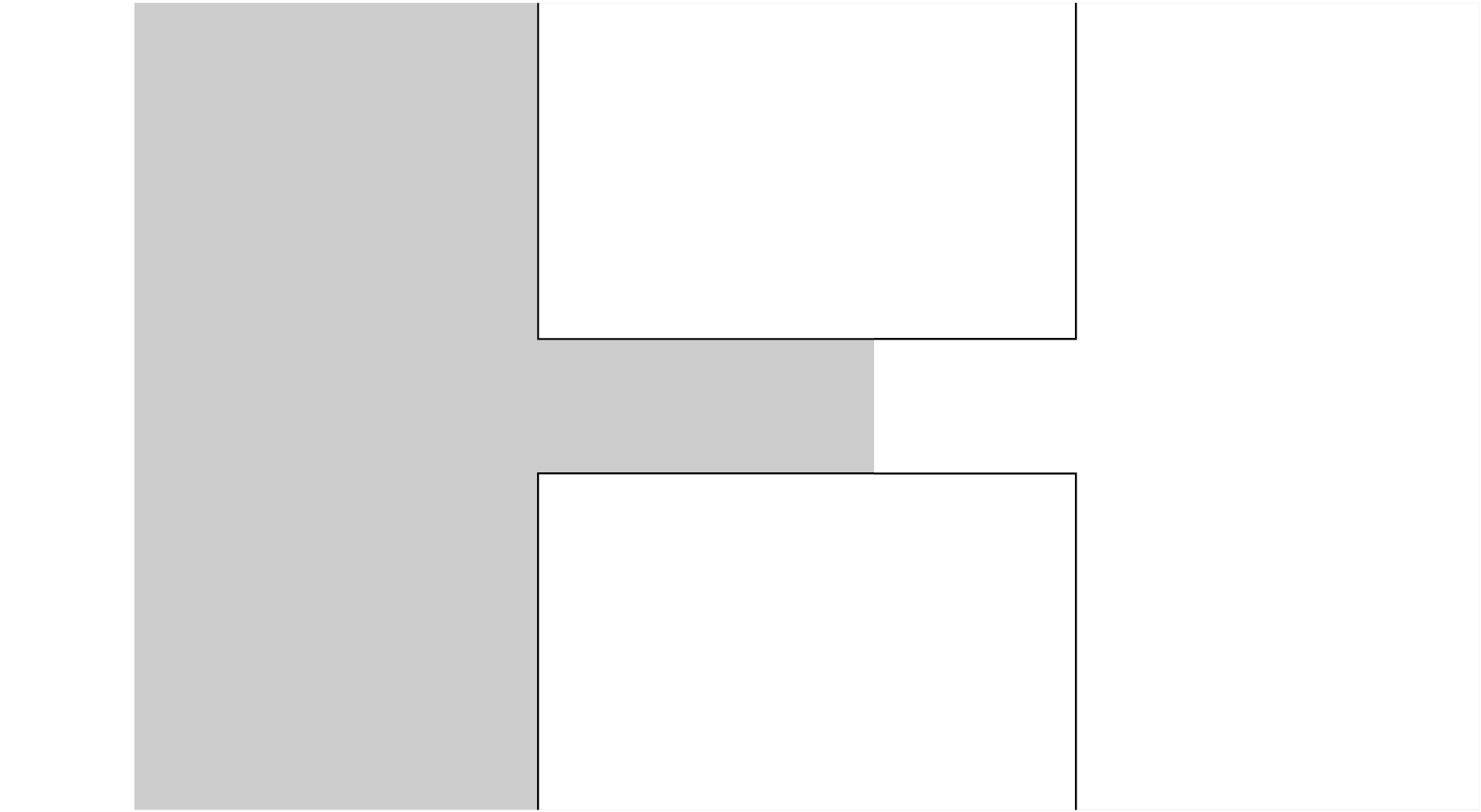}
     }
\quad \subfigure[$\Omega^\e_{r}$ with $r\geq\e$ ] {
    \includegraphics[width=4cm]{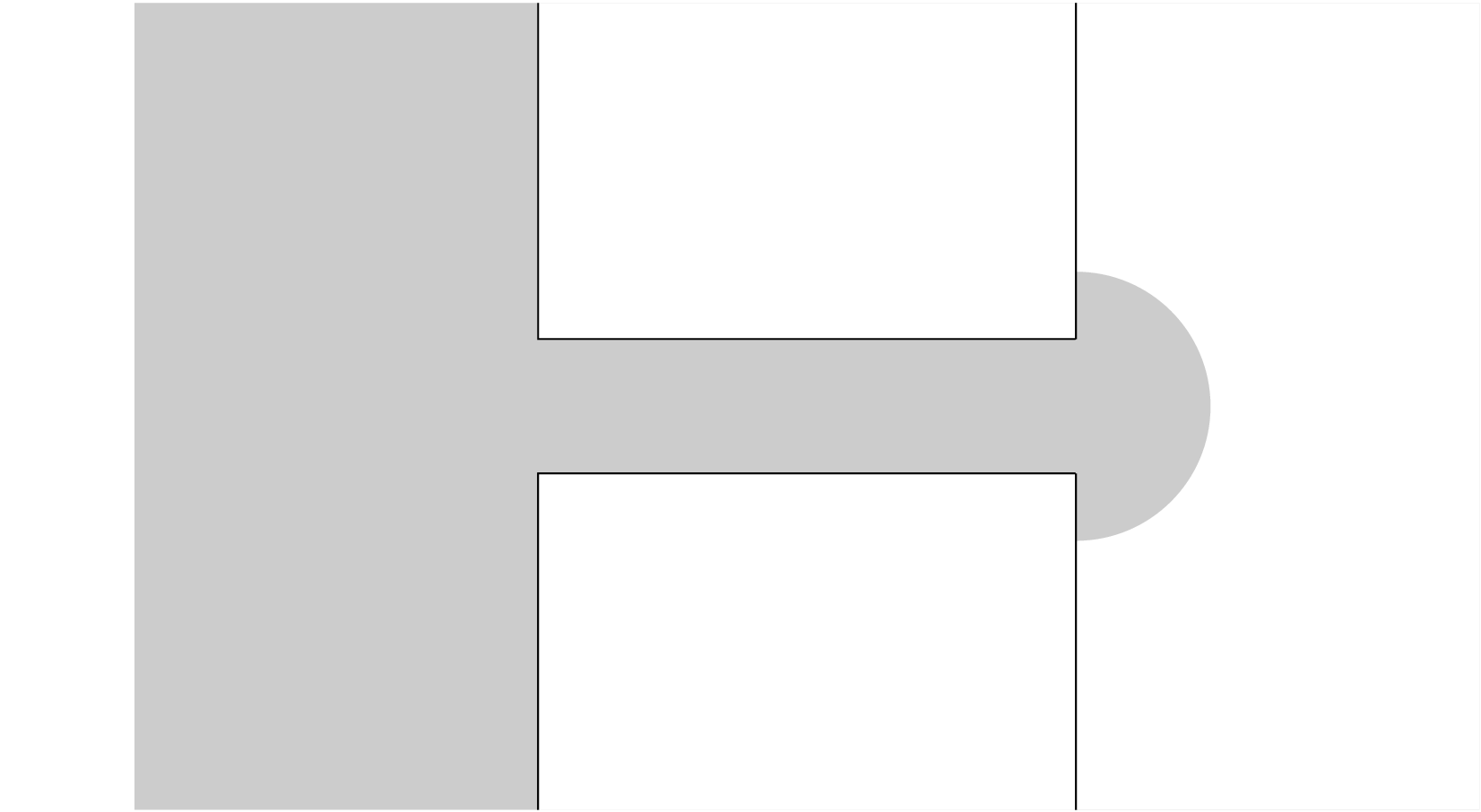}
  }
  \caption{The moving domains $\Omega^\e_{r}$ for different values of the parameters.}\label{fig:dd}
\end{figure}

A key role in the definition and in the study of the frequency
associated  to
problem \eqref{eq:24} is played by
Lemmas \ref{l:poincareN-1} and \ref{l:lemmaZERO} below, which give a
Poincar\'e type lemma on domains $\Omega_{-t}$, $t>0$, for functions
in $\mathcal H_t^-$ and, respectively, a uniform coercivity type
estimate for the quadratic form associated to equation \eqref{eq:24}
in domains $\Omega_r^\e$, $r<1$. An important ingredient for their
proof is the Kelvin transform, which is described in the following remark. 
\begin{remark}\label{rem:kelvin}
For all $R>0$,
  $v\in \mathcal H_R^-$ if and only if its Kelvin transform
$\widetilde v(x)=|x|^{-(N-2)}v\big(
\frac{x}{|x|^2}\big)$
belongs to $H^1(B_{1/R}^-)$ and has null trace on $\partial B_{1/R}^-\cap \partial D^-$; furthermore 
\begin{align*}
  &\int_{B_{1/R}^-}|\nabla \widetilde
  v(x)|^2dx+(N-2)R\int_{\Gamma_{1/R}^-}\widetilde v^2d\sigma
  =\int_{\Omega_{-R}}|\nabla v(x)|^2dx,\\
  &\int_{B_{1/R}^-}|\widetilde
  v(x)|^{2^*}dx=\int_{\Omega_{-R}}|v(x)|^{2^*}dx,\quad\text{and}\quad
  R^2\int_{\Gamma_{1/R}^-}\widetilde v^2(x)\,d\sigma
  =\int_{\Gamma_R^-}v^2(x)\,d\sigma.
\end{align*}
\end{remark}

\noindent Functions in $\mathcal H_t^-$ satisfy the following Sobolev
type inequality.

\begin{Lemma}\label{l:sob_ex}
  There exists a constant $C_{S}=C_S(N)$ depending only on the dimension $N$
  such that  for all $t>0$ and $v\in\mathcal H_t^-$ there holds
$$
C_S\bigg(\int_{\Omega_{-t}}|v(x)|^{2^*}dx\bigg)^{\!\!2/2^*}\leq 
\int_{\Omega_{-t}}|\nabla v(x)|^2dx.
$$
\end{Lemma}
\begin{pf}
  By scaling it is enough to prove the inequality for $t=1$, which, in
  view of remark \ref{rem:kelvin}, is equivalent to prove that
$$
C_S\bigg(\int_{B_1^-}|w(x)|^{2^*}dx\bigg)^{\!\!2/2^*}\leq
\int_{B_1^-}|\nabla w(x)|^2dx+(N-2)\int_{\Gamma_1^-}w^2d\sigma
$$
for all $w\in H^1(B_1^-)$ such that $w\equiv 0$ on $\partial
B_1^-\cap \partial D^-$. Such inequality follows easily from classical
Sobolev embeddings by trivially extending $w$ in $B({\mathbf 0},1)$ and
observing that
$$
\int_{B({\mathbf 0},1)}|\nabla w(x)|^2dx+(N-2)\int_{\partial B({\mathbf 0},1)}w^2d\sigma
$$
 is an equivalent norm in $H^1(B({\mathbf 0},1))$.
\end{pf}

The Poincaré inequality we will state in Lemma \ref{l:poincareN-1}
with its best constant is a consequence of the following lemma, which
is the counterpart of Lemma \ref{l:limNphi1gen} for the frequency of harmonic
functions in $\mathcal H_t^-$.

\begin{Lemma}\label{l:limNmenophi1gen}
Let  $R>0$ and $\phi\in\mathcal H_R^-\setminus\{0\}$ satisfying  
\begin{align*}
\begin{cases}
-\Delta \phi=0,&\text{ in }\Omega_{-R},\\
\phi=0,&\text{ on }\partial\Omega_{-R}\cap\partial D^-,
\end{cases}
\end{align*}
in a weak sense, and let  $N_\phi^-:(R,+\infty)\to\R$ be defined as 
$$
N_\phi^-(r):=\frac{r\int_{\Omega_{-r}}|\nabla \phi(x)|^2dx}
{\int_{\Gamma_r^-}\phi^2(x)\,d\sigma}.
$$
Then 
\begin{itemize}
\item[i)]
$N_\phi^-$ is non-increasing in $(R,+\infty)$;
\item[ii)]
there exists $K_0\in\N$, $K_0\geq1$, such that 
$$
\lim_{r\to\infty}N_\phi^-(r)=
N-2+K_0;
$$
\item[iii)] if $N_\phi^-\equiv \gamma$ for some $\gamma\in\R$ then
  $\gamma=N-2+K_0$ and $\phi(x)=|x|^{-N+2-K_0}Y(x/|x|)$ for some
  eigenfunction $Y$ of $-\Delta_{{\mathbb S}^{N-1}}$ associated to the
  eigenvalue $K_0(N-2+K_0)$, i.e.  satisfying $-\Delta_{{\mathbb
      S}^{N-1}}Y=K_0(N-2+K_0) Y$ on ${\mathbb S}^{N-1}$;
\item[iv)] if $\phi>0$  in $\Omega_{-R}$, then $K_0=1$.
\end{itemize}
\end{Lemma}
\begin{pf}
  Let $\widetilde\phi\in H^1(B_{1/R}^-)$ be the Kelvin transform of
  $\phi$, i.e.  $\widetilde \phi(x)=|x|^{-(N-2)}\phi\big(
  \frac{x}{|x|^2}\big)$. Then 
$\widetilde\phi$ satisfies 
$$
\begin{cases} 
-\Delta \widetilde\phi=0,&\text{in }B_{1/R}^-,\\
\widetilde\phi=0,&\text{ on }\partial B_{1/R}^-\cap \partial D^-,
\end{cases}
$$
and, by Remark \ref{rem:kelvin}, the
  frequency function $N_\phi^-$ can be rewritten as
\begin{equation}\label{eq:144}
N_\phi^-(r)=N-2+\widetilde N\bigg(\frac1r\bigg),
\end{equation}
where 
$$
\widetilde N:\bigg(0,\frac1R\bigg)\to\R,\quad
\widetilde N(t):=
\frac{t\int_{B_t^-}|\nabla \widetilde\phi(x)|^2dx}{\int_{\Gamma_t^-}\widetilde\phi^2d\sigma}.
$$
Let us define
$$
\widetilde \phi_0(x)=\widetilde \phi_0(x_1,x')=
\begin{cases}
\widetilde \phi(x_1,x'),&\text{if }x_1\leq0,\\
-\widetilde \phi(-x_1,x'),&\text{if }x_1>0,\\
\end{cases}
$$
and observe that $\widetilde \phi_0\in H^1(B({\mathbf 0},1/R))$ satisfies 
$\widetilde \phi_0(-x_1,x')=-\widetilde \phi_0(x_1,x')$ and weakly solves 
$$
-\Delta\widetilde \phi_0=0,\quad\text{in }B({\mathbf 0},1/R).
$$
Moreover 
\begin{equation}\label{eq:149}
\widetilde N(t)=
\frac{t\int_{B({\mathbf 0},t)}|\nabla \widetilde \phi_0(x)|^2dx}
{\int_{\partial B({\mathbf 0},t)}\widetilde \phi_0^2d\sigma}.
\end{equation}
From the classical Almgren monotonicity formula \cite{almgren}
\begin{align}\label{eq:147}
\widetilde N'(t) = \frac{2t\Big[ \Big(\int_{\partial B({\mathbf 0},t)}
\Big|\frac{\partial \widetilde\phi_0}{\partial\nu}\Big|^2
        d\sigma\Big) \Big(\int_{\partial B({\mathbf 0},t)} \widetilde\phi_0^2
        d\sigma\Big)
-\Big(\int_{\partial B({\mathbf 0},t)}\widetilde\phi_0
        \frac{\partial\widetilde\phi_0}{\partial \nu}
        d\sigma\Big)^{2} \Big]} {\Big(\int_{\partial B({\mathbf 0},t)} \widetilde\phi_0^2
      d\sigma\Big)^{2}},
\end{align}
for all $t\in(0,1/R)$, where $\nu=\nu(x)=\frac{x}{|x|}$, hence from
Schwarz's inequality $\widetilde N'\geq0$ and  the function
$t\in(0,1/R)\mapsto \widetilde N(t)$ is non-decreasing, thus implying,
in view of (\ref{eq:144}), that $N_\phi^-$ is non-increasing in
$(R,+\infty)$ and proving statement i). Furthermore from \cite[Theorem 1.3]{FFT} there exist
$K_0\in\N$ and an eigenfunction $Y$ of $-\Delta_{{\mathbb S}^{N-1}}$
associated to the eigenvalue $K_0(N-2+K_0)$, i.e.  satisfying
$-\Delta_{{\mathbb S}^{N-1}}Y=K_0(N-2+K_0) Y$ on ${\mathbb S}^{N-1}$,
such that
\begin{equation}\label{eq:146}
\lim_{t\to0^+}\widetilde N(t)=-\frac{N-2}{2}
+\sqrt{\bigg(\frac{N-2}{2}\bigg)^{\!\!2}+K_0(N-2+K_0)}=K_0
\end{equation}
and 
\begin{align}
\label{eq:62}&\lambda^{-K_0}\widetilde \phi_0(\lambda\theta)
\to Y(\theta)\quad\text{in }C^{1,\tau}({\mathbb S}^{N-1}),\\
&\label{eq:63}\lambda^{1-K_0}\nabla \widetilde \phi_0(\lambda\theta)\to
K_0 Y(\theta)\theta+\nabla_{{\mathbb S}^{N-1}}
Y(\theta)\quad \text{in }C^{0,\tau}({\mathbb S}^{N-1}),
\end{align}
as $\lambda\to 0^+$,
for every $\tau\in(0,1)$. 
Since $\widetilde \phi_0$ vanishes on 
$B({\mathbf 0},1)\cap(\{0\}\times\R^{N-1})$, from (\ref{eq:62}) we infer that  $Y$ vanishes on  the
equator ${\mathbb S}^{N-1}\cap(\{0\}\times\R^{N-1})$. Therefore, $Y$ can not be the first eigenfunction 
of $-\Delta_{{\mathbb S}^{N-1}}$ and hence 
 $K_0\geq 1$ necessarily. 
Statement ii) then follows from (\ref{eq:144}) and (\ref{eq:146}).

Let us now assume that $N_\phi^-\equiv \gamma$ for some $\gamma\in\R$, so that $\widetilde N(t)\equiv\gamma-N+2$
 in $(0,1/R)$ and hence $\widetilde N'(t)=0$ for any $t\in (0,1/R)$.  By (\ref{eq:147}) we obtain
\begin{equation*}
  \bigg(\int_{\partial B({\mathbf 0},t)} \bigg|\frac{\partial
        \widetilde\phi_0}{\partial\nu}\bigg|^2 d\sigma\bigg)
 \cdot \bigg(\int_{\partial B({\mathbf 0},t)} \widetilde\phi_0^2
    d\sigma\bigg)-\bigg(\int_{\partial B({\mathbf 0},t)} \widetilde\phi_0\frac{\partial \widetilde\phi_0}
{\partial \nu} d\sigma\bigg)^{\!\!2}=0 \quad 
  \text{for all } t\in (0,1/R),
\end{equation*}
i.e.
$ \widetilde\phi_0$ and $\frac{\partial \widetilde\phi_0}
{\partial \nu}$ 
have the same direction as vectors in $L^2(\partial B({\mathbf 0},t))$ and hence
there exists a  function $\eta=\eta(t)$ such that
$\frac{\partial \widetilde\phi_0}
{\partial \nu}(t,\theta)=\eta(t) \widetilde\phi_0(t,\theta)$ for $t\in(0,1/R)$ and $\theta\in{\mathbb S}^{N-1}$.
After integration we obtain
\begin{equation}\label{separate}
  \widetilde\phi_0(t,\theta)=e^{\int_{1/R}^t \eta(s)ds} \widetilde\phi_0\bigg(\frac1R,\theta\bigg)
  =\varphi(t) \psi(\theta), \quad  t\in(0,1/R), \ \theta\in {\mathbb S}^{N-1},
\end{equation}
where  $\varphi(t)=e^{\int_{1/R}^t \eta(s)ds}$ and $\psi(\theta)=
 \widetilde\phi_0\big(\frac1R,\theta\big)$.
Since  $-\Delta\widetilde\phi_0=0$ in $B({\mathbf 0},1/R)$,  (\ref{separate}) yields
$$
\left(-\varphi''(t)-\frac{N-1}{t} \varphi'(t) \right)\psi(\theta)
-\frac{\varphi(t)}{t^2}\Delta_{{\mathbb S}^{N-1}} \psi(\theta)=0.
$$
Taking $t$ fixed we deduce that $\psi$ is an eigenfunction of the
operator $-\Delta_{{\mathbb S}^{N-1}}$. If $K_0(N-2+K_0)$ is the corresponding
eigenvalue then $\varphi(t)$ solves the equation
$$
-\varphi''(t)-\frac{N-1}{t} \varphi(t)+\frac{K_0(N-2+K_0)}{t^2}\varphi(t)=0
$$
and hence $\varphi(t)$ is of the form
$$
\varphi(t)=c_1 t^{K_0}+c_2 t^{-(N-2)-K_0}, \quad\text{for some $c_1,c_2\in\R$}.
$$
Since the function
$|x|^{-(N-2)-K_0}\psi(\frac{x}{|x|})\notin
H^1(B_{1/R})$, then $c_2=0$ and  $\varphi(t)=c_1
t^{K_0}$. Since $\varphi(\frac1R)=1$, we obtain that $c_1=R^{K_0}$ and
then
\begin{equation}\label{eq:148}
\widetilde\phi_0(t,\theta)=R^{K_0}t^{K_0} \psi(\theta),  \quad 
\text{for all }t\in (0,1/R)\text{ and }\theta\in \SN.
\end{equation}
Therefore $\phi(y)=R^{K_0}|y|^{-N+2-K_0}\psi(\frac{y}{|y|})$ in $\Omega_{-R}$.
Substituting (\ref{eq:148}) into (\ref{eq:149}) and taking into account that
 $\widetilde N(t)\equiv\gamma-N+2$, we obtain that necessarily $\gamma-N+2=K_0$, i.e. $\gamma=N-2+K_0$.
Claim iii) is thereby proved.

If $\phi>0$  in $\Omega_{-R}$, then
$\widetilde\phi>0$ in $B_{1/R}^-$, and 
Hopf's Lemma implies that 
\begin{equation}\label{eq:61}
\frac{\partial \widetilde\phi}{\partial x_1}(0,x')<0,\quad\text{for all }x'\in\R^{N-1}\text{ s.t. }
|x'|<\frac1R.
\end{equation}
(\ref{eq:61}) and (\ref{eq:63}) 
imply that $K_0\leq 1$. Hence $K_0=1$ and statement iv) is proved.
\end{pf}

We are now ready to prove the following Poincaré type inequality.
\begin{Lemma}\label{l:poincareN-1}
For all $t>0$ and $v\in\mathcal H_t^-$ there holds
$$
\frac1{t^{N-2}}\int_{\Omega_{-t}}|\nabla
v(x)|^2dx\geq\frac{N-1}{t^{N-1}} \int_{\Gamma_t^-}v^2d\sigma,
$$
being $N-1$ the optimal constant.
\end{Lemma}
\begin{pf}
  By scaling it is enough to prove the inequality for $t=1$, i.e. the
  statement of the lemma is equivalent to prove that the infimum
\begin{equation*}
\mathcal I:=\inf_{w\in \mathcal H_1^-\setminus\{0\}}
\frac{\int_{\Omega_{-1}}|\nabla w(x)|^2dx}{\int_{\Gamma_1^-}w^2d\sigma}
\end{equation*}
is equal to $N-1$.  By standard minimization arguments and compactness
of the embedding $\mathcal H_1^-\hookrightarrow L^2(\Gamma_1^-)$, it
is easy to prove that the infimum $\mathcal I$ is strictly positive
and attained by some function $w_0\in \mathcal H_1^- \setminus\{0\}$
satisfying
$$
\begin{cases} 
-\Delta w_0=0,&\text{in }\Omega_{-1},\\
w_0>0,&\text{in }\Omega_{-1},\\
\frac{\partial w_0}{\partial \nu}=-\mathcal I w_0,&\text{on }\Gamma_1^-,\\
w_0\equiv 0,&\text{on }\partial \Omega_{-1} \cap \partial D^-,
\end{cases}
$$
being $\nu=\frac{x}{|x|}$. Then  Lemma \ref{l:limNmenophi1gen} implies that 
$$
\mathcal I=\frac{\int_{\Omega_{-1}}|\nabla w_0(x)|^2dx}{\int_{\Gamma_1^-}w_0^2d\sigma}=N^-_{w_0}(1)
\geq \lim_{r\to +\infty}N^-_{w_0}(r)\geq N-1.
$$
On the other hand the quotient $\big(\int_{\Omega_{-1}}|\nabla
w(x)|^2dx\big)\big(\int_{\Gamma_1^-}w^2d\sigma\big)^{-1}$ evaluated in
$w(x_1,x')=\frac{x_1}{|x|^N}$ is equal to $N-1$, thus implying that
$\mathcal I\leq N-1$.
\end{pf}
\begin{remark}\label{rem:poincare_inv}
 By remark \ref{rem:kelvin},
Lemma \ref{l:poincareN-1} is equivalent to
\begin{equation*}
r\int_{B_r^-}|\nabla w(x)|^2dx\geq \int_{\Gamma_r^-}w^2d\sigma
\quad\text{
for all $w\in H^1(B_r^-)$ such that $w\equiv 0$ on $\partial
B_r^-\cap \partial D^-$}.
\end{equation*}
\end{remark}
\noindent Lemma \ref{l:lemmaZERO} below provides a uniform coercivity type
estimate for the quadratic form associated to equation \eqref{eq:24},
whose validity is strongly related to the nondegeneracy
condition \eqref{eq:54}.

\begin{Lemma}\label{l:lemmaZERO}
\

\begin{itemize}
\item[i)] For every $f\in L^{N/2}(\R^N)$ and $M>0$, there exist $\tilde
  r_{M,f}>0$ and $\tilde\e_{M,f}\in(0,\e_0)$ such that for all $\e\in
  (0,\tilde\e_{M,f})$ and $r\in(\e,\tilde r_{M,f})$
$$
\int_{\Omega_{-r}}
|\nabla u_\e(x)|^2dx\geq M \int_{\Omega_{-r}}
|f(x)| u_\e^2(x)dx.
$$
\item[ii)] For every $f\in L^{N/2}(\R^N)$ and $M>0$, there exists 
 $\bar\e_{M,f}\in(0,\e_0)$ such that for all $r\in(0,1)$ and $\e\in
  (0,\bar\e_{M,f})$ 
$$
\int_{\Omega_{r}^\e}
|\nabla u_\e(x)|^2dx\geq M \int_{\Omega_{r}^\e}
|f(x)| u_\e^2(x)dx.
$$
\end{itemize}
\end{Lemma}
\begin{pf}
  To prove i), we argue by contradiction and assume that
  there exist $f\in L^{N/2}(\R^N)$, $M>0$, and sequences $\e_n\to
  0^+$, $r_n\to 0^+$, such that $r_n>\e_n$ and, denoting
  $u_n=u_{\e_n}$,
\begin{equation}\label{eq:50}
  \int_{ \Omega_{-r_n}} |\nabla u_n(x)|^2dx<
  M \int_{\Omega_{-r_n} }
  |f(x)| u_n^2(x)dx.
\end{equation}
Let us define
$$
v_n(x)=
\begin{cases}
u_n(x),&\text{if }x\in \Omega_{-r_n},\\
\big( \frac{r_n}{|x|}\big)^{N-2}u_n \big( \frac{r_n^2x}{|x|^2}\big)
,&\text{if }x\in B^-_{r_n}.
\end{cases}
$$
We notice that $v_n\in {\mathcal D}^{1,2}(D^-)$ and, by Remark \ref{rem:kelvin}, 
$$
\int_{B^-_{r_n}}|\nabla
v_n(x)|^2dx+\frac{N-2}{r_n}\int_{\Gamma_{r_n}^-} v_n^2d\sigma=\int_{\Omega_{-r_n}}|\nabla u_n(x)|^2dx,
$$
thus implying
\begin{equation}\label{eq:51}
\int_{D^-}|\nabla v_n(x)|^2dx\leq
\int_{D^-}|\nabla v_n(x)|^2dx+\frac{N-2}{r_n}\int_{\Gamma_{r_n}^-}v_n^2d\sigma =2\int_{\Omega_{-r_n}}
|\nabla u_n(x)|^2dx.
\end{equation}
From (\ref{eq:50}) and (\ref{eq:51}) it follows that, if
$$
w_n=\frac{v_n}{\big(\int_{\Omega_{-r_n} }
|f(x)| u_n^2(x)dx\big)^{1/2}},
$$
then $w_n\in {\mathcal D}^{1,2}(D^-)$ and 
$$
\int_{D^-}|\nabla w_n(x)|^2dx<2M.
$$
Hence there exists a subsequence $\{w_{n_k}\}_k$ such that 
$$
w_{n_k}\weakly w\quad\text{weakly in }{\mathcal D}^{1,2}(D^-)
$$
for some $w\in {\mathcal D}^{1,2}(D^-)$.  From $\int_{D^-} |f(x)|
v_n^2(x)dx\geq \int_{\Omega_{-r_n}} |f(x)|
u_n^2(x)dx$ we deduce that
$$
\int_{D^-} |f(x)| w_n^2(x)dx\geq 1
$$
which implies that $w\not\equiv 0$. Since $w_n$ solves 
\begin{align*}
\begin{cases}
-\Delta w_n=\lambda^{\e_n}_{\bar k} p w_n,&\text{in }\Omega_{-r_n},\\
w_n=0,&\text{on }\partial D^-,
\end{cases}
\end{align*}
and $r_n\to 0^+$, from (\ref{eq:52}) we conclude that $w$ weakly solves 
\begin{align*}
\begin{cases}
-\Delta w_=\lambda_{k_0}(D^+) p w,&\text{in }D^-,\\
w=0,&\text{on }\partial D^-,
\end{cases}
\end{align*}
thus implying $\lambda_{k_0}(D^+)\in \sigma_p(D^-)$ and contradicting
assumption (\ref{eq:54}).

Let us now prove ii).  We argue by contradiction and assume that
there exist $f\in L^{N/2}(\R^N)$, $M>0$, and sequences $\e_n\to 0^+$,
$r_n\in (0,1)$, such that denoting $u_n=u_{\e_n}$,
\begin{equation}\label{eq:80}
\int_{ \Omega_{r_n}^{\e_n} } |\nabla u_n(x)|^2dx<
M \int_{ \Omega_{r_n}^{\e_n}}|f(x)| u_n^2(x)dx.
\end{equation}
Let us define
$$
v_n(x)=
\begin{cases}
u_n(x),&\text{if }x\in \Omega_{r_n}^{\e_n},\\
u_n(2r_n{\mathbf e}_1-x)
,&\text{if }2r_n{\mathbf e}_1-x\in \Omega_{r_n}^{\e_n},\\
0,&\text{otherwise}.
\end{cases}
$$
We notice that $v_n\in {\mathcal D}^{1,2}(\R^N)$ and, by (\ref{eq:80}), 
$$
\int_{\R^N}|\nabla v_n(x)|^2dx=2\int_{\Omega_{r_n}^{\e_n}}|\nabla u_n(x)|^2dx<
2M \int_{ \Omega_{r_n}^{\e_n}}|f(x)| u_n^2(x)dx,
$$
thus implying that, letting
$$
w_n=\frac{v_n}{\big(\int_{\Omega_{r_n}^{\e_n}} 
|f(x)| u_n^2(x)dx\big)^{1/2}},
$$
then $w_n\in {\mathcal D}^{1,2}(\R^N)$ and 
$$
\int_{\R^N}|\nabla w_n(x)|^2dx<2M.
$$
Hence there exist a subsequence $\{w_{n_k}\}_k$ and some $w\in
{\mathcal D}^{1,2}(\R^N)$ such that $w_{n_k}\weakly w$ weakly in
${\mathcal D}^{1,2}(\R^N)$ and a.e. in $\R^N$. 
From
\begin{align*}
&1=\int_{\Omega_{r_n}^{\e_n}} 
|f(x)| w_n^2(x)dx
=\int_{D^-} 
|f(x)| w_n^2(x)dx+\int_{\Omega_{r_n}^{\e_n}\setminus D^-} 
|f(x)| w_n^2(x)dx,\\
&\int_{\Omega_{r_n}^{\e_n}\setminus D^-} 
|f(x)| w_n^2(x)dx
\leq \|w_n\|_{L^{2^*}(\R^N)}^2\|f\|_{L^{N/2}(\Omega_{r_n}^{\e_n}\setminus D^-)}
=o(1)\quad\text{as }n\to+\infty,
\\
&\int_{D^-} 
|f(x)| w_{n_k}^2(x)dx=\int_{D^-} 
|f(x)| w^2(x)dx+o(1)\quad\text{as }k\to+\infty,
\end{align*}
we deduce that 
$$
\int_{D^-} 
|f(x)| w^2(x)dx=1
$$
and hence $w\not\equiv 0$ in $D^-$. On the other hand, a.e. convergence
of $w_{n_k}$ to $w$ implies that $w=0$ on $\partial D^-$. Furthermore, passing to the weak limit in the equation 
$-\Delta w_{n_k}=\lambda^{\e_{n_k}}_{\bar k} p w_{n_k}$ satisfied by $w_{n_k}$ in $D^-$,  we conclude that $w$ weakly solves 
\begin{align*}
\begin{cases}
-\Delta w_=\lambda_{k_0}(D^+) p w,&\text{in }D^-,\\
w=0,&\text{on }\partial D^-,
\end{cases}
\end{align*}
thus implying $\lambda_{k_0}(D^+)\in \sigma_p(D^-)$ and contradicting
assumption (\ref{eq:54}).
\end{pf}

\noindent From Lemma \ref{l:lemmaZERO} and \eqref{eq:p2}, there exist $\check R\in(0,1)$ and
  $\check\e\in(0,\e_0)$ such that, for every $\e\in (0,\check\e)$, 
\begin{align}\label{eq:coercivity}
\int_{\Omega_r^\e} \Big(|\nabla
u_\e|^2-\lambda^\e_{\bar k} p u_\e^2\Big)dx \geq
\frac12\int_{\Omega_r^\e} |\nabla u_\e|^2dx \quad\text{for all }r\in(-\check R,-\e)\cup(0,1),
\end{align}
and 
\begin{align}\label{eq:coercivity2}
\int_{\Omega_r^\e} \Big(|\nabla
u_\e|^2-\lambda^\e_{\bar k} p u_\e^2\Big)dx \geq
\frac12\int_{\Omega_{1/2}^\e} |\nabla u_\e|^2dx \quad\text{for all }r\in(1+\e,4).
\end{align}
Estimates \eqref{eq:coercivity} and \eqref{eq:coercivity2}, together
with equation \eqref{eq:24} and classical unique continuation
principle, imply that 
$$
\int_{\Gamma_r^\e}u_\e^2(x)\,d\sigma>0\quad\text{for all }\e\in (0,\check\e)\text{ and }r\in 
 (-\check R,-\e)\cup(0,1)\cup
(1+\e,4).
$$
Therefore, for all $\e\in (0,\check\e)$, the frequency function $\mathcal N_\e:(-\check R,-\e)\cup(0,1)\cup
(1+\e,4)\to\R$,
$$
\mathcal N_\e(r)=\frac{
\Lambda_N(r,\e)\int_{\Omega_r^\e}\Big(|\nabla u_\e(x)|^2-\lambda^\e_{\bar k} p(x) u_\e^2(x)\Big)dx
}{\int_{\Gamma_r^\e}u_\e^2(x)\,d\sigma},
$$
where 
$$
\Lambda_N(r,\e)=
\begin{cases}
\big(\frac2{\omega_{N-1}}\big)^{\!\frac1{N-1}}
|\Gamma_r^\e|^{\frac1{N-1}}=\xi_\e(r),&\text{if }r\in(-\infty,-\e)\cup(1+\e,+\infty),\\[3pt]
\big(\frac {N-1}{\omega_{N-2}}\big)^{\!\frac1{N-1}}|\Gamma_r^\e|^{\frac1{N-1}}=\e,&\text{if }r\in[0,1],
\end{cases}
$$
and $|\Gamma_r^\e|$ denotes the $(N-1)$-dimensional volume of
$\Gamma_r^\e$, is well defined.

\subsection{The frequency function at the right}\label{sec:freq-funct-at}

If $\e\in (0,\check\e)$ and $r\in(1+\e,4)$, then
\begin{equation}\label{eq:N_N+}
\mathcal N_\e(r)=\mathcal N_\e^+(r-1)
\end{equation}
where, for $t\in(\e,3)$, 
\begin{align}
  \label{eq:47} \mathcal N_\e^+(t)&= 
  \frac{D_\e^+(t)}{H_\e^+(t)},\\
  \notag D_\e^+(t)&=\frac1{t^{N-2}} \int_{D^-\cup\mathcal C_\e\cup
    B^+_{t}}\Big(|\nabla
  u_\e(x)|^2-\lambda^\e_{\bar k} p u_\e^2(x)\Big)dx,\\
  \notag
  H_\e^+(t)&=\frac1{t^{N-1}}\int_{\Gamma_t^+}u_\e^2(x)\,d\sigma,
\end{align}
with $\Gamma_t^+$ as defined in \eqref{eq:gamma+}.  The behavior of
$\mathcal N_\e^+$ for small $t$ and $\e$ is described by the following
proposition.
\begin{Proposition}\label{p:lim_1}
There holds
$$
\lim_{t\to 0^+}\Big(\lim_{\e\to0^+}\mathcal N_\e^+(t)\Big)= \lim_{r\to
  1^+}\Big(\lim_{\e\to0^+}\mathcal N_\e(r)\Big)= 1.
$$
\end{Proposition}
\begin{pf}
  Let us first notice that the strong $\mathcal D^{1,2}(\R^N)$
  convergence of $u_\e$ to $\varphi_{k_0}^+$ ensured by Lemma
  \ref{l:conv_eigen} implies that, for all $t\in(0,3)$,
\begin{equation}\label{eq:56}
\lim_{\e\to 0^+}\mathcal N_\e^+(t)=
\mathcal N^+(t)
\end{equation}
where 
$$
\mathcal N^+(t)
=\frac{
 t\int_{B^+_{t}}\Big(|\nabla
  \varphi_{k_0}^+(x)|^2-\lambda_{k_0}(D^+) p(x) (\varphi_{k_0}^+(x))^2\Big)dx
}{\int_{\Gamma_t^+}(\varphi_{k_0}^+(x))^2\,d\sigma}.
$$
Let us define
$$
\varphi_0(x)=\varphi_0(x_1,x')=
\begin{cases}
\varphi_{k_0}^+(x_1+1,x'),&\text{if }x_1\geq0,\\
-\varphi_{k_0}^+(-x_1+1,x'),&\text{if }x_1<0,\\
\end{cases}
$$
and observe that $\varphi_0\in\mathcal D^{1,2}(\R^N)$ satisfies 
$\varphi_0(-x_1,x')=-\varphi_0(x_1,x')$ and weakly solves 
$$
-\Delta\varphi_0(x)=\lambda_{k_0}(D^+) p_0(x)\varphi_0(x),
$$
where 
$$
p_0(x)=p_0(x_1,x')=
\begin{cases}
p(x_1+1,x'),&\text{if }x_1\geq0,\\
p(-x_1+1,x'),&\text{if }x_1<0.\\
\end{cases}
$$
Moreover $\mathcal N^+$ can be rewritten as
$$
\mathcal N^+(t)
=\frac{
 t\int_{B({\mathbf 0},t)}\Big(|\nabla
  \varphi_0(x)|^2-\lambda_{k_0}(D^+) p_0(x) \varphi_{0}^2(x)\Big)dx
}{\int_{\partial B({\mathbf 0},t)}\varphi_0^2(x)\,d\sigma}.
$$
Hence, from \cite[Theorem 1.3]{FFT} it follows that there exist
$j_0\in\N$ and an eigenfunction $Y$ of $-\Delta_{{\mathbb S}^{N-1}}$
associated to the eigenvalue $j_0(N-2+j_0)$, i.e.  satisfying 
$-\Delta_{{\mathbb S}^{N-1}}Y=j_0(N-2+j_0) Y$ on ${\mathbb S}^{N-1}$,
such that 
\begin{equation}\label{eq:57}
\lim_{t\to0^+}\mathcal N^+(t)=-\frac{N-2}{2}
+\sqrt{\bigg(\frac{N-2}{2}\bigg)^{\!\!2}+j_0(N-2+j_0)}=j_0
\end{equation}
and 
\begin{align}
\label{eq:66}&\lambda^{-j_0}\varphi_0(\lambda\theta)
\to Y(\theta)\quad\text{in }C^{1,\tau}({\mathbb S}^{N-1}),\\
&\label{eq:67}\lambda^{1-j_0}\nabla \varphi_0(\lambda\theta)\to
j_0 Y(\theta)\theta+\nabla_{{\mathbb S}^{N-1}}
Y(\theta)\quad \text{in }C^{0,\tau}({\mathbb S}^{N-1}),
\end{align}
as $\lambda\to 0^+$,
for every $\tau\in(0,1)$.  Since the nodal set of $\varphi_0$ is
$\{0\}\times\R^{N-1}$, we infer that  $Y$ vanishes on  the
equator ${\mathbb S}^{N-1}\cap(\{0\}\times\R^{N-1})$. Therefore, $Y$ can not be the first eigenfunction 
of $-\Delta_{{\mathbb S}^{N-1}}$ and hence 
 $j_0\geq 1$ necessarily. On the other hand, (\ref{eq:13}) and (\ref{eq:67}) 
imply that $j_0\leq 1$. Hence $j_0=1$. 
 The conclusion hence follows from \eqref{eq:56} and \eqref{eq:57}.
\end{pf}

\begin{Lemma}\label{l:poho_right}
For all $\e\in(0,\check \e)$ and $t\in(2\e,3)$ there holds
\begin{multline*}
t\int_{\Gamma_t^+}|\nabla
  u_\e|^2d\sigma=2\e\int_{\Gamma_{2\e}^+}
\bigg(
|\nabla
  u_\e|^2-2\Big|\frac{\partial u_\e}{\partial
    \nu}\Big|^2\bigg)
  +
  (N-2)\int_{B_t^+\setminus B_{2\e}^+}|\nabla
  u_\e(x)|^2dx+
2t\int_{\Gamma_t^+}\bigg|\frac{\partial u_\e}{\partial \nu}\bigg|^2
  d\sigma.
\end{multline*}
\end{Lemma}
\begin{pf}
The stated identity follows from multiplication of 
 equation (\ref{eq:24}) by $(x-{\mathbf e}_1)\cdot\nabla
  u_\e$ and integration by parts over $B_t^+\setminus B_{2\e}^+$.
\end{pf}

\begin{Lemma}\label{l:der_N+}
  For all $\e\in(0,\check \e)$, $\mathcal N_\e^+\in C^1(2\e,3)$ and 
$$
(\mathcal N_\e^+)'(t)=
 \frac{ 2t\Big[
  \big(\int_{\Gamma_t^+}\big|\frac{\partial u_\e}{\partial \nu}\big|^2
  d\sigma\big) \big(\int_{\Gamma_t^+}u_\e^2 d\sigma\big)
  -\big(\int_{\Gamma_t^+}u_\e\frac{\partial u_\e}{\partial \nu}
  d\sigma\big)^2\Big]} {\big(\int_{\Gamma_t^+}u_\e^2 d\sigma\big)^2}
+\frac{R_\e^+}{\int_{\Gamma_t^+}u_\e^2 d\sigma},
$$
for all $t\in (2\e,3)$, where
\begin{equation}\label{eq:R_eps}
R_\e^+=
\int_{\Gamma_{2\e}^+}
\bigg(-(N-2)u_\e\frac{\partial u_\e}{\partial\nu}+2\e|\nabla u_\e|^2-4\e
\Big|\frac{\partial u_\e}{\partial \nu}\Big|^2\bigg)d\sigma.
\end{equation}
\end{Lemma}
\begin{pf}
  Multiplication of equation (\ref{eq:24}) by $u_\e$ and integration
  by parts over $D^-\cup\mathcal C_\e\cup B^+_{t}$ yield, for every $t>\e$,
\begin{equation}\label{eq:43}
\int_{D^-\cup\mathcal
    C_\e\cup B^+_{t}}\Big(|\nabla u_\e(x)|^2-\lambda^\e_{\bar k} p
  u_\e^2(x)\Big)dx =
  \int_{\Gamma_{t}^+}\frac{\partial
    u_\e}{\partial \nu}u_\e\,d\sigma.
\end{equation}
From Lemma \ref{l:poho_right} and (\ref{eq:43}) we deduce
\begin{align}\label{eq:48}
  &(D_\e^+)'(t)= \frac{d}{dt}\bigg(\frac1{t^{N-2}}\! \int_{D^-\cup\mathcal
    C_\e\cup B^+_{2\e}}\!\!\Big(|\nabla u_\e(x)|^2-\lambda^\e_{\bar k} p
  u_\e^2(x)\Big)dx +\frac1{t^{N-2}} \int_{B^+_{t}\setminus
    B^+_{2\e}}|\nabla
  u_\e(x)|^2dx\bigg)\\
  \notag&=-\frac{N-2}{t^{N-1}}\int_{D^-\cup\mathcal C_\e\cup
    B^+_{2\e}}\Big(|\nabla u_\e(x)|^2-\lambda^\e_{\bar k} p
  u_\e^2(x)\Big)dx\\
\notag&\qquad- \frac{N-2}{t^{N-1}}\int_{B^+_{t}\setminus
    B^+_{2\e}}|\nabla u_\e(x)|^2dx
  +\frac1{t^{N-2}}\int_{\Gamma_t^+}|\nabla
  u_\e|^2d\sigma\\
  \notag&=-\frac{N-2}{t^{N-1}}
  \int_{\Gamma_{2\e}^+}\frac{\partial
    u_\e}{\partial \nu}u_\e\,d\sigma+ \frac{2\e}{t^{N-1}}
  \int_{\Gamma_{2\e}^+}\bigg(|\nabla u_\e|^2- 2\Big|\frac{\partial
    u_\e}{\partial \nu}\Big|^2\bigg)d\sigma
  +\frac2{t^{N-2}}\int_{\Gamma_t^+}\Big|\frac{\partial u_\e}{\partial
    \nu}\Big|^2d\sigma\\
\notag&=\frac2{t^{N-2}}\int_{\Gamma_t^+}\Big|\frac{\partial u_\e}{\partial
    \nu}\Big|^2d\sigma+\frac{R_\e^+}{t^{N-1}}
\end{align}
for all $t\in(2\e,3)$. Furthermore
\begin{equation}\label{eq:45}
(H_\e^+)'(t)=\frac2{t^{N-1}}  \int_{\Gamma_{t}^+}\frac{\partial
    u_\e}{\partial \nu}u_\e\,d\sigma
\end{equation}
which, in view of (\ref{eq:43}), implies 
\begin{equation}\label{eq:46}
(H_\e^+)'(t)=\frac2t \, D_\e^+(t)\quad\text{for all }t\in(\e,3).
\end{equation}
From (\ref{eq:47}) and (\ref{eq:46}) it follows that 
$$
(\mathcal N_\e^+)'(t)=
\frac{(D_\e^+)'(t)H_\e^+(t)-\frac t2\big((H_\e^+)'(t)\big)^2}{(H_\e^+(t))^2}
$$
which yields the conclusion in view of (\ref{eq:48}) and (\ref{eq:45}).
\end{pf}

\begin{Lemma}\label{l:est_R_eps}
  For $\e\in(0,\check\e)$, let $R_\e^+$ as in (\ref{eq:R_eps}). There
  exists $C_8>0$ such that, for all $\e\in(0,\check\e)$,
$$ 
|R_\e^+|\leq C_8 \e^N.
$$
\end{Lemma}
\begin{pf}
  From (\ref{eq:R_eps}), Lemmas \ref{l:sup_sol_bound} and
  \ref{l:grad_est}, and (\ref{eq:26}), it follows that, for all
  $\e\in(0,\check\e)$,
$$
|R_\e^+|\leq {\rm const\,}\int_{\Gamma_{2\e}^+}
(\e+\Phi^\e)\,d\sigma={\rm const\,}
\bigg(2^{N-2}\e^N\omega_{N-1}+\e^N\int_{\Gamma_{2}^+}\Phi_1d\sigma+2^N\gamma_\e\e^N
\int_{\Gamma_{1}^+}\Phi_2d\sigma\bigg)
$$
thus implying the conclusion.
\end{pf}

As a consequence of the above estimates, we finally obtain the
following uniform control of the frequency close to the right junction
of the tube.
\begin{Lemma}\label{l:C9}
There exists  $C_9>0$ such that, for all $\e\in(0,\min\{\e_2, \check\e\})$ and 
$t\in (2\e,r_0)$,
\begin{equation}\label{eq:55}
(\mathcal N_\e^+)'(t)\geq -C_9\frac{\e^N}{t^{N+1}}
\end{equation}
\end{Lemma}
\begin{pf}
  From Lemma \ref{l:usot}, we deduce that, for all $t\in (2\e,r_0)$
  and $\e\in(0,\min\{\e_2, \check\e\})$,
\begin{equation}\label{eq:49}
  \int_{\Gamma_t^+}u_\e^2 d\sigma\geq \frac{C_5^2}{4}
  \int_{\Gamma_t^+}(x_1-1)^2d\sigma=
  \frac{C_5^2}{8}t^{N+1}
  \int_{\SN}|\theta\cdot {\mathbf e}_1|^2d\sigma(\theta).
\end{equation}
The conclusion follows from Lemma \ref{l:der_N+}, Schwarz's inequality, 
Lemma \ref{l:est_R_eps}, and (\ref{eq:49}).
\end{pf}

\begin{Corollary}\label{cor:stima_N+}
For all $\e\in(0,\min\{\e_2, \check\e\})$ and $r_1,r_2$ such that $1+2\e<r_1<r_2<1+r_0$
there holds 
$$
\mathcal N_\e(r_1)\leq \mathcal
N_\e(r_2)+\frac{C_9}{N}\Big(\frac{\e}{r_1-1}\Big)^{\!\!N} \leq \mathcal
N_\e(r_2)+\frac{C_9}{N 2^N}.
$$
\end{Corollary}
\begin{pf}
It follows from (\ref{eq:N_N+}) and integration of (\ref{eq:55}).
\end{pf}

\begin{Corollary}\label{cor:N+sim1}
For every $\delta>0$ there exist $\tilde r_\delta,\widetilde R_\delta>0$ such that 
$$
\mathcal N_\e(1+R\e)\leq 1+\delta
\quad\text{for all }R>\widetilde R_\delta\text{ and }
\e\in \bigg(0,\frac{\tilde r_\delta}{R}\bigg).
$$
\end{Corollary}
\begin{pf}
Let $\delta>0$. From Proposition \ref{p:lim_1} there exist $\tilde r_\delta\in(0,r_0)$ 
and $\tilde \e_\delta>0$ such that 
\begin{equation}\label{eq:100}
\mathcal N_\e(1+\tilde r_\delta)\leq
1+\frac\delta2\quad\text{for all }\e\in (0,\tilde \e_\delta).
\end{equation}
Let $\widetilde R_\delta>\max\{2,\tilde r_\delta/\min\{\e_2, \check\e\}\}$ be such that $\frac{C_9}{N}\widetilde R_\delta^{-N}<\frac\delta2$.
Then, from Corollary \ref{cor:stima_N+}, for all $R>\widetilde R_\delta$ and 
$\e\in \big(0,\frac{\tilde r_\delta}{R}\big)$ there holds
\begin{equation}\label{eq:101}
\mathcal N_\e(1+R\e)\leq \mathcal N_\e(1+\tilde r_\delta)+\frac{C_9}{N}R^{-N}
\leq \mathcal N_\e(1+\tilde r_\delta)+\frac{C_9}{N}\widetilde R_\delta^{-N}
\leq \mathcal N_\e(1+\tilde r_\delta)+\frac{\delta}2.
\end{equation}
The conclusion follows from (\ref{eq:100}) and (\ref{eq:101}).
\end{pf}

\subsection{The frequency function at the left}\label{sec:freq-funct-at-1}

If $\e\in(0,\check\e)$ and $r\in(-\check R,-\e)$, then
\begin{equation*}
\mathcal N_\e(r)=\mathcal N_\e^-(-r)
\end{equation*}
where, for $t\in(\e,\check R)$, 
\begin{align}
  \label{eq:155} \mathcal N_\e^-(t)&=
  \frac{D_\e^-(t)}{H_\e^-(t)},\\
  \label{eq:156} D_\e^-(t)&=\frac1{t^{N-2}} \int_{\Omega_{-t}}\Big(|\nabla
  u_\e(x)|^2-\lambda^\e_{\bar k} p(x) u_\e^2(x)\Big)dx,\\
  \label{eq:157}
  H_\e^-(t)&=\frac1{t^{N-1}}\int_{\Gamma_t^-}u_\e^2(x)\,d\sigma,
\end{align}
with $\Gamma_t^-$ defined in $(\ref{eq:defGamma_r-})$.

\begin{Lemma}\label{l:poho_left}
For $t>\e$ there holds
\begin{multline*}
  t\int_{\Gamma_t^-}\Big(|\nabla
  u_\e|^2-\lambda^\e_{\bar k} p u_\e^2\Big)d\sigma=
2t\int_{\Gamma_t^-}\bigg|\frac{\partial u_\e}{\partial \nu}\bigg|^2
  d\sigma-(N-2)\int_{\Omega_{-t}}|\nabla
  u_\e(x)|^2dx\\+\lambda^\e_{\bar k}\int_{\Omega_{-t}}(Np(x)+x\cdot \nabla p(x))
u_\e^2(x)dx,
\end{multline*}
where $\nu=\nu(x)=\frac{x}{|x|}$.
\end{Lemma}
\begin{pf}
  The stated identity follows from multiplication of equation
  (\ref{eq:24}) by $x\cdot\nabla u_\e$ and integration by parts over
  $\Omega_{-t}$.
\end{pf}
\begin{Lemma}\label{l:DH-'}
For $\e\in(0,\check \e)$ and $t\in(\e,\check R)$ there holds
\begin{align}
\label{eq:58}\frac{d}{dt}D_\e^-(t)=&-\frac2{t^{N-2}}
\int_{\Gamma_t^-}\bigg|\frac{\partial u_\e}{\partial \nu}\bigg|^2
  d\sigma-\frac{\lambda^\e_{\bar k}}{t^{N-1}}\int_{\Omega_{-t}}(2p(x)+x\cdot \nabla p(x))
u_\e^2(x)dx,\\
\label{eq:59}\frac{d}{dt}H_\e^-(t)=&
  \frac{2}{t^{N-1}} \int_{\Gamma_t^-} u_\e\,\frac{\partial u_\e}{\partial \nu}
  \,d\sigma=-\frac{2}{t}D_\e^-(t),\\
\label{eq:70}\frac{d}{dt}\mathcal N_\e^-(t)=&
 -2t\frac{\left(\int_{\Gamma_t^-}\big|\frac{\partial u_\e}{\partial \nu}\big|^2
  d\sigma\right)\left(\int_{\Gamma_t^-}u_\e^2(x)\,d\sigma\right)
-\left(\int_{\Gamma_t^-}u_\e\frac{\partial u_\e}{\partial \nu}
  d\sigma\right)^2}
{\left(\int_{\Gamma_t^-}u_\e^2(x)\,d\sigma\right)^2}\\[3pt]
&\notag-\lambda^\e_{\bar k}\frac{\int_{\Omega_{-t}}(2p(x)+x\cdot \nabla p(x))
u_\e^2(x)dx}{\int_{\Gamma_t^-}u_\e^2(x)  d\sigma}.
\end{align}
\end{Lemma}
\begin{pf}
Since 
\begin{equation*}
\frac{d}{dt}D_\e^-(t)=
-\frac{N-2}{t^{N-1}} \int_{\Omega_{-t}}\Big(|\nabla
  u_\e(x)|^2-\lambda^\e_{\bar k} p u_\e^2(x)\Big)dx
-\frac1{t^{N-2}} \int_{\Gamma_t^-}\Big(|\nabla
  u_\e|^2-\lambda^\e_{\bar k} p u_\e^2\Big)d\sigma,
\end{equation*}
(\ref{eq:58}) follows from Lemma \ref{l:poho_left}. From direct calculation, we obtain that 
\begin{equation*}
\frac{d}{dt}H_\e^-(t)=
  \frac{2}{t^{N-1}} \int_{\Gamma_t^-} u_\e\,\frac{\partial u_\e}{\partial \nu}
  \,d\sigma,
\end{equation*}
while testing equation (\ref{eq:24}) with $u_\e$ and integration over $\Omega_{-t}$
yield
$$
\int_{\Omega_{-t}}\Big(|\nabla
u_\e(x)|^2-\lambda^\e_{\bar k} p u_\e^2(x)\Big)dx=-\int_{\Gamma_t^-}
u_\e\,\frac{\partial u_\e}{\partial \nu} \,d\sigma,
$$
thus implying (\ref{eq:59}). Finally, (\ref{eq:70}) follows from
(\ref{eq:58}), (\ref{eq:59}), and $(\mathcal N_\e^-)'
=\frac{(D_\e^-)'H_\e^--D_\e^-(H_\e^-)'}{(H_\e^-)^2}$.
\end{pf}

\noindent The following estimates strongly rely on Lemmas \ref{l:poincareN-1}
and \ref{l:lemmaZERO}. 
\begin{Lemma}\label{l:H-lower}
For every $\delta\in(0,1)$ there exist $\bar r_\delta\in(0,\check R)$
and $\bar\e_\delta\in(0,\check \e)$ such that, 
for every $\e\in (0,\bar\e_\delta)$, 
\begin{align}\label{eq:64}
  \frac{\frac{d}{dt}H_\e^-(t)}{H_\e^-(t)}&\leq
  -\frac{2(1-\delta)(N-1)}{t}
  \quad\text{for all } t\in(\e,\bar r_\delta),\\
  \label{eq:69}\frac{\frac{d}{dt}D_\e^-(t)}{D_\e^-(t)}&\leq
  -\frac{2(1-\delta)(N-1)}{t}\quad\text{for all } t\in(\e,\bar
  r_\delta),
\end{align}
\begin{align}\label{eq:65}
  &H_\e^-(t_1)\geq \bigg(\frac{t_2}{t_1}\bigg)^{\!\!2(1-\delta)(N-1)}
  H_\e^-(t_2)\quad \text{for all }t_1,t_2\in(\e,\bar r_\delta) \text{ such that }t_1<t_2,\\
\label{eq:71}
&D_\e^-(t_1)\geq \bigg(\frac{t_2}{t_1}\bigg)^{\!\!2(1-\delta)(N-1)}
D_\e^-(t_2)\quad \text{for all }t_1,t_2\in(\e,\bar r_\delta) \text{
  such that }t_1<t_2.
\end{align}
\end{Lemma}
\begin{pf}
From Lemmas \ref{l:DH-'}, \ref{l:lemmaZERO}, and \ref{l:poincareN-1}, we deduce that, 
for every $\delta\in(0,1)$, there exist $\bar r_\delta>0$ and $\bar\e_\delta>0$ such that, 
for every $\e\in (0,\bar\e_\delta)$ and $t\in(\e,\bar r_\delta)$,  there holds
\begin{align*}
\frac{d}{dt}H_\e^-(t)&=-\frac2{t^{N-1}}
 \int_{\Omega_{-t}}\Big(|\nabla
  u_\e(x)|^2-\lambda^\e_{\bar k} p(x) u_\e^2(x)\Big)dx\\
&\leq -\frac{2(1-\delta)}{t^{N-1}}\int_{\Omega_{-t}}|\nabla
  u_\e(x)|^2dx
\leq -\frac{2(1-\delta)(N-1)}{t}H_\e^-(t)
\end{align*}
which yields (\ref{eq:64}).  From (\ref{eq:59}), we have that
$$
\int_{\Omega_{-t}}\Big(|\nabla
u_\e(x)|^2-\lambda^\e_{\bar k} p(x) u_\e^2(x)\Big)dx=-
\int_{\Gamma_t^-} u_\e\,\frac{\partial u_\e}{\partial \nu} \,d\sigma
$$
which, by Schwarz's inequality, Lemmas \ref{l:lemmaZERO} and \ref{l:poincareN-1}, 
up to shrinking   $\bar r_\delta>0$ and $\bar\e_\delta>0$,
for every $\e\in (0,\bar\e_\delta)$ and $t\in(\e,\bar r_\delta)$ yields 
\begin{align}\label{eq:68}
  \int_{\Gamma_t^-} \bigg|\frac{\partial u_\e}{\partial \nu}\bigg|^2
  \,d\sigma &\geq \frac{\int_{\Omega_{-t}}\Big(|\nabla
    u_\e(x)|^2-\lambda^\e_{\bar k} p(x)
    u_\e^2(x)\Big)dx}{\int_{\Gamma_t^-} u_\e^2 \,d\sigma }
 \int_{\Omega_{-t}}\Big(|\nabla u_\e(x)|^2-\lambda^\e_{\bar k}
  p(x) u_\e^2(x)\Big)dx\\
\notag&\geq \frac{1-\frac\delta2}{t}\frac{\frac{1}{t^{N-2}}\int_{\Omega_{-t}}|\nabla
    u_\e(x)|^2dx}{\frac1{t^{N-1}}\int_{\Gamma_t^-} u_\e^2 \,d\sigma }
 \int_{\Omega_{-t}}\Big(|\nabla u_\e(x)|^2-\lambda^\e_{\bar k}
  p(x) u_\e^2(x)\Big)dx\\
\notag&\geq\frac{(1-\frac{\delta}{2})(N-1)}{t}\int_{\Omega_{-t}}\Big(|\nabla u_\e(x)|^2-\lambda^\e_{\bar k}
  p(x) u_\e^2(x)\Big)dx.
\end{align}
From (\ref{eq:58}), (\ref{eq:68}), (\ref{eq:p}), and Lemma \ref{l:lemmaZERO},
up to shrinking   $\bar r_\delta>0$ and $\bar\e_\delta>0$, there holds
\begin{align*}
-\frac{d}{dt}&D_\e^-(t)=\frac2{t^{N-2}}
\int_{\Gamma_t^-}\bigg|\frac{\partial u_\e}{\partial \nu}\bigg|^2
  d\sigma+\frac{\lambda^\e_{\bar k}}{t^{N-1}}\int_{\Omega_{-t}}(2p(x)+x\cdot \nabla p(x))
u_\e^2(x)dx,\\
&\geq \frac{2(1-\frac{\delta}{2})(N-1)}{t^{N-1}}\int_{\Omega_{-t}}\Big(|\nabla u_\e(x)|^2-\lambda^\e_{\bar k}
  p(x) u_\e^2(x)\Big)dx\\
&\qquad-
\frac{\delta(N-1)}{t^{N-1}}\int_{\Omega_{-t}}\Big(|\nabla u_\e(x)|^2-\lambda^\e_{\bar k}
  p(x) u_\e^2(x)\Big)dx
\\
&\geq \frac{2(1-\delta)(N-1)}{t^{N-1}}\int_{\Omega_{-t}}\Big(|\nabla u_\e(x)|^2-\lambda^\e_{\bar k}
  p(x) u_\e^2(x)\Big)dx=\frac{2(1-\delta)(N-1)}{t}D_\e^-(t)
\end{align*}
thus proving (\ref{eq:69}).

Estimate (\ref{eq:65}) follows by integration of (\ref{eq:64}), while
(\ref{eq:71}) follows by integration of (\ref{eq:69}).
\end{pf}

\begin{Lemma}\label{l:stima N+}
  For every $\delta>0$ there exist $\check R_\delta\in(0,\check R)$, and
  $\check\e_\delta\in(0,\check \e)$ such that 
\begin{align}\label{eq:72}
  \frac{\frac{d}{dt}\mathcal N_\e^-(t)}{\mathcal N_\e^-(t)}\leq
\delta\quad\text{and}\quad \int_{\Omega_{-t}} \Big(|\nabla
u_\e|^2-\lambda^\e_{\bar k} p u_\e^2\Big)dx \geq
\frac12\int_{\Omega_{-t}} |\nabla u_\e|^2dx
\end{align}
for every $\e\in
  (0,\check\e_\delta)$ and  $t\in(\e,\check R_\delta)$.
\end{Lemma}
\begin{pf}
From Lemma \ref{l:H-lower}, letting $\delta_0=\frac{2N-5}{4(N-1)}\in (0,1)$, 
there holds 
\begin{align}
  \label{eq:74}
  D_\e^-(t_1)\geq \bigg(\frac{t_2}{t_1}\bigg)^{\!\!N+\frac12}
  D_\e^-(t_2)\quad \text{for every $\e\in (0,\bar\e_{\delta_0})$ and
  }t_1,t_2\in(\e,\bar r_{\delta_0}) \text{ such that }t_1<t_2.
\end{align}
Let us fix $\delta>0$. From (\ref{eq:p}), Lemma \ref{l:lemmaZERO},
\eqref{eq:coercivity}, and \eqref{eq:52}, there exist $\breve
R_\delta\in(0,\min\{\bar r_{\delta_0},\check R\})$ and $\check\e_\delta\in(0,
\min\{\bar\e_{\delta_0},\check \e\})$ such that
\begin{align}\label{eq:150}
&\|2p+x\cdot\nabla p\|_{L^{3N}\big(B^-_{\breve R_\delta}\big)}
\leq 
\bigg(\frac {2N}{\omega_{N-1}}\bigg)^{\!\!\frac5{3N}}
\frac{C_S\delta}{8\lambda_{k_0}(D^+)},\\
\label{eq:154}&\lambda^\e_{\bar k}\leq 2\lambda_{k_0}(D^+)\quad\text{for all }\e\in (0,\check\e_\delta),\\
\label{eq:76}& \int_{\Omega_{-t}} \Big(|\nabla
u_\e|^2-\lambda^\e_{\bar k} p u_\e^2\Big)dx \geq
\frac12\int_{\Omega_{-t}} |\nabla u_\e|^2dx,\text{ for all }\e\in
(0,\check\e_\delta),\ t\in(\e,\breve R_\delta)
\\
  \label{eq:75}& \int_{\Omega_{-t}} \!\!\Big(|\nabla
  u_\e|^2-\lambda^\e_{\bar k} p u_\e^2\Big)dx \geq\tfrac{4\lambda_{k_0}(D^+)}\delta \!
  \int_{\Omega_{-t}} \!\!|2p+x\cdot \nabla p| u_\e^2dx, \text{ for all }\e\in
(0,\check\e_\delta),\ t\in(\e,\breve R_\delta).
\end{align}
Let $\check R_\delta=\breve R_\delta^{5/3}$. From (\ref{eq:70}), \eqref{eq:154}, and
Schwarz's inequality, we have that, for all $\e\in
(0,\check\e_\delta)$ and $t\in(\e,\check R_\delta)$
 \begin{align}\label{eq:73}
  \frac{\frac{d}{dt}\mathcal N_\e^-(t)}{\mathcal N_\e^-(t)}&\leq
\mathcal I_\e(t)
\end{align}
where 
\begin{align}\label{eq:I(t)}
\mathcal I_\e(t)
:=\frac{2\lambda_{k_0}(D^+) }{t}\frac{\int_{\Omega_{-t}}|2p(x)+x\cdot \nabla p(x)|
u_\e^2(x)dx}{\int_{\Omega_{-t}} \Big(|\nabla
  u_\e(x)|^2-\lambda^\e_{\bar k} p(x) u_\e^2(x)\Big)dx}
=\frac{2\lambda_{k_0}(D^+)}{t}\Big(I_\e(t)+I\!I_\e(t)\Big)
\end{align}
with
\begin{align*}
  &I_\e(t)=\frac{\int_{\Omega_{-t}\setminus\Omega_{-t^{3/5}}}|2p(x)+x\cdot
    \nabla p(x)| u_\e^2(x)dx}{\int_{\Omega_{-t}} \Big(|\nabla
    u_\e(x)|^2-\lambda^\e_{\bar k} p(x) u_\e^2(x)\Big)dx},\\
  &I\!I_\e(t)= \frac{\int_{\Omega_{-t^{3/5}}}|2p(x)+x\cdot \nabla p(x)|
    u_\e^2(x)dx}{\int_{\Omega_{-t}} \Big(|\nabla
    u_\e(x)|^2-\lambda^\e_{\bar k} p(x) u_\e^2(x)\Big)dx}.
\end{align*}
By H\"older inequality, (\ref{eq:76}), Lemma \ref{l:sob_ex}, and
\eqref{eq:150}, $I_\e(t)$ can be estimated as
\begin{align}\label{eq:77}
  I_\e(t)&\leq \|2p+x\cdot\nabla p\|_{L^{3N}\big(B^-_{t^{3/5}}\big)}
  \Big|\Omega_{-t}\setminus\Omega_{-t^{3/5}}\Big|^{\frac5{3N}}
  \frac{\Big(\int_{\Omega_{-t}}|u_\e(x)|^{2^*}dx\Big)^{2/2^*}}
  {\int_{\Omega_{-t}} \Big(|\nabla
    u_\e(x)|^2-\lambda^\e_{\bar k} p(x) u_\e^2(x)\Big)dx}\\
  &\notag\leq
  \frac{2}{C_S}\bigg(\frac{\omega_{N-1}}{2N}\bigg)^{\!\!\frac5{3N}}
  \|2p+x\cdot\nabla p\|_{L^{3N}\big(B^-_{\breve R_\delta}\big)}t \leq
  \frac{\delta}{4\lambda_{k_0}(D^+)}\,t
\end{align}
 for all $t\in(\e,\check R_\delta)$ and $\e\in(0, \check\e_\delta)$.
On the other hand, from (\ref{eq:75}) and (\ref{eq:74})
\begin{align}\label{eq:78}
  I\!I_\e(t)&= \frac{\int_{\Omega_{-t^{3/5}}}|2p(x)+x\cdot \nabla p(x)|
    u_\e^2(x)dx}{\int_{\Omega_{-t^{3/5}}} \Big(|\nabla
    u_\e(x)|^2-\lambda^\e_{\bar k} p(x) u_\e^2(x)\Big)dx}
\frac{\int_{\Omega_{-t^{3/5}}} \Big(|\nabla
    u_\e(x)|^2-\lambda^\e_{\bar k} p(x) u_\e^2(x)\Big)dx}{\int_{\Omega_{-t}} \Big(|\nabla
    u_\e(x)|^2-\lambda^\e_{\bar k} p(x) u_\e^2(x)\Big)dx}\\
&\notag\leq \frac\delta{4\lambda_{k_0}(D^+)}\, t^{-\frac25(N-2)}\frac{ D_\e^-(t^{3/5})}{
  D_\e^-(t)}\leq
\frac\delta{4\lambda_{k_0}(D^+)}\, t^{-\frac25(N-2)}
 \Big(\frac{t}{t^{3/5}}\Big)^{\!N+\frac12}=\frac\delta{4\lambda_{k_0}(D^+)}\,t
\end{align}
for all $t\in(\e,\check R_\delta)$ and $\e\in(0, \check\e_\delta)$. 
 \eqref{eq:I(t)}, (\ref{eq:77}), and (\ref{eq:78}) imply that 
 \begin{equation}
   \label{eq:163}
   \mathcal I_\e(t)\leq \delta
 \end{equation}
for all $t\in(\e,\check R_\delta)$ and $\e\in(0, \check\e_\delta)$. 
Estimate (\ref{eq:72}) follows from
\eqref{eq:163} and (\ref{eq:73}).
\end{pf}

\begin{Corollary}\label{cor:stima_N-}
  For every $\delta>0$, let $\check R_\delta\in(0,1)$ and
  $\check\e_\delta>0$ as in Lemma \ref{l:stima N+}. Then, for every
  $\e\in(0,\check\e_\delta)$ and $r_1,r_2$ such that $-\check
  R_\delta<r_1<r_2<-\e$, there holds
$$
\mathcal N_\e(r_1)\leq \mathcal N_\e(r_2) e^{\delta(r_2-r_1)}.
$$
\end{Corollary}
\begin{pf}
It follows from integration of (\ref{eq:72}).
\end{pf}

\subsection{The frequency function in the corridor}\label{sec:freq-funct-corr}

If $\e\in(0,\check\e)$ and $0< r<1$, then
\begin{equation}\label{eq:N_c}
\mathcal N_\e(r)=
\frac{\e D_\e^c(r)}{H_\e^c(r)}
\end{equation}
where 
\begin{align*}
 D_\e^c(r)=\int_{\Omega_{r}^\e}\Big(|\nabla
  u_\e(x)|^2-\lambda^\e_{\bar k} p(x) u_\e^2(x)\Big)dx,
\quad H_\e^c(r)=\int_{\Gamma_r^\e}u_\e^2(x)\,d\sigma.
\end{align*}

\begin{Lemma}\label{l:poho_corr}
For all $\e\in(0,\e_0)$ and $r\in (0,1)$
\begin{align*}
  \int_{\Gamma_r^\e}\Big(|\nabla
  u_\e|^2-\lambda^\e_{\bar k} p u_\e^2\Big)d\sigma=
2\int_{\Gamma_r^\e}\bigg|\frac{\partial u_\e}{\partial x_1}\bigg|^2
  d\sigma
+\int_{S_\e}\bigg|\frac{\partial u_\e}{\partial x_1}\bigg|^2
  d\sigma
-\lambda^\e_{\bar k}\int_{\Omega_{r}^\e}\frac{\partial p}{\partial x_1}(x)
u_\e^2(x)dx,
\end{align*}
where $S_\e=\partial D^-\setminus\Gamma_0^\e=\big\{(0,x')\in\R\times\R^{N-1}:\frac{x'}\e\not\in\Sigma\big\}$.
\end{Lemma}
\begin{pf}
The stated identity follows from multiplication of 
 equation (\ref{eq:24}) by $\frac{\partial u_\e}{\partial x_1}$
 and integration by parts over $\Omega_{r}^\e$.
\end{pf}

\begin{Lemma}\label{l:DHN'corr}
For all $\e\in(0,\check \e)$ and $r\in (0,1)$ there holds
\begin{align}
  \label{eq:79}
  \frac{d}{dr}D_\e^c(r)=\,&
2\int_{\Gamma_r^\e}\bigg|\frac{\partial
    u_\e}{\partial x_1}\bigg|^2 d\sigma
  +\int_{S_\e}\bigg|\frac{\partial u_\e}{\partial x_1}\bigg|^2 d\sigma
  -\lambda^\e_{\bar k}\int_{\Omega_{r}^\e}\frac{\partial p}{\partial
    x_1}(x)
  u_\e^2(x)dx,\\
  \label{eq:81}\frac{d}{dr}H_\e^c(r)=\,&2\int_{\Gamma_r^\e}u_\e\frac{\partial
    u_\e}{\partial x_1}\,d\sigma
  =2D_\e^c(r),\\
\label{eq:82}\frac{d}{dr}\mathcal N_\e(r)=\,&
\e\!\left[
  2\frac{\Big(\int_{\Gamma_r^\e}\big|\frac{\partial u_\e}{\partial
      x_1}\big|^2
    d\sigma\Big)\Big(\int_{\Gamma_r^\e}u_\e^2\,d\sigma\Big)-\Big(\int_{\Gamma_r^\e}u_\e\frac{\partial
      u_\e}{\partial
      x_1}\,d\sigma\Big)^2}{\Big(\int_{\Gamma_r^\e}u_\e^2\,d\sigma\Big)^2}
  +\frac{\int_{S_\e}\big|\frac{\partial u_\e}{\partial x_1}\big|^2
    d\sigma}{\int_{\Gamma_r^\e}u_\e^2\,d\sigma}\right]\\
&\notag - \e\lambda^\e_{\bar
  k}\frac{\int_{\Omega_{r}^\e}\frac{\partial p}{\partial x_1}(x)
  u_\e^2(x)dx}{\int_{\Gamma_r^\e}u_\e^2\,d\sigma}.
\end{align}
\end{Lemma}
\begin{pf}
Since 
\begin{equation*}
  \frac{d}{dr}D_\e^c(r)=  \int_{\Gamma_r^\e}\Big(|\nabla
  u_\e|^2-\lambda^\e_{\bar k} p u_\e^2\Big)d\sigma,
\end{equation*}
(\ref{eq:79}) follows from Lemma \ref{l:poho_corr}. From direct calculation, we obtain that 
\begin{equation*}
\frac{d}{dr}H_\e^c(r)=2\int_{\Gamma_r^\e}u_\e\frac{\partial
    u_\e}{\partial x_1}\,d\sigma,
\end{equation*}
while, testing equation (\ref{eq:24}) with $u_\e$ and integrating over
$\Omega_r^\e$, we have that
$$
\int_{\Omega_{r}^\e}\Big(|\nabla u_\e(x)|^2-\lambda^\e_{\bar k} p
u_\e^2(x)\Big)dx= \int_{\Gamma_r^\e}u_\e\frac{\partial u_\e}{\partial
  x_1}\,d\sigma,
$$
thus implying (\ref{eq:81}). Finally, $(\mathcal N_\e)' =
\e\frac{(D_\e^c)'H_\e^c-D_\e^c(H_\e^c)'}{(H_\e^c)^2}$, (\ref{eq:79}),
and~(\ref{eq:81}) yield (\ref{eq:82}).~\end{pf}

\begin{Lemma}\label{l:stima N_c}
For every $\delta>0$ there exists $\bar\e_c^\delta\in(0,\check\e)$ such that 
\begin{align}\label{eq:72c}
  \frac{\frac{d}{dt}\mathcal N_\e(r)}{\mathcal N_\e(r)}\geq
-\delta\quad&\text{for all }\e\in (0,\bar\e_c^\delta)
\text{ and  } r\in(0,1),\\
\label{eq:118}\mathcal N_\e(r_1)\leq \mathcal N_\e(r_2) e^{\delta(r_2-r_1)}
\quad&\text{for all }\e\in (0,\bar\e_c^\delta)
\text{ and  } 
0<r_1<r_2<1.
\end{align}
\end{Lemma}
\begin{pf}
From (\ref{eq:82}) and  Schwarz's inequality we have that, for all
$\e\in(0,\check \e)$ and $r\in (0,1)$,  
\begin{align}\label{eq:83}
  \frac{d}{dr}\mathcal N_\e(r)&\geq- \e\lambda^\e_{\bar
    k}\frac{\int_{\Omega_{r}^\e}\frac{\partial p}{\partial x_1}(x)
    u_\e^2(x)dx}{\int_{\Gamma_r^\e}u_\e^2\,d\sigma}.
\end{align}
By part ii) of Lemma \ref{l:lemmaZERO}, for every $\delta>0$ there
exists $\bar\e_c^\delta\in(0,\check\e)$ such that, for every $\e\in
(0,\bar\e_c^\delta)$ and $r\in (0,1)$,
\begin{align}\label{eq:86}
 \int_{\Omega_{r}^\e} \Big(|\nabla
  u_\e(x)|^2-\lambda^\e_{\bar k} p(x) u_\e^2(x)\Big)dx \geq\frac{\lambda^\e_{\bar
    k}}{\delta}
  \int_{\Omega_{r}^\e} \bigg|\frac{\partial p}{\partial x_1}(x)\bigg| u_\e^2(x)dx.
\end{align}
Estimate (\ref{eq:72c}) follows from (\ref{eq:83}), (\ref{eq:86}), and (\ref{eq:N_c}).
(\ref{eq:118})  follows from integration of (\ref{eq:72c}).
\end{pf}

\section{Blow-up at the right}\label{sec:limiting-problem}

Throughout this section, $\widetilde u_\e$ will denote the scaling of $u_\e$ 
 introduced  in (\ref{eq:84}--\ref{eq:85}). 
For
every $R>1$ we define as 
$\mathcal H^+_R$  the completion of 
$$
\mathcal D_R^+:=\Big\{v\in C^\infty(\overline{((-\infty,1]\times\R^{N-1})\cup B_R^+}):
\mathop{\rm supp}v\Subset 
\R^N\setminus\{(1,x')\in \R\times\R^{N-1}:|x'|>R\}
\Big\}
$$
with respect to the norm $\big(\int_{((-\infty,1)\times\R^{N-1})\cup
  B_R^+}|\nabla v|^2dx\big)^{1/2}$ (which is actually equivalent to
the norm $\big(\int_{((-\infty,1)\times\R^{N-1})\cup B_R^+}|\nabla
v|^2dx+\int_{\Gamma_R^+}v^2d\sigma \big)^{1/2}$ by Poincaré
inequality), i.e.  $\mathcal H_R^+$ is the space of functions with
finite energy in $((-\infty,1]\times\R^{N-1})\cup B_R^+$ vanishing on
$\{(1,x')\in \R\times\R^{N-1}:|x'|\geq R\}$.

\begin{Lemma}\label{l:tilde_u_eps_bound}
For every sequence $\e_{n}\to 0^+$ there 
   exist a subsequence
  $\{\e_{n_k}\}_k$ and $\widetilde u\in \bigcup_{R>2}\mathcal H_R^+$ such
  that 
\begin{itemize}
\item[i)]
$\widetilde u_{\e_{n_k}}\to\widetilde u$ strongly in $\mathcal
  H_R^+$ for every $R>2$ and a.e.;
\item[ii)]
 $\widetilde u\equiv 0$ in $\R^N\setminus
  \widetilde D$;
\item[iii)]  $\widetilde u$ weakly solves
\begin{equation}\label{eq:eqtilde_lim}
\begin{cases}
-\Delta \widetilde u(x)=0,&\text{in }\widetilde D,\\
\widetilde u=0,&\text{on }\partial\widetilde D,
\end{cases}
\end{equation}
with
$\widetilde D$ as in (\ref{eq:94});
\item[iv)] $\widetilde u(x)\geq \frac{C_5}{2}(x_1-1)$ for all $x\in D^+\setminus B_2^+$.
\end{itemize}
\end{Lemma}
\begin{pf}
  Let $R>2$.  From Lemma \ref{l:sup_sol_bound} and (\ref{eq:26}),
  there exists  $C_R>0$ such that 
\begin{equation}\label{eq:87}
\int_{\Gamma_R^+}|\widetilde u_\e|^2d\sigma
=\frac{1}{\e^2}\int_{\Gamma_R^+}u_\e^2( {\mathbf e}_1+\e( x- {\mathbf
  e}_1))d\sigma\leq C_3^2\int_{\Gamma_R^+} \bigg(
\Phi_1(x)+2\gamma_\e \Phi_2\Big(\frac{x+{\mathbf e}_1}{2}\Big)\bigg)^{\!\!2}d\sigma
\leq C_R
\end{equation}
for all $\e\in (0,r_0/R)$.  By the change of variable $x={\mathbf e}_1
+\e(y-{\mathbf e}_1)$ we have that 
\begin{equation}\label{eq:88}
\mathcal N_\e(1+R\e)=
\frac{
R\int_{\widetilde{\Omega}_{R+1}^\e}\Big(|\nabla
 \widetilde u_\e(y)|^2-\lambda^\e_{\bar k}\e^2 p({\mathbf e}_1
+\e(y-{\mathbf e}_1)) \widetilde u_\e^2(y)\Big)dy}
{\int_{\Gamma_R^+}\widetilde u_\e^2(y)\,d\sigma}
\end{equation}
where 
$$
\widetilde{\Omega}_{R+1}^\e:=
\Big\{(y_1,y')\in \R\times\R^{N-1}:y_1<1-\tfrac1\e\Big\}
\cup \Big\{(y_1,y')\in \R\times\R^{N-1}:1-\tfrac1\e\leq y_1\leq 1,\ y'\in\Sigma\Big\}
\cup B_R^+.
$$
From Corollary \ref{cor:stima_N+}
\begin{equation}\label{eq:89}
\mathcal N_\e(1+R\e)\leq \mathcal N_\e(1+r_0)+\frac{C_9}{N 2^N}
\end{equation}
for all $\e\in (0,\min\{r_0/R,\e_2\})$.  From the strong $\mathcal D^{1,2}(\R^N)$
convergence of $u_\e$ to $\varphi_{k_0}^+$ ensured by Lemma
\ref{l:conv_eigen}, we deduce that there exists some positive constant
$C_{10}>0$ (depending on $r_0$ but independent of $\e$) such that 
$\mathcal N_\e(1+r_0)\leq C_{10}$ for all $\e\in (0,\e_0)$, so that 
(\ref{eq:87}--\ref{eq:89}) yield
\begin{align}\label{eq:90}
  \int_{\widetilde{\Omega}_{R+1}^\e}\Big(|\nabla \widetilde
  u_\e(y)|^2-\lambda^\e_{\bar k}\e^2 p({\mathbf e}_1 +\e(y-{\mathbf
    e}_1)) \widetilde u_\e^2(y)\Big)dy&\leq \bigg(C_{10}+\frac{C_9}{N
    2^N}\bigg)\frac{\int_{\Gamma_R^+}\widetilde u_\e^2(y)\,d\sigma}{R}\\
  \notag&\leq \bigg(C_{10}+\frac{C_9}{N 2^N}\bigg)\frac{C_R}{R}
\end{align}
for all $\e\in \big(0,\min\{r_0/R,\e_2\})$. From (\ref{eq:90}), Lemma
\ref{l:lemmaZERO}, and assumption (\ref{eq:p2}), we obtain that
\begin{align}\label{eq:91}
\int_{\widetilde{\Omega}_{R+1}^\e}|\nabla \widetilde
u_\e(y)|^2dy\leq 2\bigg(C_{10}+\frac{C_9}{N
  2^N}\bigg)\frac{C_R}{R}
\end{align}
for all $\e\in (0,\min\{r_0/R,\e_2,\bar\e_{2,2\lambda_{k_0}(D^+)p}\})$.
In view of (\ref{eq:87}) and (\ref{eq:91}), we have proved that for
every $R>2$ there exists $\e_R>0$ such that
\begin{equation}\label{eq:92}
\{\widetilde
u_\e\}_{\e\in (0,\e_R)} \text{ is bounded in }\mathcal H^+_R.
\end{equation}
Let $\e_n\to 0^+$.
From (\ref{eq:92}) and a diagonal
process, we deduce that  there exist a
subsequence $\e_{n_k}\to 0^+$  and some $\widetilde u\in
\bigcup_{R>2}\mathcal H_R^+$ such that $\widetilde
u_{\e_{n_k}}\rightharpoonup\widetilde u$ weakly in $\mathcal H_R^+$ for every
$R>2$. In particular $\widetilde u_{\e_{n_k}}\to \widetilde u$ a.e., so
that $\widetilde u\equiv 0$ in $\R^N\setminus \widetilde D$.  Passing to the weak limit in
(\ref{eq:eqtilde}), we obtain that $\widetilde u$ is a weak solution
(\ref{eq:eqtilde_lim}).
By classical elliptic estimates, we also have that $\widetilde
u_{\e_{n_k}}\to \widetilde u$ in $C^2(\overline{B_{r_2}^+\setminus
  B_{r_1}^+})$ for all $1<r_1<r_2$. Therefore, multiplying
(\ref{eq:eqtilde_lim}) by $\widetilde u$ and integrating over
$T_1^-\cup B_R^+$ with $T_1^-$ as in (\ref{eq:94}), we obtain
\begin{align}\label{eq:93}
\int_{\Gamma_R^+}\frac{\partial \widetilde u_{\e_{n_k}}}{\partial \nu}\widetilde u_{\e_{n_k}}d\sigma
\to
\int_{\Gamma_R^+}\frac{\partial \widetilde u}{\partial \nu}\widetilde u\,d\sigma
=\int_{T_1^-\cup B_R^+}|\nabla \widetilde u(x)|^2dx\quad\text{as }k\to+\infty .
\end{align}
On the other hand,  multiplication of
(\ref{eq:eqtilde}) by $\widetilde u_{\e_{n_k}}$ and integration by parts over
$\widetilde{\Omega}_{R+1}^{\e_{n_k}}$ yield
\begin{align}\label{eq:98}
  \int_{\widetilde{\Omega}_{R+1}^{\e_{n_k}}}|\nabla \widetilde u_{\e_{n_k}}(x)|^2dx=
  \int_{\Gamma_R^+}\frac{\partial \widetilde u_{\e_{n_k}}}{\partial
    \nu}\widetilde u_{\e_{n_k}}d\sigma+ \lambda^{\e_{n_k}}_{\bar k}\e_{n_k}^2
  \int_{\widetilde{\Omega}_{R+1}^{\e_{n_k}}} p({\mathbf e}_1 +{\e_{n_k}}(x-{\mathbf e}_1))
  \widetilde u_{\e_{n_k}}^2(x)\,dx.
\end{align}
We claim that 
\begin{align}\label{eq:99}
\e_{n_k}^2
  \int_{\widetilde{\Omega}_{R+1}^{\e_{n_k}}} p({\mathbf e}_1 +{\e_{n_k}}(x-{\mathbf e}_1))
  \widetilde u_{\e_{n_k}}^2(x)\,dx\to0\quad \text{as }k\to+\infty.
\end{align}
Indeed, from Lemma \ref{l:lemmaZERO}, for every $\delta>0$ there exists $k_0$ such that 
for all $k\geq k_0$ 
\begin{align*}
\int_{\Omega_{1/2}^{\e_{n_k}}}p(y)u_{\e_{n_k}}^2(y)dy\leq\delta
\int_{\Omega_{1/2}^{\e_{n_k}}}|\nabla u_{\e_{n_k}}(y)|^2dy
\end{align*}
and hence, 
from the change of variable $y={\mathbf e}_1 +{\e_{n_k}}(x-{\mathbf e}_1)$,
  assumption (\ref{eq:p2}), and (\ref{eq:91}), we deduce that 
\begin{align*}
\e_{n_k}^2
&  \int_{\widetilde{\Omega}_{R+1}^{\e_{n_k}}} p({\mathbf e}_1 +{\e_{n_k}}(x-{\mathbf e}_1))
  \widetilde u_{\e_{n_k}}^2(x)\,dx=
\e_{n_k}^{-N}\int_{\Omega_{1+R \e_{n_k}}^{\e_{n_k}}}p(y)u_{\e_{n_k}}^2(y)dy\\
&=
\e_{n_k}^{-N}\int_{\Omega_{1/2}^{\e_{n_k}}}p(y)u_{\e_{n_k}}^2(y)dy
\leq \delta\e_{n_k}^{-N}\int_{\Omega_{1+R \e_{n_k}}^{\e_{n_k}}}|\nabla u_{\e_{n_k}}(y)|^2dy\\
&=\delta \int_{\widetilde{\Omega}_{R+1}^{\e_{n_k}}}|\nabla \widetilde u_{\e_{n_k}}(x)|^2dx
\leq 2\delta \bigg(C_{10}+\frac{C_9}{N
  2^N}\bigg)\frac{C_R}{R},
\end{align*}
thus proving claim (\ref{eq:99}). Combining (\ref{eq:93}),
(\ref{eq:98}), and (\ref{eq:99}), we conclude that $\| \widetilde
u_{\e_{n_k}}\|_{\mathcal H_R^+}\to\| \widetilde u\|_{\mathcal H_R^+}$ and
then $\widetilde u_{\e_{n_k}}\to\widetilde u$ strongly in $\mathcal
H_R^+$ for every $R>2$.

To prove iv), it is enough to observe that Lemma \ref{l:usot} implies
that, for $k$ large,  
$$
\widetilde u_{\e_{n_k}}(x)\geq \frac{C_5}{2}(x_1-1)\quad\text{for all
}x\in B_{{r_0}/{\e_{n_k}}}^+\setminus B_2^+,
$$
which yields iv) thanks to  a.e convergence of $\widetilde u_{\e_{n_k}}$ to $\widetilde u$.
\end{pf}

\begin{remark}\label{rem:inftyenergy}
We notice that the function  $\widetilde u$ found in  Lemma \ref{l:tilde_u_eps_bound}
satisfies
$$
\int_{\widetilde D}|\nabla\widetilde u(x)|^2dx=+\infty.
$$
Indeed, $\int_{\widetilde D}|\nabla\widetilde u(x)|^2dx<+\infty$ would
imply, by testing (\ref{eq:eqtilde_lim}) with $\widetilde u$, that
$\widetilde u\equiv 0$ in $\widetilde D$, thus contradicting statement
iv) of Lemma \ref{l:tilde_u_eps_bound}.
\end{remark}

\begin{Lemma}\label{l:lim_freq_tildeu}
  Let $\widetilde u$ be as  in Lemma \ref{l:tilde_u_eps_bound} and,  for 
 $r\in \R\setminus (1,2)$, let $\widetilde{\mathcal N}_{\widetilde u}(r)$ 
 be  the frequency function associated to $\widetilde u$, i.e. 
\begin{equation*}
\widetilde{\mathcal N}_{\widetilde u}(r)=\frac{\Lambda_N(r)
\int_{\widetilde\Omega_r}|\nabla \widetilde u(x)|^2dx}
{\int_{\widetilde\Gamma_r}\widetilde u^2(x)\,d\sigma},
\end{equation*}
with $\widetilde\Omega_r$ and $\widetilde\Gamma_r$ defined in
(\ref{eq:115})
and $\Lambda_N(r)$ as in \eqref{eq:L_N}.
Then
\begin{itemize}
\item[i)] $\lim_{r\to+\infty}\widetilde{\mathcal N}_{\widetilde u}(r)=1$;
\item[ii)] there exists $\tilde c>0$ such that $\int_{D^+}|\nabla
  (\widetilde u-\tilde c (x_1-1))(x)|^2\,dx<+\infty$.
\end{itemize}
\end{Lemma}
\begin{pf} 
We notice that $\widetilde{\mathcal N}_{\widetilde u}$ is well defined
in $\R\setminus (1,2)$ in view of equation \eqref{eq:eqtilde_lim} and
classical unique continuation (in particular $\widetilde u\not\equiv0$
by part iv) of Lemma \ref{l:tilde_u_eps_bound})).
Let us first prove that 
\begin{equation}\label{eq:105}
\limsup_{r\to+\infty}\widetilde{\mathcal N}_{\widetilde u}(r)\leq1.
\end{equation}
Indeed, letting  $\e_{n}\to 0^+$ and 
  $\{\e_{n_k}\}_k$ as in Lemma \ref{l:tilde_u_eps_bound},
 passing to the limit as $k\to+\infty$ in (\ref{eq:88}), and using (\ref{eq:99}), we have that 
\begin{equation*}
  \lim_{k\to+\infty}\mathcal N_{\e_{n_k}}(1+R\e_{n_k})
  =\widetilde{\mathcal{N}}_{\widetilde u}(1+R)\quad\text{for every }R>0,
\end{equation*}
which, together with Corollary \ref{cor:N+sim1}, implies for every
$\delta>0$ the existence of some $\widetilde R_\delta$ such that
$$
\widetilde{\mathcal N}_{\widetilde
  u}(1+R)\leq1+\delta\quad\text{ for all $R>\widetilde
  R_\delta$},
$$
thus proving claim (\ref{eq:105}).

It is easy to prove that there exists $g\in H^1_{\rm loc}(D^+)$ such that 
\begin{align*}
\begin{cases}
-\Delta g=0,&\text{in }D^+,\\
g=\widetilde u,&\text{on }\partial D^+,\\
\int_{D^+}|\nabla g(x)|^2dx<+\infty,
\end{cases}
\end{align*}
i.e. $g$ is a finite-energy harmonic extension of $\widetilde
u\big|_{\partial D^+}$ in $D^+$.  We observe that the Kelvin transform
$\widetilde g(x)=|x-{\mathbf e}_1|^{-(N-2)}g\big( \frac{x-{\mathbf
    e}_1}{|x-{\mathbf e}_1|^2}+{\mathbf e}_1\big)$ belongs to
$H^1(B_1^+)$ and weakly satisfies
$$
\begin{cases}
-\Delta \widetilde g(x)=0,&\text{in }B_1^+,\\
\widetilde g(x)=0,&\text{on }\{(x_1,x'):x_1=1,|x'|<1\}.
\end{cases}
$$
By classical elliptic estimates, there exists
$c_{g}>0$ such that $\big|\frac{\partial \widetilde g}{\partial
  x_1}\big|\leq c_g$ in $\overline{B^+_{1/2}}$, thus implying 
$$
|\widetilde g(x_1,x')|=\bigg|\widetilde g(1,x')+\int_1^{x_1}
\frac{\partial \widetilde g}{\partial
  x_1}(s,x')\,ds\bigg|\leq \int_1^{x_1}
\bigg|\frac{\partial \widetilde g}{\partial
  x_1}(s,x')\bigg|\,ds\leq c_g(x_1-1)
$$
for all $(x_1,x')\in \overline{B^+_{1/2}}$. Then 
\begin{equation}\label{eq:104}
|g(x)|\leq c_g\,\frac{x_1-1}{|x-{\mathbf e}_1|^{N}}
\end{equation}
for all $x\in D^+\setminus B_{2}^+$.  Let us observe that the function
$v:=\widetilde u-g\in H^1_{\rm loc}(D^+)\setminus\{0\}$ satisfies
\begin{equation}\label{eq:112}
\begin{cases}
-\Delta v(x)=0,&\text{in }D^+,\\
v=0,&\text{on }\partial D^+,\\
\int_{B_r^+}|\nabla v(x)|^2dx<+\infty,&\text{for all }r>0.
\end{cases}
\end{equation}
Let us  define
$$
\mathcal N_v:(0,+\infty)\to\R,\quad
\mathcal N_v(t):= \frac{t \int_{B_t^+}|\nabla v(x)|^2dx}
{\int_{\Gamma_t^+}v^2(x)\,d\sigma}.
$$
Direct calculations yield
\begin{align}\label{eq:111}
{\mathcal N}_v'(t) = \frac{2t\bigg[ \left(\int_{\Gamma^+_t}
\left|\frac{\partial v}{\partial\nu}\right|^2
        d\sigma\right) \left(\int_{\Gamma_t^+} v^2
        d\sigma\right)
-\left(\int_{\Gamma_t^+}
        v\frac{\partial v}{\partial \nu}
        d\sigma\right)^{2} \bigg]} {\left(\int_{\Gamma_t^+} v^2
      d\sigma\right)^{2}},\quad\text{for all }t>0,
\end{align}
where $\nu=\nu(x)=\frac{x-{\mathbf e}_1}{|x-{\mathbf e}_1|}$.  In
particular, Schwarz's inequality implies that $\mathcal N_v$ is non
decreasing in $(0,+\infty)$.  From Remark \ref{rem:poincare_inv} it
follows that 
\begin{equation}\label{eq:102}
\mathcal N_v(t)\geq\lim_{r\to 0^+}\mathcal N_v(r)\geq 1
\quad\text{for all }t>0.
\end{equation}
From (\ref{eq:104}) and Lemma \ref{l:tilde_u_eps_bound}, it follows
that, if $x\in \Gamma_t^+$ and $t>2$, then 
$$
\bigg(1-\frac{2c_g}{C_5t^N}\bigg)\widetilde u(x)\leq v(x)\leq
\bigg(1+\frac{2c_g}{C_5t^N}\bigg)\widetilde u(x),
$$
so that 
$$
\bigg(1-\frac{2c_g}{C_5t^N}\bigg)^{\!\!2} \int_{\Gamma_t^+}\widetilde
u^2d\sigma \leq \int_{\Gamma_t^+}v^2d\sigma \leq
\bigg(1+\frac{2c_g}{C_5t^N}\bigg)^{\!\!2} \int_{\Gamma_t^+}\widetilde
u^2d\sigma
$$
 for all $t>\max\big\{2, (2c_g/C_5)^{1/N}\big\}$.  Let us fix
$\delta>0$. For every $R>2$ there holds
\begin{multline*}
  \int_{B^+_R}|\nabla v(x)|^2dx-(1+\delta)\int_{\widetilde \Omega_{R+1}}|\nabla \widetilde u(x)|^2dx\\
  \leq\int_{B^+_R}|\nabla g(x)|^2dx-2\int_{B^+_R}\nabla g(x)\cdot
  \nabla \widetilde u(x)dx
  -\delta \int_{\widetilde \Omega_{R+1}}|\nabla \widetilde u(x)|^2dx\\
  \leq \bigg(1+\frac2\delta\bigg)\int_{B^+_R}|\nabla
  g(x)|^2dx+\frac\delta2\int_{B^+_R}|\nabla
  \widetilde u(x)|^2dx-\delta \int_{\widetilde \Omega_{R+1}}|\nabla \widetilde u(x)|^2dx
\end{multline*}
and hence, for all $R>\max\big\{2, (2c_g/C_5)^{1/N}\big\}$,
\begin{align}\label{eq:107}
  \mathcal N_v(R)&\leq
  \frac{(1+\delta)}{\big(1-\frac{2c_g}{C_5}R^{-N}\big)^2}\frac{R
    \int_{\widetilde \Omega_{R+1}}|\nabla \widetilde u(x)|^2dx }
  {\int_{\Gamma_R^+}\widetilde u^2d\sigma}\left( 1+
    \frac{1+\frac2\delta}{1+\delta} \frac{\int_{B^+_R}|\nabla
      g(x)|^2dx}{\int_{\widetilde \Omega_{R+1}}|\nabla \widetilde
      u(x)|^2dx} \right)\\
  \notag&= \frac{(1+\delta)}{\big(1-\frac{2c_g}{C_5}R^{-N}\big)^2}
  \widetilde{\mathcal N}_{\widetilde u}(R+1) \left(1+
    \frac{1+\frac2\delta}{1+\delta} \frac{\int_{B^+_R}|\nabla
      g(x)|^2dx}{\int_{\widetilde \Omega_{R+1}}|\nabla \widetilde
      u(x)|^2dx} \right).
\end{align}
On the other hand, for every $R>2$ there holds
\begin{multline*}
  \int_{B^+_R}|\nabla v(x)|^2dx-(1-\delta)\int_{\widetilde \Omega_{R+1}}|\nabla \widetilde u(x)|^2dx\\
  =-\int_{T_1^-}|\nabla \widetilde u(x)|^2dx+ \int_{B^+_R}|\nabla
  g(x)|^2dx-2\int_{B^+_R}\nabla g(x)\cdot \nabla \widetilde u(x)dx
  +\delta \int_{\widetilde \Omega_{R+1}}|\nabla \widetilde u(x)|^2dx\\
  \geq -\int_{T_1^-}|\nabla \widetilde u(x)|^2dx+
  \bigg(1-\frac2\delta\bigg)\int_{B^+_R}|\nabla g(x)|^2dx
\end{multline*}
and hence, for all $R>\max\big\{2, (2c_g/C_5)^{1/N}\big\}$,
\begin{align}\label{eq:106}
  &\mathcal N_v(R)
\geq\frac{(1-\delta) \, \widetilde{\mathcal N}_{\widetilde u}(R+1)
}{\big(1+\frac{2c_g}{C_5}R^{-N}\big)^2}
  \left(1-\frac{1}{1-\delta}\frac{\int_{T_1^-}|\nabla \widetilde
      u(x)|^2dx} {\int_{\widetilde \Omega_{R+1}}|\nabla \widetilde
      u(x)|^2dx} + \frac{1-\frac2\delta}{1-\delta}
    \frac{\int_{B^+_R}|\nabla g(x)|^2dx}{\int_{\widetilde
        \Omega_{R+1}}|\nabla \widetilde u(x)|^2dx} \right).
\end{align}
Since $\int_{B^+_R}|\nabla g(x)|^2dx=O(1)$ and $\int_{\widetilde
  \Omega_{R+1}}|\nabla \widetilde u(x)|^2dx\to+\infty$ as
$R\to+\infty$ (see Remark \ref{rem:inftyenergy}), passing to 
$\limsup$ and $\liminf$  the in (\ref{eq:107}--\ref{eq:106})  we obtain that 
\begin{align*}
&(1-\delta)\limsup_{R\to\infty}
\widetilde{\mathcal N}_{\widetilde u}(R)\leq
\limsup_{R\to\infty}
\mathcal N_v(R)\leq 
(1+\delta)\limsup_{R\to\infty}
\widetilde{\mathcal N}_{\widetilde u}(R)
\quad\text{for all }\delta>0,\\
&(1-\delta)\liminf_{R\to\infty}
\widetilde{\mathcal N}_{\widetilde u}(R)\leq
\liminf_{R\to\infty}
\mathcal N_v(R)\leq 
(1+\delta)\liminf_{R\to\infty}
\widetilde{\mathcal N}_{\widetilde u}(R)
\quad\text{for all }\delta>0,
\end{align*}
thus implying, in view of from (\ref{eq:102}),
\begin{align}\label{eq:108}
  \liminf_{R\to\infty} \widetilde{\mathcal N}_{\widetilde u}(R)= \liminf_{R\to\infty}
  \mathcal N_v(R)\geq1,
\end{align}
and, in view of (\ref{eq:105}),
\begin{align}\label{eq:109}
1\geq  \limsup_{R\to\infty} \widetilde{\mathcal N}_{\widetilde u}(R) =\limsup_{R\to\infty}
  \mathcal N_v(R).
\end{align}
From (\ref{eq:108}) and (\ref{eq:109}) we deduce that 
\begin{equation}\label{eq:110}
 \lim_{R\to\infty} \widetilde{\mathcal N}_{\widetilde u}(R)= \lim_{R\to\infty}
  \mathcal N_v(R)=1,
\end{equation}
thus proving statement i). Furthermore (\ref{eq:110}), (\ref{eq:102}), and the fact that 
 $\mathcal N_v$ is non decreasing imply that 
\begin{equation}\label{eq:113}
  \mathcal N_v(t)\equiv 1 \text{ in }(0,+\infty).
\end{equation}
Therefore ${\mathcal
  N}_v'(t)=0$ for any $t>0$.  From (\ref{eq:111}) we obtain
\begin{equation*}
  \bigg(\int_{\Gamma^+_t}
  \bigg|\frac{\partial v}{\partial\nu}\bigg|^2
  d\sigma\bigg) \bigg(\int_{\Gamma_t^+} v^2
    d\sigma\bigg)
    =\bigg(\int_{\Gamma_t^+}
    v\frac{\partial v}{\partial \nu}
    d\sigma\bigg)^{\!\!2}
    \quad 
    \text{for all } t>0,
\end{equation*}
which implies that $v$ and $\frac{\partial v}{\partial \nu}$
are linearly dependent as vectors in $L^2(\Gamma_t^+)$, i.e.
there exists a  function $\eta=\eta(t)$ such that
$\frac{\partial v}{\partial \nu}({\mathbf e}_1+t\theta)=\eta(t) v({\mathbf e}_1+t\theta)$ for $t>0$.
After integration we obtain
\begin{equation*}
  v({\mathbf e}_1+t\theta)=e^{\int_1^t \eta(s)ds} v({\mathbf e}_1+\theta)
  =\varphi(t) \psi(\theta) \quad  
t>0, \ \theta\in {\mathbb S}^{N-1}_+,
\end{equation*}
where 
${\mathbb S}^{N-1}_+:=\{\theta=(\theta_1,\theta_2,\dots,\theta_N)\in {\mathbb S}^{N-1}: \theta_1>0 \}$,
$\varphi(t)=e^{\int_1^t \eta(s)ds}$ and $\psi(\theta)=v({\mathbf e}_1+\theta)$.
Since $v$ satisfies (\ref{eq:112}), then 
$$
\left(\varphi''(t)+\frac{N-1}{t} \varphi'(t) \right)\psi(\theta)
+\frac{\varphi(t)}{t^2} \Delta_{\SN} \psi(\theta)=0.
$$
Taking $t$ fixed, we deduce that $\psi$ is an
eigenfunction of the operator $-\Delta_{\SN}$ on ${\mathbb S}^{N-1}_+$
under null Dirichlet boundary conditions on $\partial {\mathbb
  S}^{N-1}_+$, i.e. there exists $K_0\in\N$, $K_0\geq 1$, such that 
\begin{equation}\label{eq:114}
\begin{cases}
-\Delta_{\SN}\psi=K_0(N-2+K_0)\psi,&\text{in }{\mathbb S}^{N-1}_+,\\
\psi=0,&\text{on }\partial{\mathbb S}^{N-1}_+.
\end{cases}
\end{equation}
Then $\varphi(t)$ solves the equation
$$
\varphi''(t)+\frac{N-1}{t} \varphi(t)-\frac{K_0(N-2+K_0)}{t^2}\varphi(t)=0
$$
and hence $\varphi$ is of the form
$$
\varphi(r)=c_1 t^{K_0}+c_2 t^{-(N-2)-K_0},
$$
for some $c_1,c_2\in\R$.
Since, by elliptic regularity theory, $v$ is smooth in $\overline{D^+}$,  $c_2$ must be $0$
 and  $\varphi(t)=c_1
t^{K_0}$. Since $\varphi(1)=1$, we obtain that $c_1=1$ and
then
\begin{equation} \label{expw}
  v({\mathbf e}_1+t\theta)=t^{K_0}\psi(\theta),  \quad 
\text{for all }t>0\text{ and }\theta\in {\mathbb S}^{N-1}_+.
\end{equation}
Substituting (\ref{expw}) into (\ref{eq:113}), we find that
$1\equiv\mathcal N_v(t)\equiv K_0$ and therefore
$K_0=1$. Being $N-1$ the first eigenvalue of problem (\ref{eq:114}),
$\psi$ is simple. Hence there exists $\tilde c\in\R$ such that
$\psi(\theta)=\tilde c\theta_1^+$ and $v(x)=\tilde c(x_1-1)^+$. Lemma
\ref{l:tilde_u_eps_bound} part iv) and estimate \eqref{eq:104} imply
that $\tilde c>0$, 
thus proving ii).
\end{pf}

\begin{Corollary}\label{cor:uniqueness_limit}
Let $\widetilde u$ be as in Lemma \ref{l:tilde_u_eps_bound}
and $\tilde c$ as in Lemma \ref{l:lim_freq_tildeu}. Then 
$$
\widetilde u=
\tilde c\,\mathcal T(x_1-1)=
\tilde c\, \Phi_1
$$
where $\Phi_1$ is defined in (\ref{eq:Phi1}).
\end{Corollary}
\begin{pf}
It follows from Lemmas \ref{l:tilde_u_eps_bound} and \ref{l:lim_freq_tildeu}, taking into account Lemma 
\ref{l:sup_right} and the fact that $\mathcal T(c\psi)=c\mathcal T(\psi)$.
\end{pf}

\begin{Lemma}\label{l:limNepscorr}
For every $R>0$ 
$$
\lim_{\e\to0^+}\mathcal N_\e(1-R\e)=\widetilde{\mathcal N}(1-R),
$$
with $\widetilde{\mathcal N}$ as in (\ref{eq:freq_Phi1}).
\end{Lemma}
\begin{pf}
  Fix $R>0$. Let $\e_n\to0^+$. From Lemma \ref{l:tilde_u_eps_bound}
  and Corollary \ref{cor:uniqueness_limit}, there exist a subsequence
  $\{\e_{n_k}\}_k$ and $\tilde c>0$ such that $\widetilde
  u_{\e_{n_k}}\to\tilde c\,\Phi_1$ strongly in $\mathcal H_r^+$ for
  every $r>2$.
 By the change of variable $x={\mathbf e}_1
+\e(y-{\mathbf e}_1)$, we have that, for $\e<\frac1R$, 
\begin{equation}\label{eq:117}
\mathcal N_\e(1-R\e)=
\frac{
\int_{\widetilde{\Omega}_{1-R}^\e}\Big(|\nabla
 \widetilde u_\e(y)|^2-\lambda^\e_{\bar k}\e^2 p({\mathbf e}_1
+\e(y-{\mathbf e}_1)) \widetilde u_\e^2(y)\Big)dy}
{\int_{\widetilde \Gamma_{1-R}}\widetilde u_\e^2(y)\,d\sigma}
\end{equation}
where 
$\widetilde \Gamma_{1-R}$ is defined in (\ref{eq:115}) and 
$$
\widetilde{\Omega}_{1-R}^\e:=
\Big\{(y_1,y')\in \R\times\R^{N-1}:y_1<1-\tfrac1\e\Big\}
\cup \Big\{(y_1,y')\in \R\times\R^{N-1}:1-\tfrac1\e\leq y_1\leq 1-R,\ y'\in\Sigma\Big\}.
$$
From strong convergence of $\widetilde
  u_{\e_{n_k}}$ to $\tilde c\,\Phi_1$ in $\mathcal H_r^+$ for
  every $r>2$, passing to the limit in (\ref{eq:117}) along the subsequence $\{\e_{n_k}\}_k$ 
and using (\ref{eq:99}), we obtain that 
\begin{equation*}
  \lim_{k\to+\infty}\mathcal N_{\e_{n_k}}(1-R\e_{n_k})=
  \frac{
    \int_{\widetilde{\Omega}_{1-R}}|\nabla
    (\tilde c\,\Phi_1)(y)|^2dy}
  {\int_{\widetilde \Gamma_{1-R}}(\tilde c\,\Phi_1)^2(y)\,d\sigma}=
  \frac{
    \int_{\widetilde{\Omega}_{1-R}}|\nabla
    \Phi_1(y)|^2dy}
  {\int_{\widetilde \Gamma_{1-R}}\Phi_1^2(y)\,d\sigma}=\widetilde{\mathcal N}(1-R),
\end{equation*}
where $\widetilde{\Omega}_{1-R}$ is defined in \eqref{eq:115}.  Since
the limit depends neither on the sequence $\{\e_n\}_{n\in\N}$ nor on
its subsequence $\{\e_{n_k}\}_{k\in\N}$, we conclude that the
convergence actually holds as $\e\to 0^+$ thus proving the lemma.
\end{pf}

\begin{Lemma}\label{l:walk_on_corr}
For every $R>0$ and $\delta>0$, there exists
$\hat\e_{R,\delta}\in(0,\check \e)$ such that 
\begin{align*}
\mathcal N_\e(r)< (1+\delta)\sqrt{\lambda_1(\Sigma)}
\quad\text{for all }r\in(0,R\e]\text{ and }\e\in(0,\hat\e_{R,\delta}).
\end{align*}
\end{Lemma}
\begin{pf}
Let $\delta>0$ and  choose $\delta'>0$  sufficiently small such that $(1+\delta')^2e^{\delta'}<1+\delta$. 
From Corollary \ref{l:limNphi1}, there exists $\widehat R_\delta>0$ such that 
\begin{equation}\label{eq:119}
\widetilde{\mathcal N}(1-\widehat R_\delta)<(1+\delta')\sqrt{\lambda_1(\Sigma)}.
\end{equation}
From Lemma \ref{l:limNepscorr}, there exists $\e_\delta>0$ such that 
\begin{equation}\label{eq:120}
  \mathcal N_\e(1-\widehat R_\delta\, \e)< (1+\delta')\,
  \widetilde{\mathcal N}(1-\widehat R_\delta)
  \quad\text{for all }\e\in(0,\e_\delta).
\end{equation}
Let $R>0$. Letting $\bar\e_c^{\delta'}$ as in Lemma \ref{l:stima N_c}
and using (\ref{eq:118}), (\ref{eq:119}), and (\ref{eq:120}), for all
$r\in(0,R\e)$ and
$0<\e<\min\left\{\e_\delta,\bar\e_c^{\delta'},\frac1{R+\widehat
    R_\delta}\right\}$ we obtain
$$
\mathcal N_\e(r)\leq\mathcal N_\e(1-\widehat R_\delta\, \e)
e^{\delta'(1-\widehat R_\delta \e-r)}
<(1+\delta')^2e^{\delta'}\sqrt{\lambda_1(\Sigma)}<
(1+\delta)\sqrt{\lambda_1(\Sigma)}.
$$
The lemma is thereby proved.
\end{pf}

\section{Blow-up at the left}\label{sec:blow-up-at-left}

\noindent Let us define
\begin{align}\label{eq:121}
\widehat{u}_\e:\widehat \Omega^\e\to\R,\quad
\widehat u_\e(x)=\frac{u_\e(\e x)}{\sqrt{\e^{1-N}\int_{\Gamma^\e_{\e}}u_\e^2d\sigma}}
\end{align}
where 
\begin{equation*}
\widehat \Omega^\e:=\{x\in\R^N:\e x\in\Omega^\e\}
=D^-\cup\{(x_1,x')\in T_1:0\leq x_1\leq 1/\e\}\cup
\{(x_1,x'):x_1> 1/\e\}. 
\end{equation*} 
We observe that $\widehat u_\e$ solves
\begin{equation}\label{eq:eqhat}
\begin{cases}
-\Delta \widehat u_\e(x)=\e^2\lambda^\e_{\bar k} p(\e x) \widehat u_\e(x),&\text{in }\widehat \Omega^\e,\\
\widehat u_\e=0,&\text{on }\partial\widehat \Omega^\e.
\end{cases}
\end{equation}
We denote 
\begin{align*}
T_1^+=\{(x_1,x'):x'\in\Sigma,x_1\geq 0\},\quad \widehat D=D^-\cup T_1^+.
\end{align*}
For
every $R>0$ we define 
\begin{align}
\label{eq:124}
\widehat\Omega_R=D^-\cup \{(x_1,x')\in T_1^+:x_1<R\},\quad
\widehat\Gamma_R=\Gamma_R=\{(x_1,x')\in T_1^+:x_1=R\},
\end{align}
and $\mathcal H_R$ as the completion of
$$
\mathcal D_R:=\Big\{v\in C^\infty(\overline{\widehat\Omega_R}):
\mathop{\rm supp}v\Subset \widehat D\Big\}
$$
with respect to the norm $\big(\int_{\widehat\Omega_R}|\nabla
v|^2dx\big)^{\!1/2}$ 
 (which is  equivalent to $\big(\int_{\widehat\Omega_R}|\nabla
v|^2dx+\int_{\widehat\Gamma_R}v^2d\sigma \big)^{\!1/2}$),
 i.e.
$\mathcal H_R$ is the space of functions with finite energy in
$\widehat\Omega_R$ vanishing on $\{(x_1,x')\in \partial \widehat\Omega_R:
x_1< R\}$.

The change of variable $y'=\e x'$ yields
\begin{equation}\label{eq:normalization_eqhat}
\int_{\widehat\Gamma_1}\widehat u_\e^2d\sigma=1.
\end{equation}

\begin{Lemma}\label{l:hatgammabound}
  For every $R>1$, there exists  $\hat\e_{R}>0$ such that 
$$
\int_{\widehat\Gamma_R}\widehat u_\e^2d\sigma\leq
e^{4\sqrt{\lambda_1(\Sigma)}(R-1)}\quad\text{for all
}\e\in(0,\hat\e_{R}).
$$
\end{Lemma}
\begin{pf}
For $R>1$, let  $\hat\e_{R}=\hat\e_{R,1}>0$ as in Lemma \ref{l:walk_on_corr}. 
From Lemma \ref{l:walk_on_corr}, (\ref{eq:81}), and (\ref{eq:N_c}) it follows that 
$$
 \frac{\frac{d}{dr} H_\e^c(r)}{H_\e^c(r)}=\frac{2}{\e}\mathcal N_\e(r)
\leq \frac{4}{\e}\sqrt{\lambda_1(\Sigma)}
\quad\text{for all }r\in(0,R\e]\text{ and }\e\in(0,\hat\e_{R}),
$$
which after integration between $\e$ and $R\e$ yields
$$
H_\e^c(R\e)\leq H_\e^c(\e)e^{4\sqrt{\lambda_1(\Sigma)}(R-1)}
\quad\text{for all }\e\in(0,\hat\e_{R}).
$$
 (\ref{eq:121}) and the change of variable $y'=\e x'$ yield  
$$
\int_{\widehat\Gamma_R}\widehat u_\e^2d\sigma=
\frac{H_\e^c(R\e)}{H_\e^c(\e)}
$$
thus implying  the conclusion.
\end{pf}

\begin{Lemma}\label{l:hat_u_eps_bound}
For every sequence $\e_{n}\to 0^+$ there 
   exist a subsequence
  $\{\e_{n_k}\}_k$ and $\widehat u\in \bigcup_{R>1}\mathcal H_R$ such
  that
\begin{itemize}
\item[i)]
$\widehat u_{\e_{n_k}}\to\widehat u$ strongly in $\mathcal
  H_R$ for every $R>1$ and a.e.;
\item[ii)]
  $\widehat u\not\equiv 0$  in $\widehat D$;
\item[iii)]  $\widehat u$ weakly solves
\begin{equation}\label{eq:eqhat_lim}
\begin{cases}
-\Delta \widehat u(x)=0,&\text{in }\widehat D,\\
\widehat u=0,&\text{on }\partial\widehat D.
\end{cases}
\end{equation}
\end{itemize}
\end{Lemma}
\begin{pf}
  Let $R>1$. By the change of variable $x=\e y$ we have that, for
  $\e\in(0,\min\{1/R,\check \e\})$, 
\begin{equation}\label{eq:125}
\mathcal N_\e(R\e)=
\frac{
\int_{\widehat{\Omega}_{R}}\Big(|\nabla
 \widehat u_\e(y)|^2-\lambda^\e_{\bar k}\e^2 p(\e y) \widehat u_\e^2(y)\Big)dy}
{\int_{\widehat \Gamma_R}\widehat u_\e^2(y)\,d\sigma}.
\end{equation}
From Lemma \ref{l:walk_on_corr}, for every $\delta>0$ there exists $\hat\e_{R,\delta}>0$ such that 
\begin{align}\label{eq:126}
\mathcal N_\e(R\e)< (1+\delta)\sqrt{\lambda_1(\Sigma)}
\quad\text{for all }\e\in(0,\hat\e_{R,\delta}).
\end{align}
Choosing $\delta=1$, from (\ref{eq:125}), (\ref{eq:126}), and Lemma \ref{l:hatgammabound}, we have that 
\begin{equation}\label{eq:127}
\int_{\widehat{\Omega}_{R}}\Big(|\nabla
 \widehat u_\e(y)|^2-\lambda^\e_{\bar k}\e^2 p(\e y) \widehat u_\e^2(y)\Big)dy
\leq
2\sqrt{\lambda_1(\Sigma)}
\int_{\widehat \Gamma_R}\widehat u_\e^2(y)\,d\sigma
\leq 
2\sqrt{\lambda_1(\Sigma)}
e^{4\sqrt{\lambda_1(\Sigma)}(R-1)}
\end{equation}
for all $\e\in (0,\hat\e_{R})$, where $\hat\e_{R}=\hat\e_{R,1}>0$
(accordingly with the notation of Lemma \ref{l:hatgammabound}).  From
(\ref{eq:127}) and Lemma \ref{l:lemmaZERO}, we obtain that for all
$\e\in (0,\min\{\hat\e_{R},\bar\e_{2,2\lambda_{k_0}(D^+)p}\})$
\begin{align}\label{eq:128}
\int_{\widehat{\Omega}_{R}}|\nabla
 \widehat u_\e(y)|^2dy\leq 4\sqrt{\lambda_1(\Sigma)}
e^{4\sqrt{\lambda_1(\Sigma)}(R-1)}.
\end{align}
In view of (\ref{eq:128}) and Lemma \ref{l:hatgammabound}, we have  that for
every $R>1$ there exists $\e_R>0$ such that 
\begin{equation}\label{eq:129}
\{\widehat
u_\e\}_{\e\in (0,\e_R)} \text{ is bounded in }\mathcal H_R.
\end{equation}
Let $\e_n\to 0^+$.  From (\ref{eq:129}) and a diagonal process, we
deduce that there exist a subsequence $\e_{n_k}\to 0^+$ and some
$\widehat u\in \bigcup_{R>1}\mathcal H_R$ such that $\widehat
u_{\e_{n_k}}\rightharpoonup\widehat u$ weakly in $\mathcal H_R$ for
every $R>1$ and almost everywhere. From compactness of the embedding
$\mathcal H_R\hookrightarrow L^2(\widehat\Gamma_1)$ and
(\ref{eq:normalization_eqhat}) we deduce that
$\int_{\widehat\Gamma_1}\widehat u^2d\sigma=1$; in particular
$\widehat u\not\equiv 0$.

Passing to the weak limit in
(\ref{eq:eqhat}), we obtain that $\widehat u$ is a weak solution
(\ref{eq:eqhat_lim}).
By classical elliptic estimates, we also have that $\widehat
u_{\e_{n_k}}\to \widehat u$ in 
$C^2(\{(x_1,x')\in T_1^+:r_1\leq x_1\leq r_2\})$
 for all $0<r_1<r_2$. Therefore, multiplying
(\ref{eq:eqhat_lim}) by $\widehat u$ and integrating over
$\widehat{\Omega}_{R}$, we obtain
\begin{align}\label{eq:93left}
\int_{\widehat\Gamma_R}\frac{\partial \widehat u_{\e_{n_k}}}{\partial x_1}\widehat u_{\e_{n_k}}d\sigma
\to
\int_{\widehat\Gamma_R}\frac{\partial \widehat u}{\partial x_1}\widehat u\,d\sigma
=\int_{\widehat{\Omega}_{R}}|\nabla \widehat u(x)|^2dx.
\end{align}
On the other hand,  multiplication of
(\ref{eq:eqhat}) by $\widehat u_{\e_{n_k}}$ and integration by parts over
$\widehat{\Omega}_{R}$ yield
\begin{align}\label{eq:98left}
  \int_{\widehat{\Omega}_{R}}|\nabla \widehat u_{\e_{n_k}}(x)|^2dx=
  \int_{\widehat\Gamma_R}\frac{\partial \widehat u_{\e_{n_k}}}{\partial
    x_1}\widehat u_{\e_{n_k}}d\sigma+ \lambda^{\e_{n_k}}_{\bar k}\e_{n_k}^2
  \int_{\widehat{\Omega}_{R}} p(\e_{n_k}x)
  \widehat u_{\e_{n_k}}^2(x)\,dx.
\end{align}
We claim that, for every $R>1$, 
\begin{align}\label{eq:99left}
\e_{n_k}^2
  \int_{\widehat{\Omega}_{R}} p(\e_{n_k}x)
  \widehat u_{\e_{n_k}}^2(x)\,dx\to0\quad \text{as }n\to+\infty.
\end{align}
Indeed, from Lemma \ref{l:lemmaZERO}, for every $\delta>0$ there exists $k_0$ such that 
for all $k\geq k_0$ 
\begin{align*}
\int_{\Omega_{R\e_{n_k}}^{\e_{n_k}}}p(y)u_{\e_{n_k}}^2(y)dy\leq\delta
\int_{\Omega_{R\e_{n_k}}^{\e_{n_k}}}|\nabla u_{\e_{n_k}}(y)|^2dy
\end{align*}
and hence, 
from the change of variable $y=\e_{n_k}x$  and (\ref{eq:128}), we deduce that 
\begin{align*}
  \e_{n_k}^2 & \int_{\widehat{\Omega}_{R}} p(\e_{n_k}x) \widehat
  u_{\e_{n_k}}^2(x)\,dx =
  \frac{\e_{n_k}}{\int_{\Gamma^{\e_{n_k}}_{\e_{n_k}}}u_{\e_{n_k}}^2d\sigma}
  \int_{\Omega_{R\e_{n_k}}^{\e_{n_k}}}p(y)u_{\e_{n_k}}^2(y)dy\\
  &\leq
  \frac{\e_{n_k}\delta}{\int_{\Gamma^{\e_{n_k}}_{\e_{n_k}}}u_{\e_{n_k}}^2d\sigma}
  \int_{\Omega_{R\e_{n_k}}^{\e_{n_k}}}|\nabla u_{\e_{n_k}}(y)|^2dy\\
  &=\delta \int_{\widehat{\Omega}_{R}}|\nabla \widehat
  u_{\e_{n_k}}(x)|^2dx \leq 4\delta \sqrt{\lambda_1(\Sigma)}
  e^{4\sqrt{\lambda_1(\Sigma)}(R-1)},
\end{align*}
thus proving claim (\ref{eq:99left}). Combining (\ref{eq:93left}),
(\ref{eq:98left}), and (\ref{eq:99left}), we conclude that $\| \widehat
u_{\e_{n_k}}\|_{\mathcal H_R}\to\| \widehat u\|_{\mathcal H_R}$ and
then $\widehat u_{\e_{n_k}}\to\widehat u$ strongly in $\mathcal
H_R$ for every $R>1$.\end{pf}

\begin{remark}\label{rem:inftyenergy_hat}
We notice that the function  $\widehat u$ found in  Lemma \ref{l:hat_u_eps_bound}
satisfies
$$
\int_{\widehat D}|\nabla\widehat u(x)|^2dx=+\infty.
$$
Indeed, $\int_{\widehat D}|\nabla\widehat u(x)|^2dx<+\infty$ would
imply, by testing (\ref{eq:eqhat_lim}) with $\widehat u$, that
$\widehat u\equiv 0$ in $\widehat D$, thus contradicting statement
ii) of Lemma \ref{l:hat_u_eps_bound}.

We also observe that, denoting as $\widehat H(r)=\int_{\widehat \Gamma_r}\widehat u^2d\sigma$ for all $r>0$,
multiplication of  (\ref{eq:eqhat_lim}) by $\widehat u$ and integration over $\widehat\Omega_r$ yield
$$
\frac{d}{dr}\widehat H(r)=2\int_{\widehat \Gamma_r}\widehat
u\frac{\partial \widehat u}{\partial x_1}d\sigma
=2\int_{\widehat\Omega_r}|\nabla\widehat u(x)|^2dx\to \int_{\widehat D}|\nabla\widehat u(x)|^2dx=+\infty
\quad\text{as }r\to+\infty,
$$
thus implying that 
$$
\lim_{r\to+\infty}\widehat H(r)=
\lim_{r\to+\infty}\int_{\Gamma_r}\widehat u^2d\sigma=+\infty.
$$
\end{remark}

\begin{Lemma}\label{l:lim_freq_hatu}
  Let $\widehat u$ as  in Lemma \ref{l:hat_u_eps_bound} and,  for 
 $r>0$, let $\widehat{\mathcal N}_{\widehat u}(r)$ 
 be  the frequency function associated to $\widehat u$, i.e. 
\begin{equation*}
\widehat{\mathcal N}_{\widehat u}(r)=\frac{
\int_{\widehat\Omega_r}|\nabla \widehat u(x)|^2dx}
{\int_{\widehat\Gamma_r}\widehat u(x)\,d\sigma},\quad r>0,
\end{equation*}
with $\widehat\Omega_r$ and $\widehat\Gamma_r$ defined in (\ref{eq:124}).
Then
\begin{itemize}
\item[i)] $\lim_{r\to+\infty}\widehat{\mathcal N}_{\widehat u}(r)=\sqrt{\lambda_1(\Sigma)}$;
\item[ii)] there exists $\hat c\in\R\setminus\{0\}$ such that $\int_{T_1}|\nabla
  (\widehat u-\hat c h)(x)|^2\,dx<+\infty$,
\end{itemize}
where 
\begin{equation}\label{eq:145}
h:T_1\to\R,
\quad h(x_1,x')=f(1-x_1,x')=e^{\sqrt{\lambda_1(\Sigma)}x_1}\psi_1^\Sigma(x'),
\end{equation}
being $f$ defined in (\ref{eq:29}).
\end{Lemma}
\begin{pf}
Letting  $\e_{n}\to 0^+$ and 
  $\{\e_{n_k}\}_k$ as in Lemma \ref{l:hat_u_eps_bound},
 passing to the limit as $k\to+\infty$ in (\ref{eq:125}), and using (\ref{eq:99left}), we have that 
\begin{equation*}
\lim_{k\to+\infty}\mathcal N_{\e_{n_k}}(R\e_{n_k})=\widehat{\mathcal{N}}_{\widehat u}(R)\quad\text{for every }R>0,
\end{equation*}
which, together with (\ref{eq:126}), implies that, for every
$\delta>0$ and $R>0$,
$$
\widehat{\mathcal N}_{\widehat u}(R)\leq (1+\delta)\sqrt{\lambda_1(\Sigma)}.
$$
Therefore
\begin{equation}\label{eq:105left}
  \widehat{\mathcal N}_{\widehat u}(R)\leq\sqrt{\lambda_1(\Sigma)}\quad\text{for every }R>0.
\end{equation}
It is easy to prove that there exists $\zeta\in H^1_{\rm loc}(T_1)\cap L^\infty(T_1)$ such that 
\begin{align*}
\begin{cases}
-\Delta \zeta=0,&\text{in }T_1,\\
\zeta=\widehat u,&\text{on }\partial T_1,\\
\int_{T_1}|\nabla \zeta(x)|^2dx<+\infty,
\end{cases}
\end{align*}
i.e. $\zeta$ is a finite-energy harmonic extension of $\widehat
u\big|_{\partial T_1}$ in $T_1$.
Since $w(x_1,x')=e^{-\sqrt{\lambda_1(\Sigma)}\frac{x_1}2}\psi_1^\Sigma\big(\frac{x'}2\big)$
is harmonic and strictly positive in $T_1$, bounded from below away
from $0$ in $\{(x_1,x')\in T_1:x_1\leq0\}$, and $\int_{\{(x_1,x')\in T_1:\,x_1\geq r\}}(|\nabla
w|^2+|w|^{2^*})<+\infty$ for all $r$,  from the Maximum Principle we deduce
that $|\zeta|\leq {\rm const\,}w$ in $T_1$, thus implying that, 
for some $c_\zeta>0$,
\begin{equation*}
|\zeta(x)|\leq c_\zeta e^{-\sqrt{\lambda_1(\Sigma)}\frac{x_1}2}\quad\text{for all }x\in T_1.
\end{equation*}
 Let us observe that the function
$\widehat v:=\widehat u-\zeta\in H^1_{\rm loc}(T_1)\setminus\{0\}$ satisfies
\begin{equation*}
\begin{cases}
-\Delta \widehat v(x)=0,&\text{in }T_1,\\
\widehat v=0,&\text{on }\partial T_1.
\end{cases}
\end{equation*}
We notice that $\widehat v\not\equiv0$ in view of Remark \ref{rem:inftyenergy_hat}.
Let 
$$
N_{\widehat v}:\R\to\R,\quad
N_{\widehat v}(r):=\frac{\int_{T_{1,r}}|\nabla \widehat v(x)|^2dx}
{\int_{\Gamma_r}\widehat v^2(x)\,d\sigma},
$$
be as in Lemma \ref{l:limNphi1gen}, where, for all $r\in\R$, $T_{1,r}$
and $\Gamma_r$ are defined in (\ref{eq:141}). From Lemma
\ref{l:limNphi1gen} it follows that $N_{\widehat v}$ is non
decreasing in $\R$ and
\begin{equation}\label{eq:102left}
  N_{\widehat v}(t)\geq\lim_{r\to -\infty} N_{\widehat v}(r)\geq \sqrt{\lambda_1(\Sigma)}
  \quad\text{for all }t\in\R.
\end{equation}
For all $R>0$, $\delta\in(0,1)$,
\begin{multline*}
  \int_{\Gamma_R}\widehat v^2d\sigma-(1-\delta)
  \int_{\Gamma_R}\widehat u^2d\sigma
  =
\int_{\Gamma_R}\zeta^2d\sigma-2\int_{\Gamma_R}\zeta\widehat u d\sigma
+\delta\int_{\Gamma_R}\widehat u^2d\sigma\\
\geq  \bigg(1-\frac2\delta\bigg)\int_{\Gamma_R}\zeta^2d\sigma+\frac\delta2\int_{\Gamma_R}\widehat u^2d\sigma
\geq  \bigg(1-\frac2\delta\bigg)\int_{\Gamma_R}\zeta^2d\sigma
\geq  \bigg(1-\frac2\delta\bigg)\frac{\omega_{N-2}}{N-1}\,c_\zeta^2 e^{-\sqrt{\lambda_1(\Sigma)}R}
\end{multline*}
and
\begin{multline*}
  \int_{T_{1,R}}|\nabla \widehat v(x)|^2dx-(1+\delta)\int_{\widehat \Omega_{R}}|\nabla \widehat u(x)|^2dx\\
  \leq\int_{T_{1,R}}|\nabla \zeta(x)|^2dx-2\int_{T_{1,R}}\nabla \zeta(x)\cdot
  \nabla \widehat u(x)dx
  -\delta \int_{\widehat \Omega_{R}}|\nabla \widehat u(x)|^2dx\\
  \leq \bigg(1+\frac2\delta\bigg)\int_{T_{1,R}}|\nabla
  \zeta(x)|^2dx+\frac\delta2\int_{T_{1,R}}|\nabla
  \widehat u(x)|^2dx-\delta \int_{\widehat \Omega_{R}}|\nabla \widehat u(x)|^2dx\\
  \leq \bigg(1+\frac2\delta\bigg)\int_{T_{1,R}}|\nabla
  \zeta(x)|^2dx,
\end{multline*}
thus implying 
\begin{align}\label{eq:142}
N_{\widehat v}(R)\leq \frac{1+\delta}{1-\delta}
\,\widehat{\mathcal N}_{\widehat u}(R)
\frac{1+\dfrac{1+\frac2\delta}{1+\delta}\dfrac{\int_{T_{1,R}}|\nabla
  \zeta(x)|^2dx}{\int_{\widehat \Omega_{R}}|\nabla \widehat u(x)|^2dx}}
{1+\dfrac{(1-\frac2\delta)\omega_{N-2}c_\zeta^2 }
{(1-\delta)(N-1) \int_{\Gamma_R}\widehat u^2d\sigma}
e^{-\sqrt{\lambda_1(\Sigma)}R}
}.
\end{align}
On the other hand, for all $R>0$, $\delta\in(0,1)$,
\begin{multline*}
  \int_{\Gamma_R}\widehat v^2d\sigma-(1+\delta)
  \int_{\Gamma_R}\widehat u^2d\sigma
  =
\int_{\Gamma_R}\zeta^2d\sigma-2\int_{\Gamma_R}\zeta\widehat u d\sigma
-\delta\int_{\Gamma_R}\widehat u^2d\sigma\\
\leq  \bigg(1+\frac2\delta\bigg)\int_{\Gamma_R}\zeta^2d\sigma-\frac\delta2\int_{\Gamma_R}\widehat u^2d\sigma
\leq  \bigg(1+\frac2\delta\bigg)\int_{\Gamma_R}\zeta^2d\sigma
\leq  \bigg(1+\frac2\delta\bigg)\frac{\omega_{N-2}}{N-1}\,c_\zeta^2 e^{-\sqrt{\lambda_1(\Sigma)}R}
\end{multline*}
and 
\begin{multline*}
  \int_{T_{1,R}}|\nabla \widehat v(x)|^2dx-(1-\delta)\int_{\widehat \Omega_{R}}|\nabla \widehat u(x)|^2dx\\
  =-\int_{D^-\setminus T_1} |\nabla \widehat u(x)|^2dx+
  \int_{T_{1,R}}|\nabla \zeta(x)|^2dx-2\int_{T_{1,R}}\nabla \zeta(x)\cdot
  \nabla \widehat u(x)dx+\delta\int_{\widehat \Omega_{R}}|\nabla \widehat u(x)|^2dx\\
\geq
-\int_{D^-} |\nabla \widehat u(x)|^2dx+
\bigg(1-\frac2\delta\bigg)  \int_{T_{1,R}}|\nabla \zeta(x)|^2dx-\frac\delta2\int_{T_{1,R}}
|\nabla \widehat u(x)|^2dx+\delta\int_{\widehat \Omega_{R}}|\nabla \widehat u(x)|^2dx\\
\geq
-\int_{D^-} |\nabla \widehat u(x)|^2dx+
\bigg(1-\frac2\delta\bigg)  \int_{T_{1,R}}|\nabla \zeta(x)|^2dx,
\end{multline*}
thus implying 
\begin{align}\label{eq:143}
N_{\widehat v}(R)\geq \frac{1-\delta}{1+\delta}
\,\widehat{\mathcal N}_{\widehat u}(R)
\,\frac{1
-\dfrac{\int_{D^-} |\nabla \widehat u(x)|^2dx}{(1-\delta)\int_{\widehat \Omega_{R}}|\nabla \widehat u(x)|^2dx}
+\dfrac{1-\frac2\delta}{1-\delta}\dfrac{\int_{T_{1,R}}|\nabla
  \zeta(x)|^2dx}{\int_{\widehat \Omega_{R}}|\nabla \widehat u(x)|^2dx}}
{1+\dfrac{(1+\frac2\delta)\omega_{N-2}c_\zeta^2 }
{(1+\delta)(N-1) \int_{\Gamma_R}\widehat u^2d\sigma}
e^{-\sqrt{\lambda_1(\Sigma)}R}
}.
\end{align}
Since $\int_{T_{1,R}}|\nabla \zeta(x)|^2dx=O(1)$,  $\int_{\widehat
  \Omega_{R}}|\nabla \widehat u(x)|^2dx\to+\infty$, and 
$\int_{\Gamma_R}\widehat u^2d\sigma\to+\infty$
 as
$R\to+\infty$ (see Remark \ref{rem:inftyenergy_hat}), passing to 
$\limsup$ and $\liminf$  in (\ref{eq:142}--\ref{eq:143})  we obtain that 
\begin{gather*}
\frac{1-\delta}{1+\delta}\limsup_{R\to\infty}
\widehat{\mathcal N}_{\widehat u}(R)\leq
\limsup_{R\to\infty}
 N_{\widehat v}(R)\leq 
\frac{1+\delta}{1-\delta}\limsup_{R\to\infty}
\widehat{\mathcal N}_{\widehat u}(R)
\quad\text{for all }\delta>0,\\
\frac{1-\delta}{1+\delta}\liminf_{R\to\infty}
\widehat{\mathcal N}_{\widehat u}(R)\leq
\liminf_{R\to\infty}
N_{\widehat v}(R)\leq 
\frac{1+\delta}{1-\delta}\liminf_{R\to\infty}
\widehat{\mathcal N}_{\widehat u}(R)
\quad\text{for all }\delta>0,
\end{gather*}
thus implying, in view of (\ref{eq:102left}),
\begin{align}\label{eq:108left}
  \liminf_{R\to\infty} \widehat{\mathcal N}_{\widehat u}(R)= \liminf_{R\to\infty}
   N_{\widehat v}(R)\geq\sqrt{\lambda_1(\Sigma)}
\end{align}
and, in view of (\ref{eq:105left}),
\begin{align}\label{eq:109left}
\sqrt{\lambda_1(\Sigma)}\geq  \limsup_{R\to\infty} \widehat{\mathcal N}_{\widehat u}(R) 
=\limsup_{R\to\infty}
   N_{\widehat v}(R).
\end{align}
From (\ref{eq:108left}) and (\ref{eq:109left}) we deduce that 
\begin{equation}\label{eq:110left}
  \lim_{R\to\infty} \widehat{\mathcal N}_{\widehat u}(R)= \lim_{R\to\infty}
  N_{\widehat v}(R)=\sqrt{\lambda_1(\Sigma)},
\end{equation}
thus proving statement i). Furthermore (\ref{eq:110left}),
(\ref{eq:102left}), and the fact that $N_{\widehat v}$ is non
decreasing imply that
\begin{equation*}
N_{\widehat v}(t)\equiv \sqrt{\lambda_1(\Sigma)} \quad\text{in }\R.
\end{equation*}
From Lemma \ref{l:limNphi1gen} iii), it follows that there exists
$\hat c\in\R\setminus\{0\}$ such that $\widehat v(x_1,x')=\hat c h(x_1,x')$ with $h$
as in (\ref{eq:145}).  Since $\int_{T_1}|\nabla (\widehat u-\hat c
h)(x)|^2\,dx= \int_{T_1}|\nabla \zeta(x)|^2\,dx <+\infty$, also claim ii) is proved. 
\end{pf}

\begin{Corollary}\label{cor:uniqueness_limit_left}
Let $\widehat u$ be as in Lemma \ref{l:hat_u_eps_bound}
and $\hat c$ as in Lemma \ref{l:lim_freq_hatu}. Then 
$$
\widehat u(x_1,x')= \hat c\,\Phi_2 (1-x_1,x')
$$
where $\Phi_2$ is as in Lemma \ref{l:Phi2}.
\end{Corollary}
\begin{pf}
It follows from Lemmas \ref{l:hat_u_eps_bound} and \ref{l:lim_freq_hatu}, taking into account Lemma 
\ref{l:Phi2}.
\end{pf}

Let us define $\widehat \Phi(x_1,x'):= \Phi_2(1-x_1,x')$ and, for all $r<-1$ 
\begin{equation*}
\widehat{\mathcal N}(r)=
\widehat{\mathcal N}_{\widehat \Phi}(r)=
\frac{(-r)
\int_{\Omega_r}|\nabla \widehat \Phi(x)|^2dx}
{\int_{\Gamma_{-r}^-}\widehat \Phi(x)\,d\sigma},
\end{equation*}
with $\Omega_r$ as in \eqref{eq:defOmega_r} and $\Gamma^-_{-r}$ as in
\eqref{eq:defGamma_r-}, so that, according to notation of Lemma
\ref{l:limNmenophi1gen}, $\widehat{\mathcal N}(r)=N^-_{\widehat \Phi}(-r)$
for all $r<-1$.

\begin{Lemma}\label{l:limitN_hat}
$\lim_{r\to-\infty}\widehat{\mathcal N}(r)=N-1$.
\end{Lemma}
\begin{pf}
The proof follows from Lemma \ref{l:limNmenophi1gen} and Remark \ref{rem:phi2pos}.
\end{pf}

\begin{Lemma}\label{l:limNepsleft}
For every $R>1$ 
$$
\lim_{\e\to0^+}{\mathcal N}_\e(-R\e)=\widehat{\mathcal N}(-R).
$$
\end{Lemma}
\begin{pf}
  Fix $R>1$. Let $\e_n\to0^+$. From Lemma \ref{l:hat_u_eps_bound}
  and Corollary \ref{cor:uniqueness_limit_left}, there exist a subsequence
  $\{\e_{n_k}\}_k$ and $\hat c\neq0$ such that $\widehat
  u_{\e_{n_k}}\to\hat c\,\widehat\Phi$ strongly in $\mathcal H_r$ for
  every $r>1$.
 By the change of variable $x=\e y$ we have that, for $\e\in(0,\check \e)$ and
 $R>1$, 
\begin{equation}\label{eq:117left}
\mathcal N_\e(-R\e)=
\frac{R
\int_{\Omega_{-R}}\Big(|\nabla
 \widehat u_\e(y)|^2-\lambda^\e_{\bar k}\e^2 p(\e y) \widehat u_\e^2(y)\Big)dy}
{\int_{\Gamma_{R}^-}\widehat u_\e^2(y)\,d\sigma}
\end{equation}
with $\Omega_{-R}$ and $\Gamma_R^-$ as in \eqref{eq:defOmega_r} and
\eqref{eq:defGamma_r-} respectively.  From strong convergence of
$\widehat u_{\e_{n_k}}$ to $\hat c\,\widehat\Phi$ in $\mathcal H_r$
for every $r>1$, passing to the limit in (\ref{eq:117left}) along the
subsequence $\{\e_{n_k}\}_k$ and using (\ref{eq:99left}) we obtain that
\begin{equation*}
\lim_{k\to+\infty}\mathcal N_{\e_{n_k}}(-R\e_{n_k})=
\frac{R
\int_{\Omega_{-R}}|\nabla
(\hat c\,\widehat\Phi)(y)|^2dy}
{\int_{\Gamma_{R}^-}(\hat c\,\widehat\Phi)^2(y)\,d\sigma}=
\frac{R
\int_{\Omega_{-R}}|\nabla
\widehat\Phi(y)|^2dy}
{\int_{\Gamma_{R}^-}\widehat\Phi^2(y)\,d\sigma}=
\widehat{\mathcal N}(-R).
\end{equation*}
Since the limit 
depends neither on the sequence
$\{\e_n\}_{n\in\N}$ nor on its subsequence
$\{\e_{n_k}\}_{k\in\N}$, we conclude that  the convergence 
 actually holds as $\e\to 0^+$
thus  proving the lemma.
\end{pf}

\begin{Lemma}\label{l:Neps_bounded_left}
  For every $\delta>0$ there exist $K_\delta>1$, $k_\delta\in(0,1)$,
  and $\rho_\delta\in\big(0,\frac{k_\delta}{K_\delta}\big)$, such that
\begin{equation}\label{eq:161}
\mathcal N_\e(r)\leq N-1+\delta
\quad\text{and}\quad \int_{\Omega_{r}} \Big(|\nabla
u_\e|^2-\lambda^\e_{\bar k} p u_\e^2\Big)dx \geq
\frac12\int_{\Omega_{r}} |\nabla u_\e|^2dx
\end{equation}
for all 
$r\in(-k_\delta,-K_\delta\e)$ and $\e\in(0,\rho_\delta)$.
\end{Lemma}
\begin{pf}
  Let $\delta>0$ and fix $\delta'\in(0,1)$ such that 
$$
(N-1+2\delta') e^{\delta'}=N-1+\delta.
$$
From Lemma \ref{l:limitN_hat} there
  exists some $K_\delta>1$ such that $\widehat{\mathcal
    N}(-K_\delta)<N-1+\delta'$. From Lemma \ref{l:limNepsleft}
  there exists some $\e'_\delta>0$ such that, for all $\e\in(0,\e'_\delta)$,
  $\mathcal N_\e(-K_\delta\e)<\widehat{\mathcal N}(-K_\delta)+\delta'<N-1+2\delta'$.
Letting $\check R_{\delta'}$, $\check\e_{\delta'}$ as in Lemma
\ref{l:stima N+} and Corollary
\ref{cor:stima_N-}, we have that for all $\e\in
\big(0,\min \big\{\e_\delta',\check\e_{\delta'}, \check R_{\delta'}/K_\delta
\big\}\big)$ and $r\in (-\check R_{\delta'},-K_\delta\e)$
$$
\mathcal N_{\e}(r)\leq \mathcal N_{\e}(-K_\delta
\e)e^{\delta' \check R_{\delta'}}
\leq (N-1+2\delta') e^{\delta' \check R_{\delta'}}\leq N-1+\delta
$$
and $\int_{\Omega_{r}} \big(|\nabla u_\e|^2-\lambda^\e_{\bar k} p
u_\e^2\big)dx \geq \frac12\int_{\Omega_{r}} |\nabla u_\e|^2dx$.  Then the
lemma follows choosing $k_\delta=\check R_{\delta'}$ and
$\rho_\delta=\min \big\{\e_\delta',\check\e_{\delta'}, \check
R_{\delta'}/K_\delta \big\}$.
\end{pf}

\section{Asymptotics at the left junction}\label{sec:asymptotics-at-left}

Throughout this section, we fix $\delta\in(0,1)$ so that 
$N-1+\delta<N$. Let us denote $\widetilde
K=K_\delta>1$, $\tilde h=k_\delta\in(0,1)$, and $\tilde
\rho=\rho_\delta\in\big(0,\frac{\tilde h}{\widetilde K}\big)$ with
$K_\delta,k_\delta,\rho_\delta$ as in Lemma \ref{l:Neps_bounded_left},
so that
\begin{equation}\label{eq:160}
\mathcal N_\e(r)\leq N-1+\delta< N\quad\text{for all }
r\in(-\tilde h,-\widetilde K\e)\text{ and }\e\in(0,\tilde \rho).
\end{equation}
Let us denote 
\begin{equation}
  \label{eq:151}
  U_\e(x)=\frac{u_\e(x)}{\sqrt{\int_{\Gamma^-_{\tilde h}}u_\e^2d\sigma}}
\end{equation}
with $\Gamma^-_{\tilde h}$ as in \eqref{eq:defGamma_r-}. Let us
notice that, for $\e\in(0,\e_0)$, $U_\e$ solves
\begin{align}\label{eq:24norm}
\begin{cases}
-\Delta U_\e=\lambda^\e_{\bar k} p U_\e,&\text{in }\Omega^\e,\\
U_\e=0,&\text{on }\partial \Omega^\e,
\end{cases}
\end{align}
and 
\begin{equation}
  \label{eq:152}
\int_{\Gamma^-_{\tilde h}}U_\e^2d\sigma=1.
\end{equation}

\begin{Proposition}\label{p:conUeps}
  For every sequence $\e_{n}\to 0^+$ there exist a subsequence
  $\{\e_{n_k}\}_k$ and a function $U\in C^2(D^-)\cup\big(\bigcup_{t>0}\mathcal H_t^-\big)$ such
  that
\begin{itemize}
\item[i)]
$U_{\e_{n_k}}\to U$ strongly in $\mathcal
  H_t^-$ for every $t>0$
and in $C^2(\overline{B_{t_2}^-\setminus
  B_{t_1}^-})$ for all $0<t_1<t_2$;
\item[ii)]
  $U\not\equiv 0$  in $D^-$;
\item[iii)]  $U$ solves
\begin{equation}\label{eq:eqU}
\begin{cases}
-\Delta U(x)=\lambda_{k_0}(D^+)p(x)U(x),&\text{in }D^-,\\
U=0,&\text{on }\partial D^-;
\end{cases}
\end{equation}
\item[iv)] if $\mathcal N_U:(-\infty,0)\to\R$ is defined as 
  \begin{equation}\label{eq:159}
    \mathcal N_U(r):=
    \frac{(-r) \int_{\Omega_{r}}\Big(|\nabla
      U(x)|^2-\lambda_{k_0}(D^+)p(x) U^2(x)\Big)dx}{\int_{\Gamma_{-r}^-}U^2(x)\,d\sigma},
      \end{equation}
then 
\begin{equation}
  \label{eq:153}
  \mathcal N_U(r)\leq N-1+\delta \quad\text{for all
  }r\in(-\tilde h,0).
\end{equation}
\end{itemize}
\end{Proposition}
\begin{pf}
  Letting $H_\e^-(t)$ as in \eqref{eq:157}, from \eqref{eq:59},
  \eqref{eq:155}--\eqref{eq:157}, and Lemma \ref{l:Neps_bounded_left}
  it follows that
\begin{align*}
  \frac{\frac{d}{dt}H_\e^-(t)}{H_\e^-(t)}=
  -\frac{2}{t}\mathcal N_\e(-t)\geq -\frac{2N}t
\end{align*}
for all $t\in (\widetilde K\e,\tilde h)$ and $\e\in(0,\tilde \rho)$, which after
integration yields
\begin{equation}
  \label{eq:158}
  H_\e^-(t)\leq \tilde h^{2N}H_\e^-(\tilde h)t^{-2N}\quad \text{for all }t\in (\widetilde K\e,\tilde h).
\end{equation}
From \eqref{eq:151}, \eqref{eq:155}--\eqref{eq:157}, \eqref{eq:158},
and Lemma  \ref{l:Neps_bounded_left}, we deduce that
\begin{multline*}
  \frac12 \int_{\Omega_{-t}}|\nabla U_\e(x)|^2dx \leq
  \int_{\Omega_{-t}}\Big(|\nabla U_\e(x)|^2 -\lambda^\e_{\bar k} p(x)
  U_\e^2(x)\Big)dx\\
  =\frac{\int_{\Omega_{-t}}\Big(|\nabla u_\e(x)|^2 -\lambda^\e_{\bar
      k} p(x) u_\e^2(x)\Big)dx}{\int_{\Gamma^-_{\tilde h}}u_\e^2d\sigma}
  = \frac{t^{N-2}}{\tilde h^{N-1}}\mathcal N_\e(-t)\frac{ H_\e^-(t)}{
    H_\e^-(\tilde h)}\leq N\tilde h^{N+1}t^{-N-2}
\end{multline*}
for all $t\in (\widetilde K\e,\tilde h)$  and $\e\in(0,\tilde \rho)$.
Hence for every $t>0$ 
\begin{equation}\label{eq:92meno}
\{U_\e\}_{\e\in (0,\min\{\tilde \rho,t/\widetilde K\})} \text{ is bounded in }\mathcal H^-_t.
\end{equation}
Let $\e_n\to 0^+$.  From (\ref{eq:92meno}) and a diagonal process,
there exist a subsequence $\e_{n_k}\to 0^+$ and some $U\in
\bigcup_{R>0}\mathcal H_t^-$ such that $U_{\e_{n_k}}\rightharpoonup U$
weakly in $\mathcal H_t^-$ for every $t>0$ and a.e. in $D^-$. From
compactness of the embedding $\mathcal H_t^-\hookrightarrow
L^2(\Gamma_t^-)$, passing to the limit in \eqref{eq:152} we obtain
that $\int_{\Gamma^-_{\tilde h}}U^2d\sigma=1$; in particular $U\not\equiv 0$  in $D^-$.
Passing to the weak limit in \eqref{eq:24norm}, we obtain that
$U$ is a weak solution to (\ref{eq:eqU}).  By
classical elliptic estimates, we also have that $U_{\e_{n_k}}\to U$ in $C^2(\overline{B_{t_2}^-\setminus
  B_{t_1}^-})$ for all $0<t_1<t_2$. Therefore, multiplying \eqref{eq:eqU} by $U$ and integrating over
$\Omega_{-t}$, we obtain
\begin{align}\label{eq:93meno}
\int_{\Gamma_t^-}\frac{\partial U_{\e_{n_k}}}{\partial \nu}U_{\e_{n_k}}d\sigma
\to
\int_{\Gamma_t^-}\frac{\partial U}{\partial \nu}U\,d\sigma
=-\int_{\Omega_{-t}}\Big(|\nabla U(x)|^2-\lambda_{k_0}(D^+)p(x) U^2(x)\Big)
dx,
\end{align}
being $\nu=\nu(x)=\frac{x}{|x|}$. On the other hand,  multiplication of
(\ref{eq:24norm}) by $U_{\e_{n_k}}$ and integration by parts over
$\Omega_{-t}$ yield
\begin{align}\label{eq:98meno}
  \int_{\Omega_{-t}}\Big(|\nabla U_{\e_{n_k}}(x)|^2-\lambda^\e_{\bar k} p(x) U_\e^2(x)
\Big)dx=-
  \int_{\Gamma_t^-}\frac{\partial U_{\e_{n_k}}}{\partial
    \nu}U_{\e_{n_k}}d\sigma.
\end{align}
Since weak $\mathcal H_t^-$-convergence of $U_{\e_{n_k}}$ to $U$
implies that 
\begin{align}\label{eq:99meno}
  \int_{\Omega_{-t}}p(x) U_{\e_{n_k}}^2(x)dx\to
  \int_{\Omega_{-t}}p(x) U^2(x)dx
\quad \text{as }k\to+\infty,
\end{align}
combining (\ref{eq:93meno}),
(\ref{eq:98meno}), and (\ref{eq:99meno}), we conclude that 
$\|U_{\e_{n_k}}\|_{\mathcal H_R^-}\to\| U\|_{\mathcal H_R^-}$ and
then $U_{\e_{n_k}}\to U$ strongly in $\mathcal
H_t^-$ for every $t>0$.

Finally, we notice that strong $\mathcal H_t^-$-convergence of  $U_{\e_{n_k}}$
to $U$ implies that, for every $r<0$,
\begin{equation*}
  \mathcal N_{\e_{n_k}}(r)=
  \frac{(-r) \int_{\Omega_{r}}\Big(|\nabla
    U_{\e_{n_k}}(x)|^2-\lambda^\e_{\bar k}p(x)
    U_{\e_{n_k}}^2(x)\Big)dx}{\int_{\Gamma_{-r}^-}U_{\e_{n_k}}^2(x)\,d\sigma}
  \to \mathcal N_U(r)\quad\text{as }k\to+\infty,
\end{equation*}
hence, passing to the limit in \eqref{eq:160} as $\e=\e_{n_k}\to0$, we
obtain \eqref{eq:153} and complete the proof. 
\end{pf}

\begin{Lemma}\label{l:asympU}
  Let $U$ be as in Proposition \ref{p:conUeps} and let ${\mathcal
    N}_U:(-\infty,0)\to\R$ be the frequency function associated to $U$
  defined in \eqref{eq:159}. Then
 \begin{enumerate}[\rm(i)]
 \item $\lim_{r\rightarrow 0^-} {\mathcal N}_U(r)=N-1$;
\item   for every sequence $\lambda_{n}\to 0^+$ there exist a subsequence
  $\{\lambda_{n_k}\}_k$ and some constant $c\in\R\setminus\{0\}$ such
  that 
$$
\frac{U(\lambda_{n_k} x)}{\sqrt{H_U(\lambda_{n_k})}}
\mathop{\longrightarrow}\limits_{k\to+\infty}c\,\frac{x_1}{|x|^N}
$$
strongly in $\mathcal H_t^-$ for every $t>0$ and in
$C^2(\overline{B_{t_2}^-\setminus B_{t_1}^-})$ for all $0<t_1<t_2$,
where
\begin{equation}
  \label{eq:HU}
H_U(\lambda):=\frac1{\lambda^{N-1}}\int_{\Gamma_\lambda^-}U^2(x)\,d\sigma.  
\end{equation}
\end{enumerate}
\end{Lemma}
\begin{pf}
  We first notice that, letting $\e_{n}\to 0^+$ and $\{\e_{n_k}\}_k$
  as in Proposition \ref{p:conUeps}, passing to the limit as
  $k\to+\infty$, from \eqref{eq:161} and strong $\mathcal H_t^-$-convergence of
  $U_{\e_{n_k}}$ to $U$ we obtain that 
 \begin{equation}\label{eq:168}
 \int_{\Omega_{r}} \Big(|\nabla
U|^2-\lambda_{k_0}(D^+) p U^2\Big)dx \geq
\frac12\int_{\Omega_{r}} |\nabla U|^2dx
\end{equation}
for all $r\in(-\tilde h,0)$. In particular 
\begin{equation}
  \label{eq:166}
  \mathcal N_U(r)\geq 0\quad\text{for all }r\in(-\tilde h,0).
\end{equation}
  Arguing as in the proof of Lemma
\ref{l:DH-'}, we can prove that, for all $r<0$,
  \begin{equation}\label{eq:176}
    \frac{d}{dr}\mathcal N_U(r)=\nu_1(r)+\nu_2(r),
  \end{equation}
where 
\begin{align}
\label{eq:177}\nu_1(r)&=-2r\frac{\left(\int_{\Gamma_{-r}^-}\big|\frac{\partial
        U}{\partial \nu}\big|^2
      d\sigma\right)\left(\int_{\Gamma_{-r}^-}U^2(x)\,d\sigma\right)
    -\left(\int_{\Gamma_{-r}^-}U\frac{\partial U}{\partial \nu}
      d\sigma\right)^2}
  {\left(\int_{\Gamma_{-r}^-}U^2(x)\,d\sigma\right)^2}\\[3pt]
\label{eq:178}\nu_2(r)  &= \lambda_{k_0}(D^+)\frac{\int_{\Omega_{r}}(2p(x)+x\cdot \nabla p(x))
    U^2(x)dx}{\int_{\Gamma_{-r}^-}U^2(x) d\sigma}.
\end{align}
 Schwarz's inequality implies that 
 \begin{equation}
   \label{eq:164}
\nu_1(r)\geq0 \quad\text{for all }r<0.   
 \end{equation}
Furthermore 
\begin{align*}
  \frac{|\nu_2(r)|}{\mathcal N_U(r)}\leq
  \lambda_{k_0}(D^+)\frac{\int_{\Omega_{r}}|2p(x)+x\cdot \nabla p(x)|
    U^2(x)dx}{(-r)\int_{\Omega_{r}} \Big(|\nabla
    U(x)|^2-\lambda_{k_0}(D^+) p(x) U^2(x)\Big)dx}\leq\delta,\quad
  \text{for all }r\in(-\tilde h,0),
\end{align*}
where the last inequality is obtained passing to the limit as
$\e=\e_{n_k}\to0^+$ in \eqref{eq:163}. Hence from \eqref{eq:153} we
obtain that 
\begin{equation}\label{eq:165}
|\nu_2(r)|\leq \delta \big(N-1+\delta\big),\quad \text{for all
  }r\in(-\tilde h,0).
\end{equation}
From \eqref{eq:164} and \eqref{eq:165} it follows that
$\frac{d}{dr}{\mathcal N}_U$ is the sum of a nonnegative function and
of a bounded function on $(-\tilde h,0)$.  Therefore ${\mathcal
  N}_U(r)={\mathcal N}_U(-\tilde h)+\int_{-\tilde h}^r
(\nu_1(s)+\nu_2(s))\, ds $ admits a limit as $r\rightarrow 0^+$ which
is necessarily finite in view of (\ref{eq:153}) and (\ref{eq:166}).
  More precisely, denoting as 
\begin{equation}\label{eq:172}
\gamma:=\lim_{r\to 0^-}\mathcal N_U(r),
\end{equation}
(\ref{eq:153}) and (\ref{eq:166}) ensure that 
\begin{equation}\label{eq:174}
\gamma\in [0,N-1+\delta]\subset[0, N).
\end{equation}
For all $x\in D^-$ and $\lambda>0$, let us consider
\begin{equation}
  \label{eq:188}
U^\lambda(x):=\frac{U(\lambda x)}{\sqrt{H_U(\lambda)}}  
\end{equation}
where $H_U(\lambda)$ is defined in \eqref{eq:HU}.
We notice that 
\begin{equation}\label{eq:normalizationU}
  \int_{\Gamma_1^-}U_\lambda^2d\sigma=1.
\end{equation}
Furthermore, by direct calculation (see also the proof of Lemma
\ref{l:DH-'} which is analogous), we have that 
\begin{align}\label{eq:116}
  \frac{H_U'(\lambda)}{H_U(\lambda)}=-\frac2{\lambda}\mathcal
  N_U(-\lambda)\geq -\frac2{\lambda}\big(N-1+\delta\big)
\quad\text{for all }\lambda\in(0,\tilde h),
\end{align}
which after integration yields
\begin{equation}
  \label{eq:167}
  H_U(\lambda_1)\leq H_U(\lambda_2)\bigg(\frac{\lambda_2}{\lambda_1}\bigg)^{\!\!2(N-1+\delta)}
  \quad\text{for all }0<\lambda_1<\lambda_2<\tilde h. 
\end{equation}
From \eqref{eq:168}, \eqref{eq:167}, and \eqref{eq:153}, for every
$t\in(0,1)$ and $\lambda\in(0,\tilde h/t)$, we have that
\begin{align}\label{eq:169}
  \int_{\Omega_{-t}}&|\nabla
  U^\lambda(x)|^2dx=t^{N-2}\frac{H_U(\lambda t)}{H_U(\lambda)}
  \frac{\lambda t\int_{\Omega_{-\lambda t}}|\nabla
    U(x)|^2dx}{\int_{\Gamma^-_{\lambda t}}U^2d\sigma} \leq 2t^{N-2}
  t^{-2(N-1+\delta)}
  \mathcal N_U(-\lambda t)\\
  \notag&\leq 2 t^{-N-2\delta} (N-1+\delta).
\end{align}
Hence for
every $t\in(0,1)$ there exists $\lambda_t>0$ such that
\begin{equation}\label{eq:162}
  \{U^\lambda\}_{\lambda\in (0,\lambda_t)} \text{ is bounded in }\mathcal H^-_t.
\end{equation}
Let $\lambda_n\to 0^+$.
From (\ref{eq:162}) and a diagonal
process, we deduce that  there exist a
subsequence $\lambda_{n_k}\to 0^+$  and some $\widetilde U\in
\bigcup_{t>0}\mathcal H_t^-$ such that
$U^{\lambda_{n_k}}\rightharpoonup\widetilde U$
 weakly in $\mathcal H_t^-$ for every
$t>0$ and a.e. in $D^-$. 
Since $U^\lambda$ solves 
\begin{equation}\label{eq:eqUlambda}
\begin{cases}
-\Delta U^\lambda(x)=\lambda^2\lambda_{k_0}(D^+)p(\lambda x)U^\lambda(x),&\text{in }D^-,\\
U^\lambda=0,&\text{on }\partial D^-,
\end{cases}
\end{equation}
passing to the weak limit in
(\ref{eq:eqUlambda}), we obtain that $\widetilde U$ satisfies 
\begin{equation}\label{eq:171}
\begin{cases}
-\Delta \widetilde U(x)=0,&\text{in }D^-,\\
\widetilde U=0,&\text{on }\partial D^-.
\end{cases}  
\end{equation}
By compactness of the embedding $\mathcal H_1^-\hookrightarrow
L^2(\Gamma_1^-)$, passing to the limit in \eqref{eq:normalizationU},
we have that $\int_{\Gamma_1^-} \widetilde U^2d\sigma=1$. In
particular $\widetilde U\not\equiv 0$.

From Lemma \ref{l:lemmaZERO}, for every $\alpha>0$ there exists
$k_\alpha\in\N$ and $t_\alpha>0$ such that   
for all $k>k_\alpha$ and $t\in(\e_{n_k},t_\alpha)$
$$
  \int_{\Omega_{-t}}
|p(x)| U^2_{\e_{n_k}}(x)dx\leq \alpha\int_{\Omega_{-t}}
|\nabla U_{\e_{n_k}}(x)|^2dx.
$$
Strong $\mathcal H_t^-$-convergence of
  $U_{\e_{n_k}}$ to $U$ then implies that
$$
  \int_{\Omega_{-t}}
|p(x)| U^2(x)dx\leq \alpha\int_{\Omega_{-t}}
|\nabla U(x)|^2dx,\quad\text{for all }t\in(0,t_\alpha).
$$
Hence, by the change of variable $x=\lambda y$ and \eqref{eq:169}, we obtain that, for
every $s>0$, 
$$
\lambda^2 \int_{\Omega_{-s}} |p(\lambda y)| |U^\lambda(y)|^2 dy\leq
\alpha\int_{\Omega_{-s}} |\nabla U^\lambda(y)|^2dy \leq 2\alpha
s^{-N-2\delta} (N-1+\delta) ,\text{ for all
}\lambda<\min\bigg\{\frac{t_{\alpha}}{s},\frac{\tilde h}{s}\bigg\},
$$
thus implying that, for every $s>0$,
\begin{equation}
  \label{eq:170}
\lambda^2 \int_{\Omega_{-s}} |p(\lambda y)| |U^\lambda(y)|^2 dy
 \to 0\quad\text{as }\lambda\to 0^+.
\end{equation}
By classical elliptic estimates, we also have that
$U^{\lambda_{n_k}}\to \widetilde U$ in
$C^2(\overline{B_{r_2}^-\setminus B_{r_1}^-})$ for all
$0<r_1<r_2$. Therefore, multiplying~\eqref{eq:171} by
$\widetilde U$ and integrating over $\Omega_{-t}$, we obtain
\begin{align}\label{eq:93U}
\int_{\Gamma_t^-}\frac{\partial U^{\lambda_{n_k}}}{\partial \nu}U^{\lambda_{n_k}}d\sigma
\to
\int_{\Gamma_t^-}\frac{\partial \widetilde U}{\partial \nu}\widetilde U\,d\sigma
=-\int_{\Omega_{-t}}|\nabla \widetilde U(x)|^2dx,
\end{align}
while  multiplication of
(\ref{eq:eqUlambda}) by $U^{\lambda_{n_k}}$ and integration by parts over
$\Omega_{-t}$ yield
\begin{align}\label{eq:98U}
  \int_{\Omega_{-t}}|\nabla U^{\lambda_{n_k}}|^2dx=
 - \int_{\Gamma_t^-}\frac{\partial U^{\lambda_{n_k}}}{\partial
    \nu}U^{\lambda_{n_k}}d\sigma+ \lambda_{n_k}^2
\lambda_{k_0}(D^+)  \int_{\Omega_{-t}} p(\lambda_{n_k} x) |U^{\lambda_{n_k}}(x)|^2\,dx.
\end{align}
Combining (\ref{eq:93U}),
(\ref{eq:98U}), and (\ref{eq:170}), we conclude that $\|
U^{\lambda_{n_k}}\|_{\mathcal H_t^-}\to
\| \widetilde U\|_{\mathcal H_t^-}$ and
then $U^{\lambda_{n_k}}\to\widetilde U$ strongly in $\mathcal
H_t^-$ for every $t>0$.

From \eqref{eq:159}, strong convergence
$U^{\lambda_{n_k}}\to\widetilde U$ in $\mathcal H_t^-$, and \eqref{eq:170}, we have that,
for every $t>0$,
\begin{align}\label{eq:173}
  \mathcal N_U(-t\lambda_{n_k})&=
  \frac{t\lambda_{n_k} \int_{\Omega_{-t\lambda_{n_k} }}\Big(|\nabla
    U(x)|^2-\lambda_{k_0}(D^+)p(x) U^2(x)\Big)dx}{\int_{\Gamma_{t\lambda_{n_k}}^-}U^2(x)\,d\sigma}\\
\notag  &= \frac{t\int_{\Omega_{-t}}\Big(|\nabla
    U^{\lambda_{n_k}}(x)|^2-\lambda_{n_k}^2\lambda_{k_0}(D^+)p(\lambda_{n_k}x)
    |U ^{\lambda_{n_k}} (x)|^2\Big)dx} {\int_{\Gamma_{t}^-}|U
    ^{\lambda_{n_k}} (x)|^2\,d\sigma}\\
\notag&\longrightarrow 
\frac{t\int_{\Omega_{-t}}|\nabla \widetilde U(x)|^2dx}
{\int_{\Gamma_{t}^-}\widetilde U^2(x)\,d\sigma}\quad\text{as }k\to+\infty.
\end{align} 
Combining  \eqref{eq:172} and \eqref{eq:173} we conclude that 
\begin{equation*}
  \frac{t\int_{\Omega_{-t}}|\nabla \widetilde U(x)|^2dx}
{\int_{\Gamma_{t}^-}\widetilde U^2(x)\,d\sigma}=
\gamma\quad\text{for all }t>0.
\end{equation*}
From Lemma \ref{l:limNmenophi1gen} there exists $K_0\in \N$,
$K_0\geq1$, such that 
\begin{equation}
  \label{eq:175}
\gamma=N-2+K_0  
\end{equation}
and $\widetilde U(x)=|x|^{-N+2-K_0}Y(x/|x|)$ for some
  eigenfunction $Y$ of $-\Delta_{{\mathbb S}^{N-1}}$ associated to the
  eigenvalue $K_0(N-2+K_0)$, i.e.  satisfying $-\Delta_{{\mathbb
      S}^{N-1}}Y=K_0(N-2+K_0) Y$ on ${\mathbb S}^{N-1}$. From
  \eqref{eq:174} and \eqref{eq:175} we infer that necessarily $K_0=1$,
  so that 
$$
\gamma=N-1\quad\text{and}\quad
 \widetilde U(x)=|x|^{-N+1}Y(x/|x|).
$$
From $\widetilde U=0$ on $\partial D^-$, we deduce that $Y\equiv 0$ on
$\{\theta=(\theta_1,\theta_2,\dots,\theta_N)\in{\mathbb
  S}^{N-1}:\theta_1=0\}$, hence $Y$ is an eigenfunction of
$-\Delta_{{\mathbb S}^{N-1}}$ on ${\mathbb
  S}^{N-1}_-=\{\theta=(\theta_1,\theta_2,\dots,\theta_N)\in{\mathbb
  S}^{N-1}:\theta_1<0\}$ under null Dirichlet boundary conditions
associated to the eigenvalue $N-1$. It is easy to verify that $N-1$ is
the first eigenvalue of such eigenvalue problem and hence it is simple;
furthermore an eigenfunction associated to the eigenvalue $N-1$ is
$\theta=(\theta_1,\theta_2,\dots,\theta_N)\in{\mathbb
  S}^{N-1}_-\mapsto \theta_1$. Therefore we conclude that there exists
some constant $c\in\R\setminus\{0\}$ such that $Y(\theta)=c\theta_1$
and then 
$$
\widetilde U(x)=c\,\frac{x_1}{|x|^N}.
$$
The proof is thereby completed.
\end{pf}

\begin{Lemma}\label{l:limitHU}
   Let $U$ as  in Proposition \ref{p:conUeps} and let $H_U:(0,+\infty)\to\R$ 
 be  defined in
 \eqref{eq:HU}. Then 
\begin{enumerate}[\rm(i)]
 \item $H_U(\lambda)\leq e^{2\delta(N-1+\delta)\tilde h}\tilde
   h^{2(N-1)}H_U(\tilde h)\lambda^{-2(N-1)}$ for all $\lambda\in
   (0,\tilde h)$;
\item for every $\varrho>0$ there exists $\lambda_\varrho>0$ such that 
$H_U(\lambda)\geq H_U(\lambda_\varrho)
\lambda_\varrho^{2(N-1-\varrho)}\lambda^{-2(N-1-\varrho)}$ 
for all $\lambda\in
   (0,\lambda_\varrho)$; 
\item $\lim_{\lambda\to 0^+}\lambda^{2(N-1)}H_U(\lambda)$ exists and
is finite.
\end{enumerate}
\end{Lemma}
\begin{pf}
From Lemma \ref{l:asympU} (i), \eqref{eq:176}, \eqref{eq:164},
\eqref{eq:165}, we obtain that 
\begin{align*}
  N-1-\mathcal N_U(-\lambda)=\int_{-\lambda}^0\mathcal N_U'(s)ds\geq 
\int_{-\lambda}^0\nu_2(s)\,ds\geq
-\delta(N-1+\delta)\lambda\quad\text{for all }\lambda\in(0,\tilde h)
\end{align*}
where $\nu_2$ is defined in \eqref{eq:178}, and then 
$$
\mathcal N_U(-\lambda)\leq N-1+\delta(N-1+\delta)\lambda
\quad\text{for all }\lambda\in(0,\tilde h),
$$
which, together with \eqref{eq:116}, yields
\begin{align*}
     \frac{H_U'(\lambda)}{H_U(\lambda)}=-\frac2{\lambda}\mathcal
  N_U(-\lambda)\geq -\frac{2(N-1)}{\lambda}-2\delta (N-1+\delta)
\quad\text{for all }\lambda\in(0,\tilde h).
\end{align*}
Integration of the above inequality between $\lambda$ and $\tilde h$ proves
estimate (i).

 From Lemma \ref{l:asympU} (i), for any $\rho>0$ there exists $\lambda_{\varrho}>0$ such
that ${\mathcal N}_U(r)>N-1-\rho$ for any $r\in (-\lambda_{\varrho} ,0)$ and
hence
$$
\frac{H_U'(\lambda)}{H_U(\lambda)}=-\frac2{\lambda}\mathcal
  N_U(-\lambda)<-\frac{2(N-1-\rho)}{\lambda}\quad \text{for all } 
\lambda\in (0,\lambda_{\varrho}).
$$
Integration over the interval $(\lambda,\lambda_{\varrho})$ yields (ii).

In view of (i), to prove (iii) it is sufficient to show that the limit exists. 
  From \eqref{eq:116}, Lemma \ref{l:asympU} (i), and \eqref{eq:176} it
  follows that 
  \begin{align*}
    \frac{d}{d\lambda}\Big(\lambda^{2(N-1)}H_U(\lambda)\Big)&=
    2\lambda^{2N-3}H_U(\lambda)\bigg(N-1+\frac\lambda2
    \frac{H_U'(\lambda)}{H_U(\lambda)}\bigg)=
    2\lambda^{2N-3}H_U(\lambda)(N-1-\mathcal N_U(-\lambda))\\
    &= 2\lambda^{2N-3}H_U(\lambda)\int_{-\lambda}^0\mathcal N'_U(s)ds=
    2\lambda^{2N-3}H_U(\lambda)\int_{-\lambda}^0\big(\nu_1(s)+\nu_2(s)\big)ds
  \end{align*}
where $\nu_1$ and $\nu_2$ are defined in \eqref{eq:177} and
\eqref{eq:178} respectively. By integration of the above identity we
obtain that, for all $\lambda \in (0,\tilde h)$,
\begin{align}\label{eq:60}
  \lambda^{2(N-1)}H_U(\lambda)-\tilde h^{2(N-1)}H_U(\tilde
  h)=&-2\int_{\lambda}^{\tilde
    h}s^{2N-3}H_U(s)\bigg(\int_{-s}^{0}\nu_1(t)\,dt\bigg)ds \\
  \notag&-2\int_{\lambda}^{\tilde
    h}s^{2N-3}H_U(s)\bigg(\int_{-s}^{0}\nu_2(t)\,dt\bigg)ds.
\end{align}
From \eqref{eq:164} the limit 
$$
\lim_{\lambda\to 0^+}\int_{\lambda}^{\tilde
    h}s^{2N-3}H_U(s)\bigg(\int_{-s}^{0}\nu_1(t)\,dt\bigg)ds
$$ exists.  On the other hand from (i) and \eqref{eq:165} it
follows that 
\begin{equation}
  \label{eq:180}
s^{2N-3}H_U(s)\bigg(\int_{-s}^{0}\nu_2(t)\,dt\bigg)=
O(1)\quad\text{as }s\to 0^+  
\end{equation}
thus proving in particular that $s\mapsto
s^{2N-3}H_U(s)\big(\int_{-s}^{0}\nu_2(t)\,dt\big) \in L^1(0,\tilde
h)$. We conclude that both terms at the right hand side of \eqref{eq:60}
admit a limit as $\lambda\to0^+$, the second one being finite in view
of \eqref{eq:180}, thus completing the proof of the lemma.
\end{pf}

\begin{Lemma}\label{l:limitHUpos}
   Let $U$ be as  in Proposition \ref{p:conUeps},
$Y_1$ as in \eqref{eq:Y1}, and let $H_U:(0,+\infty)\to\R$ 
 be  defined in
 \eqref{eq:HU}. Then 
\begin{align}
\tag{i}
& \int_{{\mathbb S}^{N-1}_-}U(\lambda
  \theta)Y_1(\theta)\,d\sigma(\theta)\\
\notag&\qquad=\lambda^{1-N}\bigg[
\int_{\Gamma_1^-}U(x)Y_1(\tfrac{x}{|x|})\,dx
-\tfrac{\lambda_{k_0}(D^+)}N\int_{D^-}p(x)U(x)Y_1(\tfrac{x}{|x|})
\Big(|x|\alchi_{B^-_1}(x)+\tfrac{\alchi_{\Omega_{-1}}(x)}{|x|^{N-1}}\Big) dx\bigg]\\
&\notag\quad\qquad+O(\lambda^{3-N})\quad\text{as }
\lambda\to 0^+,\\[5pt]
\tag{ii}
& \lim_{\lambda\to0^+}\lambda^{2(N-1)}H_U(\lambda)>0.
\end{align}
\end{Lemma}
\begin{pf}
Let us  define, for all $\lambda>0$,
\begin{align*}
   \mu(\lambda)=\int_{{\mathbb S}^{N-1}_-}U(\lambda
  \theta)Y_1(\theta)\,d\sigma(\theta),\quad
\varsigma(\lambda)=\lambda_{k_0}(D^+)\int_{{\mathbb S}^{N-1}_-}p(\lambda
  \theta)U(\lambda
  \theta)Y_1(\theta)\,d\sigma(\theta).
\end{align*}
From \eqref{eq:eqU} $\mu$ satisfies 
\begin{equation*}
-\mu''(\lambda)-\frac{N-1}{\lambda}\mu'(\lambda)+
\frac{N-1}{\lambda^2}\mu(\lambda)=\varsigma(\lambda),\quad\text{in }(0,+\infty).
\end{equation*}
Hence there exist $c_1,c_2\in\R$ such that 
\begin{equation}
  \label{eq:181}
\mu(\lambda)=\lambda\bigg(c_1+\frac1N\int_\lambda^1\varsigma(t)dt\bigg)+
\lambda^{1-N}\bigg(c_2-\frac1N\int_\lambda^1t^N\varsigma(t)dt\bigg)
\quad\text{for all }\lambda\in(0,+\infty).
\end{equation}
Since $p\in L^{N/2}(\R^N)$ and $U\in\mathcal H^-_1$ ensure that
$\frac{p(x)U(x)Y_1(x/|x|)}{|x|^{N-1-\alpha}}\in L^1(\Omega_{-1})$ for
all $\alpha\in[0,\frac N2)$ and, for all $\lambda>1$,
\begin{align*}
&\int_\lambda^1t^\alpha\varsigma(t)dt=-\lambda_{k_0}(D^+)\int_{B^-_\lambda\setminus
  B^-_1}\frac{p(x)U(x)Y_1(x/|x|)}{|x|^{N-1-\alpha}}\,dx,\\
&\int_1^\lambda t^\alpha|\varsigma(t)|dt\leq \lambda_{k_0}(D^+)\int_{B^-_\lambda\setminus
  B^-_1}\frac{p(x)|U(x)|Y_1(x/|x|)}{|x|^{N-1-\alpha}}\,dx,
\end{align*}
we deduce that $\int_\lambda^1t^\alpha\varsigma(t)dt$ admits a finite
limit 
and $\int_1^\lambda t^\alpha|\varsigma(t)|dt=O(1)$
as $\lambda\to+\infty$ for every $\alpha\in[0,\frac N2)$.
In particular  $\int_\lambda^1\varsigma(t)dt$ admits a finite
limit as $\lambda\to+\infty$ and 
$$
\bigg|\int_\lambda^1t^N\varsigma(t)dt\bigg|
\leq
\int_1^\lambda t^{N-1}t|\varsigma(t)|dt\leq
\lambda^{N-1}\int_1^\lambda t
|\varsigma(t)|dt=O(\lambda^{N-1})\quad\text{as }\lambda\to+\infty.
$$
Hence from \eqref{eq:181} we deduce
\begin{equation}
  \label{eq:179}
\mu(\lambda)=\lambda\bigg(c_1-\frac1N\int_1^{+\infty}\varsigma(t)dt+o(1)\bigg)+O(1)
\quad\text{as }\lambda\to+\infty.  
\end{equation}
Since $U\in\mathcal H^-_1$ yields
$\int_1^{+\infty}t^{N-1}|\mu(t)|^{2^*}dt<+\infty$, \eqref{eq:179} 
necessarily implies that $c_1=\frac1N\int_1^{+\infty}\varsigma(t)dt$.
Then \eqref{eq:181} becomes 
\begin{equation}\label{eq:182}
  \mu(\lambda)=\frac{\lambda}N\int_\lambda^{+\infty}\varsigma(t)dt+
  \lambda^{1-N}\bigg(c_2-\frac1N\int_\lambda^1t^N\varsigma(t)dt\bigg)
  \quad\text{for all }\lambda\in(0,+\infty).
\end{equation}
The above formula at $\lambda=1$ yields
\begin{equation}\label{eq:184}
c_2=\mu(1)-\frac1N\int_1^{+\infty}\varsigma(t)dt =\int_{{\mathbb
    S}^{N-1}_-}U(\theta)Y_1(\theta)\,d\sigma(\theta)
-\frac{\lambda_{k_0}(D^+)}N\int_{\Omega_{-1}}\frac{p(x)U(x)Y_1(\frac{x}{|x|})}{|x|^{N-1}}dx.
\end{equation}
Since 
$$
|\varsigma(\lambda)|\leq \lambda_{k_0}(D^+)\sup_{B_1^-}|p|\sqrt{H_U(\lambda)}
\quad\text{for all }\lambda\in (0,1),
$$
from Lemma \ref{l:limitHU} (i)  we deduce that 
$$
\varsigma(\lambda)=O(\lambda^{1-N})\quad\text{as }\lambda\to 0^+.
$$
Hence 
\begin{equation}\label{eq:183}
\frac{\lambda}N\int_\lambda^{+\infty}\varsigma(t)dt=O(\lambda^{3-N})\quad\text{as
}\lambda\to 0^+,
\end{equation}
 $t^N\varsigma(t)\in L^1(0,1)$, and $t^N\varsigma(t)=O(t)$ as $t\to0^+$,
 so that 
\begin{align}\label{eq:185}
-\frac1N\int_\lambda^1t^N\varsigma(t)dt&=
-\frac1N\int_0^1t^N\varsigma(t)dt +\frac1N\int_0^\lambda t^N\varsigma(t)dt\\
\notag&=
-\frac{\lambda_{k_0}(D^+)}N\int_{B^-_{1}}|x|p(x)U(x)Y_1(\tfrac{x}{|x|})dx
+O(\lambda^2) \text{ as
}\lambda\to 0^+.
\end{align}
Combining (\ref{eq:182}--\ref{eq:185}) we obtain statement (i).

To prove (ii), let us assume by contradiction that 
$\lim_{\lambda\to0^+}\lambda^{2(N-1)}H_U(\lambda)=0$. Since, by 
 Schwarz's inequality, $H_U(\lambda)=\int_{{\mathbb S}^{N-1}_-}U^2(\lambda
  \theta)\,d\sigma(\theta)\geq |\mu(\lambda)|^2$, it would follow
  that 
  \begin{equation*}
\lim_{\lambda\to0^+}\lambda^{N-1}\mu(\lambda)=0.    
  \end{equation*}
Hence (i) would imply that 
\begin{equation*}
  \int_{\Gamma_1^-}U(x)Y_1(\tfrac{x}{|x|})\,dx
-\tfrac{\lambda_{k_0}(D^+)}N\int_{D^-}p(x)U(x)Y_1(\tfrac{x}{|x|})
\Big(|x|\alchi_{B^-_1}(x)+\tfrac{\alchi_{\Omega_{-1}}(x)}{|x|^{N-1}}\Big) dx=0
\end{equation*}
and 
\begin{equation*}
  \int_{{\mathbb S}^{N-1}_-}U(\lambda
  \theta)Y_1(\theta)\,d\sigma(\theta)=O(\lambda^{3-N})\quad\text{as }
\lambda\to 0^+.
\end{equation*}
Therefore, letting $U^\lambda$ as in \eqref{eq:188} and using Lemma
\ref{l:limitHU} (ii) with $\varrho<2$, we obtain that 
\begin{equation}
  \label{eq:190}
\int_{{\mathbb S}_-^{N-1}}U^\lambda(\theta)
Y_1(\theta)\,d\sigma(\theta)=O(\lambda^{2-\varrho}) \quad\text{as }
\lambda\to 0^+.   
\end{equation}
From Lemma \ref{l:asympU} (ii), for every sequence $\lambda_{n}\to
0^+$ 
there exist a subsequence
  $\{\lambda_{n_k}\}_k$ and some constant $c\in\R\setminus\{0\}$ such
  that 
  \begin{equation}
    \label{eq:189}
    U^{\lambda_{n_k}}\to c\,Y_1\text{ in } L^2({\mathbb S}_-^{N-1}) .
  \end{equation}
From \eqref{eq:190} and \eqref{eq:189} we infer that 
$$
0=\lim_{k\to+\infty}\int_{{\mathbb S}_-^{N-1}}U^{\lambda_{n_k}}(\theta)
Y_1(\theta)\,d\sigma(\theta)=c \int_{{\mathbb S}_-^{N-1}}
Y_1^2(\theta)\,d\sigma(\theta)=c
$$
thus reaching a contradiction and proving statement (ii).
\end{pf}

\begin{Proposition}\label{p:asymptoticU}
  Let $U$ be as in Proposition \ref{p:conUeps}. Then
$$
\lambda^{N-1} U(\lambda x)
\mathop{\longrightarrow}\limits_{\lambda \to 0^+}\beta\,\frac{x_1}{|x|^N}
$$
strongly in $\mathcal H_t^-$ for every $t>0$ and in
$C^2(\overline{B_{t_2}^-\setminus B_{t_1}^-})$ for all $0<t_1<t_2$,
where 
\begin{equation}\label{eq:beta}
\beta=-\frac{
\int_{\Gamma_1^-}U(x)Y_1(\tfrac{x}{|x|})\,dx
-\tfrac{\lambda_{k_0}(D^+)}N\int_{D^-}p(x)U(x)Y_1(\tfrac{x}{|x|})
\Big(|x|\alchi_{B^-_1}(x)+\tfrac{\alchi_{\Omega_{-1}}(x)}{|x|^{N-1}}\Big) dx}{\Upsilon_N}
\neq0
\end{equation}
and $\Upsilon_N$ is defined in \eqref{eq:upsilonN}.
\end{Proposition}
\begin{pf}
  Let $\{\lambda_n\}_{n\in\N}\subset (0,+\infty)$ such that
$\lim_{n\to+\infty}\lambda_n=0$. Then, from part (ii) of Lemma
\ref{l:asympU} and part (ii) of Lemma \ref{l:limitHUpos}, 
there exist a subsequence $\{\lambda_{n_k}\}_{k\in\N}$
and 
some constant $\beta\in\R\setminus\{0\}$ such
  that 
  \begin{equation}
    \label{eq:192}
\lambda_{n_k}^{N-1}U(\lambda_{n_k} \theta)
\mathop{\longrightarrow}\limits_{k\to+\infty}\beta\,\frac{x_1}{|x|^N}    
  \end{equation}
strongly in $\mathcal H_t^-$ for every $t>0$ and in
$C^2(\overline{B_{t_2}^-\setminus B_{t_1}^-})$ for all $0<t_1<t_2$. In particular
$$
\lambda_{n_k}^{N-1}U(\lambda_{n_k} \theta)
\mathop{\longrightarrow}\limits_{k\to+\infty}\beta\,\theta_1
\quad \text{in } 
C^{2}({\mathbb S}_-^{N-1}) \quad \text{as }k\to+\infty.
$$
From Lemma 
\ref{l:limitHUpos}
\begin{multline*}
\lim_{k\to+\infty}
\lambda_{n_k}^{N-1} \int_{{\mathbb S}^{N-1}_-}U(\lambda_{n_k}
  \theta)Y_1(\theta)\,d\sigma(\theta)\\
=\int_{\Gamma_1^-}U(x)Y_1(\tfrac{x}{|x|})\,dx
-\tfrac{\lambda_{k_0}(D^+)}N\int_{D^-}p(x)U(x)Y_1(\tfrac{x}{|x|})
\Big(|x|\alchi_{B^-_1}(x)+\tfrac{\alchi_{\Omega_{-1}}(x)}{|x|^{N-1}}\Big)
dx  ,
\end{multline*}
thus implying that 
\begin{align*}
\beta&=\frac{\int_{\Gamma_1^-}U(x)Y_1(\tfrac{x}{|x|})\,dx
-\tfrac{\lambda_{k_0}(D^+)}N\int_{D^-}p(x)U(x)Y_1(\tfrac{x}{|x|})
\Big(|x|\alchi_{B^-_1}(x)+\tfrac{\alchi_{\Omega_{-1}}(x)}{|x|^{N-1}}\Big)
dx}{\int_{{\mathbb S}^{N-1}_-}\theta_1Y_1(\theta)d\sigma(\theta)}\\
&=-\frac{\int_{\Gamma_1^-}U(x)Y_1(\tfrac{x}{|x|})\,dx
-\tfrac{\lambda_{k_0}(D^+)}N\int_{D^-}p(x)U(x)Y_1(\tfrac{x}{|x|})
\Big(|x|\alchi_{B^-_1}(x)+\tfrac{\alchi_{\Omega_{-1}}(x)}{|x|^{N-1}}\Big)
dx}{\Upsilon_N}.
\end{align*}
Hence we have proved
that $\beta$ depends neither on the sequence
$\{\lambda_n\}_{n\in\N}$ nor on its subsequence
$\{\lambda_{n_k}\}_{k\in\N}$,  thus implying that the convergence in
(\ref{eq:192}) actually holds as $\lambda\to 0^+$
and proving the proposition.
\end{pf}

The following lemmas investigate the sign of the $\beta$ in
\eqref{eq:beta}, thus allowing the study of the nodal properties of
$u_\e$ close to the left junction. 
\begin{Lemma}\label{l:nonodal}
  Let $U$ be as in Proposition \ref{p:conUeps} and $\beta\neq0$ as in
  \eqref{eq:beta}. If $\beta>0$ (respectively $\beta<0$) then there
  exists $R>0$ such that
\begin{align*}
  &\text{for every }r\in(0,R)\text{ there exists }\e_r>0\text{ such
    that }\\
  &u_\e<0\text{ (respectively $u_\e>0$)}\text{ in } \Gamma_r^-\text{
    for all }\e\in(0,\e_r).
\end{align*}
\end{Lemma}
\begin{pf}
 Let us prove the lemma under the assumption $\beta>0$ (under the
 assumption $\beta<0$ the argument is exactly the same). We claim that 
 \begin{equation}
   \label{eq:186}
   \text{there
 exists $R>0$ such that }U<0\text{ in }B_R^-.
 \end{equation}
To prove \eqref{eq:186}, let us assume by contradiction that there
exist $\lambda_n\to0^+$, $\theta_n\in {\mathbb S}^{N-1}_-$, 
$\bar\theta\in \overline{{\mathbb S}^{N-1}_-}$ such that 
$\theta_n\to\bar\theta$ and $U(\lambda_n\theta_n)\geq0$. If
$\bar\theta\in {\mathbb S}^{N-1}_-$ then from Proposition
\ref{p:asymptoticU} we obtain that 
\begin{align*}
0\leq  \lambda^{N-1}U(\lambda_n\theta_n)=\Big(
  \lambda^{N-1}U(\lambda_n\theta_n)-\beta(\theta_n)_1\Big)+\beta(\theta_n)_1
=\beta \bar\theta_1+o(1)\quad\text{as }n\to+\infty
\end{align*}
which yields a contradiction. On the other hand, if $\bar\theta\in
\partial{\mathbb S}^{N-1}_-$, i.e. if $\bar\theta_1=0$, then, letting
$s>0$ sufficiently small to have
$|x|^N-N|x|^{N-2}x_1^2>c>0$ for all $x\in
A_s:=\{x\in 
B^-_1\setminus B^-_{1/2}:x_1>-s\}$, we have that 
$(\frac{t}{\lambda_n}, \theta'_n)\in A_s$
for all $t\in
(\lambda_n(\theta_n)_1,0)$ and $n$ large. Since from Proposition
\ref{p:asymptoticU}
$\lambda^N\frac{\partial U}{\partial x_1}(\lambda x)\to
\beta\frac{|x|^N-N|x|^{N-2}x_1^2}{|x|^{2N}}$ in $C^1(\overline{A_s})$, we deduce
that $\frac{\partial U}{\partial x_1}(\lambda_n x)>0$ for all $x\in
A_s$ and $n$ large. Hence 
$$
U(\lambda_n \theta_n)=-\int_{\lambda_n(\theta_n)_1}^{0}
\frac{\partial U}{\partial x_1}(t,\lambda_n \theta'_n)\,dt<0
$$
thus giving a contradiction. Claim \eqref{eq:186} is thereby proved.
It remains to prove that
\begin{equation}
  \label{eq:187}
  \text{for every }r\in(0,R)\text{ there exists }\e_r>0\text{ such
    that }u_\e<0\text{ in }\Gamma_r^-\text{ for all }\e\in(0,\e_r).
\end{equation}
To prove \eqref{eq:187}, let us assume by contradiction that there
exist $r\in (0,R)$,  $\e_n\to0^+$, $\theta_n\in {\mathbb S}^{N-1}_-$, 
$\bar\theta\in \overline{{\mathbb S}^{N-1}_-}$ such that 
$\theta_n\to\bar\theta$ and $u_{\e_n}(r\theta_n)\geq0$ (and hence $U_\e(r\theta_n)\geq0$). If
$\bar\theta\in {\mathbb S}^{N-1}_-$ then from Proposition
\ref{p:conUeps} it follows that 
\begin{align*}
0\leq  U_{\e_n}(r\theta_n)=\Big(U_{\e_n}(r\theta_n)-U(r\theta_n)
\Big)+U(r\theta_n)
=U(r\bar\theta)+o(1)\quad\text{as }n\to+\infty
\end{align*}
which contradicts \eqref{eq:186}. On the other hand, if $\bar\theta\in
\partial{\mathbb S}^{N-1}_-$, then by Hopf's Lemma 
$\frac{\partial U}{\partial x_1}(r\bar\theta)>0$.
If $t\in
(r(\theta_n)_1,0)$,  Proposition
\ref{p:conUeps} yields
$$
\frac{\partial U_{\e_n}}{\partial x_1}(t, r\theta_n')
=\bigg(\frac{\partial U_{\e_n}}{\partial x_1}(t, r\theta_n')-
\frac{\partial U}{\partial x_1}(t, r\theta_n')\bigg)+\frac{\partial
  U}{\partial x_1}(t, r\theta_n')=\frac{\partial U}{\partial x_1}(r\bar\theta)+o(1)
$$
as $n\to+\infty$, so that 
$$
\frac{\partial U_{\e_n}}{\partial x_1}(t, r\theta_n')>0
$$
provide $n$ is sufficiently large. Therefore 
$$
U_{\e_n}(r \theta_n)=-\int_{r(\theta_n)_1}^{0}
\frac{\partial U_{\e_n}}{\partial x_1}(t,r\theta'_n)\,dt<0
$$
leads to  a contradiction proving claim \eqref{eq:187}.
\end{pf}

\noindent In fact, condition \eqref{eq:13} forces the sign of $\beta$ to be
negative, as we show below.
\begin{Lemma}\label{l:sign_beta}
  Let $U$ be as in Proposition \ref{p:conUeps} and $\beta\neq0$ as in
  \eqref{eq:beta}. Then
$$
\beta<0.
$$ 
\end{Lemma}
\begin{pf} Let us assume by contradiction that $\beta>0$. From Lemma
  \ref{l:nonodal}, for every $n$ (sufficiently large), there
exists $\e_n\in(0,1/n)$ such that 
\begin{equation}\label{eq:191}
  u_{\e_n}<0 \quad\text{on }\Gamma^-_{1/n}.  
\end{equation}
Let us denote $u_{\e_n}^-:=\max\{0,-u_{\e_{n}}\}$. From  Lemma
\ref{l:usot}, $u_{\e_n}^-=0$ on $\partial
\Omega^{\e_n}_{1+2\e_n}$. Therefore, letting  
$$
v_n:=
\begin{cases}
u_{\e_n},&\text{in }\Omega_{-1/n},\\
-u_{\e_n}^-,&\text{in }\Omega^{\e_n}_{1+2\e_n}\setminus\Omega_{-1/n},\\
0,&\text{in }\R^N\setminus \Omega^{\e_n}_{1+2\e_n},
\end{cases}
$$ 
\eqref{eq:191} ensures that $v_n\in {\mathcal
  D}^{1,2}(\Omega^{\e_n}_{1+2\e_n})\subset {\mathcal
  D}^{1,2}(\R^N)$, $v_n\not\equiv 0$ in $D^-$. 

Testing equation 
$-\Delta u_{\e_n}=\lambda^{\e_n}_{\bar k} p u_{\e_n}$ with $v_n$, we
obtain $\int_{\Omega^{\e_n}_{1+2\e_n}}|\nabla v_n|^2dx= 
\lambda^{\e_n}_{\bar k}\int_{\Omega^{\e_n}_{1+2\e_n}}pv_n^2dx$, hence,
defining 
$$
w_n:=\frac{v_n}{\sqrt{\int_{\Omega^{\e_n}_{1+2\e_n}}pv_n^2dx}},
$$
we have that $w_n\in {\mathcal
  D}^{1,2}(\Omega^{\e_n}_{1+2\e_n})\subset {\mathcal
  D}^{1,2}(\R^N)$, 
\begin{align*}
&\int_{\R^N}pw_n^2dx=\int_{\Omega^{\e_n}_{1+2\e_n}}pw_n^2dx=1,\quad
\int_{\R^N}|\nabla w_n|^2dx=\int_{\Omega^{\e_n}_{1+2\e_n}}|\nabla w_n|^2dx=
\lambda^{\e_n}_{\bar k}.
\end{align*}
Hence $\{w_{n}\}_n$ is bounded in ${\mathcal D}^{1,2}(\R^N)$ and there
exists a subsequence $\{w_{n_k}\}_k$ such that
$w_{n_k}\weakly w$ weakly in ${\mathcal D}^{1,2}(\R^N)$ and 
$w_{n_k}\to w$ a.e. in $\R^N$, 
for some $w\in {\mathcal D}^{1,2}(\R^N)$.  
Since $\mathop{\rm supp}w_n\subset
\Omega^{\e_n}_{1+2\e_n}$, a.e. convergence implies that  $\mathop{\rm
  supp}w\subset D^-$ so that $w\in {\mathcal D}^{1,2}(D^-)$.
From $\int_{\R^N}pw_n^2dx=1$ we deduce that
$\int_{D^-}pw^2dx=1$
which implies that $w\not\equiv 0$. Since $w_n$ solves 
\begin{align*}
\begin{cases}
-\Delta w_n=\lambda^{\e_n}_{\bar k} p w_n,&\text{in }\Omega_{-1/n},\\
w_n=0,&\text{on }\partial\Omega_{-1/n}\cap\partial D^-,
\end{cases}
\end{align*}
weak convergence and  (\ref{eq:52}) imply that $w$ weakly solves 
\begin{align*}
\begin{cases}
-\Delta w=\lambda_{k_0}(D^+) p w,&\text{in }D^-,\\
w=0,&\text{on }\partial D^-,
\end{cases}
\end{align*}
thus implying $\lambda_{k_0}(D^+)\in \sigma_p(D^-)$ and contradicting
assumption (\ref{eq:54}).
\end{pf}

\noindent The proofs of  the main results of the paper follow
by combining the previous results.

\begin{pfn}{Theorem \ref{t:main}}
It follows by combining Propositions \ref{p:conUeps},
\ref{p:asymptoticU} and Lemma \ref{l:sign_beta}.
\end{pfn}

\begin{pfn}{Corollary \ref{c:cor}}
It follows from Lemmas \ref{l:nonodal} and \ref{l:sign_beta}.
\end{pfn}

\end{document}